\documentclass[11pt]{amsart}   
\usepackage[margin=.5in]{geometry}

\usepackage{amssymb,amsfonts}   
\usepackage{graphicx}	        
\usepackage{hyperref}           
\usepackage{comment}
\usepackage{amsmath}
\usepackage{pst-node}
\usepackage{tikz-cd} 
\usepackage{bm}

\newtheorem{theorem}{Theorem}[section]

\newtheorem{lemma}[theorem]{Lemma}
\newtheorem{definition}[theorem]{Definition}

\newtheorem{remark}[theorem]{Remark}

\allowdisplaybreaks[1]

\begin{document}
\title[Meromorphic Forms]{ON THE ALGEBRO-GEOMETRIC ANALYSIS OF MEROMORPHIC $(1,0)$-FORMS} 
 
\author{Sergio Charles}
\date{July 10, 2017}              

\begin{abstract}
In this paper, we analyze the theory of meromorphic $(1,0)$-forms $\omega\in\mathcal{M}\Omega^{(1,0)}(\mathbb{CP}^1).$ Hence, we show that on a compact Riemann surface of genus $g=0,$ isomorphic to $\mathbb{CP}^1,$ every non-constant meromorphic function $f:X\to\mathbb{CP}^1$ has as many zeros as poles, where each is counted according to multiplicities. Such an analysis gives rise to the following result. Invoking the Riemann-Roch theorem for a compact Riemann $X$ with canonical divisor $K,$ it follows that $deg(f)=0$ for any principal divisor $(f):=D$ on $X.$ More precisely, $\ell(D)-\ell(K-D)=deg(D)+1=1$ or $\ell(D)-\ell(K-D)-1=0.$ Furthermore, for a diffeomorphism $\eta:X\to\mathbb{CP}^1$ of a certain kind, a multistep program is implemented to show $X$ is a compact algebraic variety of dimension one, i.e. a non-singular projective variety. Hence, we adopt a group-theoretic approach and provide a useful heuristic, that is, a set of technical conditions to facilitate the algebro-geometric analysis of simply connected Riemann surfaces $X.$
\end{abstract}
\maketitle
\section{Introduction\label{sec:intro}}
We introduce the theory of meromorphic $(1,0)$-forms $\omega\in\mathcal{M}\Omega^{(1,0)}$ on a class of Riemann surfaces $X.$ In the first part of the paper, we analyze the diffeomorphism $\eta:X\to\mathbb{CP}^1$ which alludes to topological considerations of complex projective space ${\mathbb{CP}}^n.$ We thereby resolve that if $\eta^{*}\omega$ is a $(1,0)$-form belonging to the sheaf of differential $(1,0)$-forms $\Omega^{(1,0)}(X)$ on the compact Riemann surface $X,$ endowed with a complex Hermitian inner product, for  $\eta^{*}\omega$ the pullback of $\omega\in\mathcal{M}\Omega^{(1,0)}(\mathbb{CP}^1)$ modulo diffeomorphism $\eta:X\to\mathbb{CP}^1,$ then the Fourier-like transform $\sigma_j(\chi)=\int_{\partial X_j}e^{-i\pi\langle \xi, \chi \rangle}{\eta}^*\omega$ cannot be of compact support for $\xi,\chi\in X_j\subset X$ unless $\omega$ is identically zero. Consequently, $\int_{\partial U_j}\omega=0$ for $\omega$ a meromorphic $(1,0)$-form, $U_j\subset\mathbb{CP}^1,$ $\sup_{\phi:=(z_1,z_2)}|g|<+\infty,$ and a globally defined chart $\phi:X\to\mathbb{C}.$  Such a result is proved via measure theory and cobordism theory. It is then shown that $\eta:X\to\mathbb{CP}^1$ is biholomorphic, where $X$ is a compact Riemann surface realized as the quotient of its universal covering by a subgroup of deck transformations, i.e. $X:=\tilde X/\Gamma.$ Thus, $X$ is then isomorphic to the complex projective line $\mathbb{CP}^1$ such that it inherits an elliptic geometry. 

Thereafter, we consider a Lie group interpretation of $\int_{\partial X_j}e^{-i\pi\langle \xi, \chi \rangle}{\eta}^*\omega\ne 0$ for which a representation of the M\"{o}bius group $Aut(\hat{\mathbb{C}})$ is considered. An explicit construction is given by considering M\"{o}bius transformations and actions of $SO(3)$ on $\mathbb{S}^2.$ Using a linear approximation, we readily obtain an expression for $\sigma_j(\chi)$ in terms of group actions on $\mathbb{S}^2.$ To prove the isomorphism  $X\cong\mathbb{CP}^1$ we develop a three step program. In particular, we show that the map $\eta$ is, in fact, a biholomorphism. Secondly, it is shown that the Riemann surface $X$ has vanishing first singular homology group $H_1(X;\mathbb{C}).$ Lastly, it is shown that the Riemann surface $X$ given by the vanishing polynomial equation
\begin{equation}
z_0^{deg(p)}\left[\eta\left(\frac{z_1}{z_0},\frac{z_2}{z_0}\right)-\hat\Phi\left(\frac{z_1}{z_0}\right)\Theta\left(\frac{z_2}{z_0}\right)\right]=0
\end{equation} is an algebraic variety of dimension one such that $X$ is necessarily compact by Griffiths and Harris \cite[Pg. 215]{Griffiths}. More precisely, this means that $X$ is not a Stein manifold.  

An analysis is carried out for the case in which the local coordinates $(z_i)$ on the manifold are isothermal, such that the metric is conformally equivalent to the (constant curvature) Euclidean metric, belonging to the equivalence class $[g]=\{g|h=\lambda^2 g \text{ for }\lambda\text{ a real-valued smooth function}\}.$ In particular, it is demonstrated that $(z_i)$ are necessarily local isothermal coordinates on the manifold if and only if \begin{equation}
\left(\frac{d\phi}{dz_1}\right)^2+\phi^2(z_1)\left[\left(\frac{d\alpha}{dz_1}\right)^2-\left(\frac{d\beta}{dz_2}\right)^2+2i\frac{d\beta}{dz_2}\frac{d\alpha}{dz_1}\right]=0,
\end{equation} obtained without invoking the Beurling transform to solve the Beltrami equation. Alternatively, if the sufficient condition $g^2=-{\bar g}^2$ is satisfied for $g\in\mathcal{M}\Omega^{(0,0)}(\mathbb{CP}^1),$ then the Beltrami equation is satisfied and $(z_1,z_2)$ must necessarily be local isothermal coordinates on $X,$ such that \begin{equation}
\chi(X)=-\frac{1}{4\pi}\iint\limits_X\frac{\Delta\rho}{e^{\rho}}dS=-\frac{1}{4\pi}\iint\limits_X\frac{\Delta log\Gamma}{\Gamma}dS=2
\end{equation} implies that the Riemann surface has genus zero. Likewise, we analyze the Beltrami equation for the complex local coordinates $w=z_1+iz_2$ and $z=x_1+ix_2$ imposed on $X.$ We also show that for $w$ assumed to have continuous partial derivatives, then $w$ is a $\mu$-quasiconformal mapping provided it satisfies the Beltrami equation
\begin{equation}
\frac{\partial w}{\partial \bar z}=\mu(z)\frac{\partial w}{\partial z}
\end{equation} for a complex-valued Lebesgue measurable $\mu$ satisfying the norm condition $|\mu|^2<1$ or $sup\text{ }|\mu|<1.$ 

By invoking Hodge theory, it is shown that the condition on genera can be strengthened. Using this we obtain a statement on homology, whereby the singular homology groups assume the form \begin{equation}
H_k(X;\mathbb{C})=\begin{cases} \mathbb{C} & \text{ if } k=0,2,\\ 0 & \text{ otherwise}  \end{cases}
\end{equation} if and only if the condition for $g=0,$ i.e. 
\begin{equation*} \int_X\frac{\Sigma\Delta\rho+\Gamma(z_1)\kappa}{\Gamma(z_1)\Sigma\Phi(z_1,z_2)}dz_1\wedge dz_2=0, \end{equation*}
in isothermal coordinates $(z_1,z_2)$ is satisfied. Finally, the paper concludes with the analysis of cohomology theory where we realize $\sigma_j(\chi)=\int_{\partial X_j}e^{-i\pi\langle\xi,\chi\rangle}\eta^*\omega$ as a de Rham homomorphism $I:H_{dR}^p(X)\to H^p(X;\mathbb{C})$ to obtain an algebro-geometric result: If $X$ is a Riemann surface belonging to the category $C,$ and if $\tilde X/\Gamma=X\cong\mathbb{CP}^1$ for $\Gamma=Aut_C(X)\equiv\pi_1(X)$ trivial where $X$ is simply connected with genus $g=0,$ then the first singular cohomology group is nontrivial if and only if the first de Rham cohomology group is nontrivial. However, since $X$ is simply connected, the first singular cohomology group $H^1(X;\mathbb{C})$ vanishes, which implies that for the induced homomorphism $\left(f_1\right)^*:H_{dR}^1(X)\to H^1(X;\mathbb{C})$
\begin{equation*}
ker\left((f_1)^*\right)=\{[\omega]\in H_{dR}^1(X):\text{ } (f_1)^*([\omega])=e_{H^1(X;\mathbb{C})}=0\}=H_{dR}^1(X)
\end{equation*} since $H^1(X;\mathbb{C})$ is the trivial group, i.e. $ker\left(\int_{c^1}\theta^1\right)=H_{dR}^1(X)$ where $c^1$ denotes a $1$-cycle in $[c^1].$
By considering the logarithmic $(1,0)$-form $\omega\in\Omega^{(1,0)}_X(logD)$ for $D$ a principal divisor on $X,$ we finally prove a special case of a well-known theorem: On any compact Riemann surface $X$ every non-constant meromorphic function $f:X\to\mathbb{CP}^1$ has as many zeros as poles, where each is counted according to multiplicities.

\section{Topology of Complex Projective Space\label{sec:topology}}
We interpret ${\mathbb{CP}}^1$ as being diffeomorphic to $\mathbb{S}^2$ to characterize the compact Riemann surface $X,$ thus facilitating later Lie group analysis. Let $\eta:X\to\mathbb{CP}^1$ be a map from the compact Riemann surface $X$ (algebraic variety of dimension one) onto the complex projective line. Similarly let $U_j\subset\mathbb{CP}^1,X_j\subset X$ be compact subsets of $\mathbb{CP}^1$ and $X,$ respectively, such that $\bigcup_jU_j\subsetneq\mathbb{C}$ and $\bigcup_jX_j\subsetneq X$ are compact. Consider a meromorphic $(1,0)$-form $\omega:=g(z)dz\in\mathcal{M}\Omega^{(1,0)}(\mathbb{CP}^1)$ defined globally with $g$ a meromorphic function on $U_j,$ that is, holomorphic on $U_j\setminus D_j$ where $g$ has singular points on $D_j.$ Then the pullback ${\eta}^{*}\omega,$ modulo projection, belongs to $\mathcal{M}\Omega^{(1,0)}(X).$ Note that the space $\Omega^{(1,0)}$ is stable under holomorphic coordinate transformations such that its elements transform tensorially. Thus, in general, the spaces $\Omega^{(1,0)}$ and $\Omega^{(0,1)}$ determine complex vector bundles on an arbitrary complex manifold.

We now give a preliminary introduction to the theory of logarithmic differential forms, originally due to Deligne \cite[Pg. 89-101]{Peters}. 
\begin{definition}
Let $X$ be a complex manifold, $D\in X$ a principal divisor, and $\omega$ a holomorphic $p$-form on $X-D.$ If $\omega$ and $d\omega$ have a pole of order at most one on $D,$ then $\omega$ is said to have a logarithmic pole along $D,$ whereby $\omega$ is formally defined as a logarithmic $p$-form. The sheaf of logarithmic $p$-forms on the manifold $X$ make up a subsheaf of the meromorphic $p$-forms on $X$ with a pole along the principal divisor $D,$ denoted by $\Omega_X^{p}(logD).$
\end{definition}
For the present discussion of Riemann surfaces, logarithmic $1$-forms have local expressions given by $\omega=\frac{df}{f}=\left(\frac{m}{z}+\frac{g'(z)}{g(z)}\right)dz$ for a meromorphic function $f(z)=z^m g(z)$ of order $m$ at $0,$ where $g$ is a non-vanishing holomorphic function at $0.$ The order $m$ of $f$ at $0$ will henceforth be $m=0,$ such that logarithmic 1-forms $\omega\in\Omega_X^{(1,0)}(logD)\subset\mathcal{M}\Omega^{(1,0)}(X)$ have local expressions $\omega=\frac{df}{f}=\frac{g'(z)}{g(z)}dz$ for $g$ a non-vanishing holomorphic function at $0.$ By definition of $\Omega_X^p(logD),$ where the exterior derivative satisfies the boundary condition $d^2=0,$ it follows that $d\Omega_X^p(logD)\subset\Omega_X^{(p+1)}(logD).$ For clarification, throughout the paper we refer to meromorphic $(1,0)$-forms on $X$; however, as shown in the concluding example, the theory can be specialized to $\Omega_X^{(1,0)}(log D)\subset\mathcal{M}\Omega^{(1,0)}(X).$
\begin{definition}
The construction $d\Omega_X^p(logD)\subset\Omega_X^{(p+1)}(logD)$ leads to the complex of sheaves $(\Omega_X^{\bullet}(logD),d^{\bullet})$ defined to be the holomorphic log complex with the corresponding divisor $D.$
\end{definition}
The complex $(\Omega_X^{\bullet}(logD),d^{\bullet})$ is a subcomplex of $j_{*}\Omega_{X-D}^{\bullet},$ for which $j:X-D\hookrightarrow X$ is the inclusion of $\Omega_{X-D}^{\bullet},$ the sheaf of holomorphic forms on $X-D.$ Let the divisor $D$ have simple normal crossings, such that $D=\sum_{\nu}D_{\nu}$ for $D_\nu$ smooth, irreducible, mutually transverse components. It follows that the divisor is locally given by the union of hyperplanes $z_1...z_n=0$ in local holomorphic coordinates. It can be shown that the stalk at $p$ of $\Omega_X^1(logD)$ satisfies \cite[Pg. 90]{Peters}, 
\begin{equation*}
\Omega_X^1(logD)_p=\mathcal{O}_{X,p}\frac{dz_1}{z_1}\otimes ...\otimes\mathcal{O}_{X,p}\frac{dz_k}{z_k}\otimes\mathcal{O}_{X,p}dz_{k+1}\otimes ...\otimes \mathcal{O}_{X,p}dz_n
\end{equation*}
for $\mathcal{O}_X$ the sheaf of structure rings on $X,$ with $\Omega_X^k(logD)_p=\bigwedge_{j=1}^k\Omega_X^1(logD)_p.$ 

The Riemann surface $X$ admits the complex inner product $(\alpha,\beta)=\int_X \alpha \wedge\star \overline{\beta}$ for $\alpha$ and $\beta$ in $\mathcal{M}\Omega^{(1,0)}_C(X)$, the sheaf of meromorphic $(1,0)$-forms of compact support on $X.$ Consider the following integral, $\tilde\sigma_j:=\int_{\partial U_j}\omega=\int_{\partial U_j} gdz\equiv\int_{U_j}d\omega,$ by Stokes' theorem, wherein $U_j$ is a compact subset of the complex plane with boundary consisting of piecewise smooth rectifiable Jordan curves, i.e., it has Lebesgue measure $\mu(U_j)<+\infty$ with $\sup_{\mathcal{P}}\sum_{k=1}^n|z_{j,k}-z_{j,k-1}|< +\infty$ for $\omega$ a meromorphic $(1,0)$-form defined in a neighborhood of the closure of $U_j$ and $\partial U_j$ given parametrically by $z_j(t)$ on the interval $a\le t\le b$ such that $\mathcal{P}=\{t_0,...,t_n\}$ is a partition of the interval $[a,b].$

Before parameterizing the Riemann surface, we begin with elementary topological considerations. In particular, we endow the topological space $\mathbb{CP}^1$ with a chart, and thereby a globally defined atlas. The following construction is due to Forster \cite[Pg. 3-4]{Forster}. Let ${\mathbb{CP}}^1:=\widehat{\mathbb{C}}$ be the one point compactification of $\mathbb{C},$ where the singleton $\{\infty\}$ is not contained in $\mathbb{C}.$ Thus, one introduces the following topology on ${\mathbb{CP}}^1.$ The open subsets of the space are the conventional sets $U\subset\mathbb{C}$ together with the modified sets $V\cup\{\infty\}$ where $V\subset\mathbb{C}$ is the complement of $K\subset\mathbb{C}$ for $K$ compact. Hence, with this topology, ${\mathbb{CP}}^1$ becomes a complex manifold homeomorphic to $\mathbb{S}^2.$ Let $U_1:={\mathbb{CP}}^1\setminus\{\infty\}=\mathbb{C},$ $U_2:={\mathbb{CP}}^1\setminus\{0\}={\mathbb{C}}^{*} \cup \{\infty\}$ and define the maps $\phi_i:U_i\to\mathbb{C},$ for $i=1,2,$ whereby $\phi_1$ is the identity and $\phi_2$ is given by 
\begin{equation*}
\phi_2(z)=\begin{cases} 1/z & \text{for } z\in{\mathbb{C}}^*, \\ 0 & \text{for } z=\infty. \end{cases}
\end{equation*} Such maps are homeomorphisms and therefore ${\mathbb{CP}}^1$ is a real two-dimensional manifold. The coverings $U_1$ and $U_2$ are connected and have non-empty intersection, meaning that ${\mathbb{CP}}^1$ is connected. The complex structure on ${\mathbb{CP}}^1$ must consequently be defined by the atlas consisting of the charts $\phi_i:U_i\to\mathbb{C},$ for $i=1,2.$ Lastly, it must be shown that the two charts are holomorphically compatible. In particular, $\phi_1(U_1\cap U_2)=\phi_2(U_1\cap U_2)={\mathbb{C}}^*$ and $\phi_2\circ\phi_1^{-1}:{\mathbb{C}}^*\to {\mathbb{C}}^*, z \mapsto 1/z$ is biholomorphic. More generally, we consider complex projective space ${\mathbb{CP}}^n,$ which is the set of all lines through the origin of ${\mathbb{C}}^{n+1}.$ Equivalently, it is defined as 
\begin{equation*}
{\mathbb{CP}}^n:=({\mathbb{C}}^{n+1}\setminus \{0\})/ {\mathbb{C}}^*\equiv({\mathbb{C}}^{n+1}\setminus \{0\})/ \sim,
\end{equation*} where ${\mathbb{C}}^*$ acts by scalar multiplication on the complex vector space ${\mathbb{C}}^{n+1}.$ If $(z_0,z_1,...,z_n)$ is a point in ${\mathbb{CP}}^n,$ then for $\lambda\in{\mathbb{C}}^*$ the two points $(\lambda z_0,\lambda z_1,...,\lambda z_n)$ and $(z_0,z_1,...,z_n)$ define the same point, inducing an equivalence class denoted by $[z_0:z_1:...:z_n]$ in homogeneous coordinates, in the sense of projective geometry. The origin $(0,0,...,0)$ does not define a point in ${\mathbb{CP}}^n$. To endow the space with a topology, let $U_i$ be the open set $U_i:=\{[z_0:...:z_n] | z_i \ne 0\}\subset {\mathbb{CP}}^n.$ Define the bijective maps \cite[Pg. 1-5]{Ionel}
\begin{equation*}
\tau_i:U_i\to{\mathbb{C}}^n, [z_0:...:z_n] \mapsto \left(\frac{z_0}{z_i},...,\frac{z_{i-1}}{z_i},\frac{z_{i+1}}{z_i},...,\frac{z_n}{z_i}\right).
\end{equation*}
Then the transition maps 
\begin{equation*}
\tau_{ij}=\tau_i \circ \tau_j^{-1}:\tau_j(U_i\cap U_j)\to \tau_i(U_i \cap U_j),
\end{equation*}
\begin{equation*}
(w_1,...,w_n)\mapsto \left(\frac{w_1}{w_i},...,\frac{w_{i-1}}{w_i},\frac{w_{i+1}}{w_i},...,\frac{w_{j-1}}{w_i},\frac{1}{w_i},\frac{w_{j+1}}{w_i},...,\frac{w_n}{w_i}\right)
\end{equation*}
are biholomorphic. That is, 
\begin{equation*}
\tau_{ij}(w_1,...,w_n)=\tau_i\circ\tau_j^{-1}(w_1,...,w_n)
\end{equation*}
\begin{equation*}
\begin{split}
&=\tau_i([w_1:...:w_{j-1}:1:w_{j+1}:...:w_n]) \\
&=\tau_i\left(\left[\frac{w_1}{w_i}:...:\frac{w_{i-1}}{w_i}:1:\frac{w_{i+1}}{w_i}:...:\frac{w_{j-1}}{w_i}:\frac{1}{w_i}:\frac{w_{j+1}}{w_i}:...:\frac{w_n}{w_i}\right]\right) \\
&=\left(\frac{w_1}{w_i},...,\frac{w_{i-1}}{w_i},\frac{w_{i+1}}{w_i},...,\frac{w_{j-1}}{w_i},\frac{1}{w_i},\frac{w_{j+1}}{w_i},...,\frac{w_n}{w_i}\right).
\end{split}
\end{equation*} Therefore, ${\mathbb{CP}}^n$ carries the structure of a complex manifold of complex dimension $n.$ 
In terms of universal coverings, ${\mathbb{CP}}^n$ can be realized as the quotient of the unit $2n+1$ sphere in ${\mathbb{C}}^{n+1}$ under the action of $U(1),$ i.e. ${\mathbb{CP}}^n=\mathbb{S}^{2n+1}/U(1).$ One obtains ${\mathbb{CP}}^n$ by first projecting onto the unit sphere, whereby every line in ${\mathbb{C}}^{n+1}$ intersects the sphere in a circle $\mathbb{S}^1,$ and then identifying the object by the natural action of $U(1).$ In particular, for $n=1,$ ${\mathbb{CP}}^1=U_0\cup U_1$ where \\ 
$U_0=\{[z_0:z_1] | z_0\ne 0\}=\left\{\left[1:\frac{z_1}{z_0}\right] \biggr\rvert z_0 \ne 0\right\}=\{[1:w] | w\in\mathbb{C}\}\cong \mathbb{S}^2\setminus \{\infty\}$ and $U_1=\{[z_0:z_1] | z_1\ne 0]\}=\left\{\left[\frac{z_0}{z_1}:1\right] \biggr\rvert z_1 \ne 0\right\}=\{[w:1] | w\in\mathbb{C}\}\cong \mathbb{S}^2\setminus \{0\}.$ It follows that $\tau_{01}=\tau_0\circ\tau_1^{-1}(w)=\tau_0([w:1])=\frac{1}{w},$ with $\tau_{10}=\tau_{01}^{-1}.$ Note that ${\mathbb{CP}}^1=\mathbb{S}^3/U(1),$ obtained by projecting to the unit sphere $\mathbb{S}^2$ and then identifying under the action of $U(1).$ Such a map induces the classical Hopf fibration $\mathbb{S}^3 \hookrightarrow \mathbb{S}^2.$

\section{Compact Support\label{sec: compact support}}
We recall the definition of an $m$-current \cite{Whitney}. 
\begin{definition} 
Let $\Omega_C^m(M)$ denote the space of smooth $m$-forms with compact support on a smooth manifold $M.$ A current is a continuous linear functional on $\Omega_C^m(M)$ in the sense of distributions. The linear functional $T:\Omega_C^m(M)\to\mathbb{R}$ is an $m$-current if it is continuous in the sense of distributions. In particular, if $\omega_n$ is a sequence of smooth forms with compact support in the same set, constructed in such a way so that all of the derivatives of their coefficients tend to $0,$ uniformly, when $n$ tends to infinity, then $T(\omega_n)$ tends to zero.
\end{definition}
The following theorem is of critical importance in the theory of differential meromorphic forms and, by extension, logarithmic forms.
\begin{theorem}\label{theorem: compact support}
If $\eta^{*}\omega$ is a $(1,0)$-form belonging to the sheaf of differential $(1,0)$-forms $\Omega^{(1,0)}(X)$ on the compact Riemann surface $X,$ for  $\eta^{*}\omega$ the pullback of $\omega\in\mathcal{M}\Omega^{(1,0)}(\mathbb{CP}^1)$ modulo diffeomorphism $\eta:X\to\mathbb{CP}^1,$ then the Fourier-like transform $\sigma_j(\chi):=\int_{\partial X_j}e^{-i\pi\langle \xi, \chi \rangle}{\eta}^*\omega$ cannot be of compact support for $\xi,\chi\in X_j\subset X$ unless $\omega$ is identically zero, i.e. $supp(\sigma_j(\chi))=\overline{\{\chi\in X_j|\sigma_j(\chi)\ne 0\}}$ is not compact. 
\end{theorem} 

\proof
The proof of this theorem will be bifurcated into two parts: one in which the $m$-current is in the dual space $\left(\Omega_C^{(m,0)}(X)\right)^*_{\mathbb{R}}$ and more generally, one in which the $m$-current is in the dual space $\left(\Omega_C^{(m,0)}(X)\right)^*_{\mathbb{C}}.$ Such an argument will require both elementary measure and cobordism theory. 

(Case 1). 
Let $\omega\in\Omega_C^{(m,0)}(X),$ where the space of currents is naturally endowed with the weak*-topology, simply called weak convergence. For $X$ a complex manifold equipped with a Hermitian inner product then the $m$-current is defined as a linear functional from $\Omega_C^{(m,0)}(X)$ to the base field, with respect to $\mathbb{R},$ for compact $X=\bigcup_j X_j$ endowed with a strong topology. Such a current $T(\omega)$ is in the dual space of $\Omega_C^{(m,0)}(X).$ Hence, $f:\Omega_C^m(X_j)\to I_j$ is in the dual space of $\Omega_C^m(X_j)$ with respect to $I_j\subset\mathbb{R}$ a compact subset (interval) of $\mathbb{R}$ for $f:=g^{-1}\circ T$ with $\left(\Omega_C^m(X_j)\right)^*_{\mathbb{R}}\ni T:\Omega_C^m(X_j)\to\mathbb{R}$ and $g:I_j\hookrightarrow\mathbb{R}$ an inclusion map. Consequently, the integral $\int_{X_j}\omega$ is an $m$-current for $dim(X_j)=m;$ i.e., an $m$-current can be defined as 
\begin{equation*}
[[\partial X_j]](\omega):=\int_{\partial X_j} \omega=\int_{X_j}d\omega\equiv[[X_j]](d\omega).\end{equation*} Similarly, the space of all $m$-currents on $X$, denoted by $D_m(X),$ is a real-valued vector space by hypothesis of $[[\partial X_j]](\omega)\in\left(\Omega_C^{(m,0)}(X)\right)^*_{\mathbb{R}}$ with operations defined by $(T+S)(\omega):=T(\omega)+S(\omega), (\lambda T)(\omega):=\lambda T(\omega).$ The support of a current $T\in\mathcal{D}_m(X)$ is the complement of the largest open subset $U\subset X$ such that $T(\omega)=0$ whenever $\omega\in\Omega_C^m(U).$ The linear subspace of $\mathcal{D}_m(X)$ consisting of currents with compact support (in the above sense) is denoted by $\mathcal{E}_m(X).$ Therefore $[[\partial X_j]](\omega):=\int_{\partial X_j} \omega=\int_{X_j}d\omega$ is an $m$-current defined for homological integration. Consider the case in which $X$ is a Riemann surface, whereby homological integration defines a $1$-current $[[\partial X_j]](\omega):=\int_{\partial X_j} \omega=\int_{X_j}d\omega$ for $\omega\in\Omega_C^{(1,0)}(X)$, in the sense of $m=1$ the complex dimension of the manifold. It follows that for the map $f_j:I_j\to X_j,$ locally defined on compact $X_j\subset X$ with $I_j\subset\mathbb{R}$ a compact subset of $\mathbb{R},$ the induced homomorphism on sheaves (by naturality of the pullback) is $f_j^*:\Omega_C^1(X_j)\to\Omega_C^1(I_j).$ For $\omega\in\Omega_C^{(1,0)}(X_j)$ (the restriction of the globally defined $(1,0)$-form, $\omega,$ to $X_j\subset X$), let $\omega=\eta^*(g(z)dz)$ for $\eta^*g\in C^{\infty}_C(X_j),$ $gdz\in\Omega_C^{(1,0)}(U_j)$ (for $U_j\subsetneq\mathbb{CP}^1$ compact and $\mathbb{CP}^1:=\bigcup_jU_j$) and $\eta:X\to\mathbb{CP}^1$ so that $\theta:=e^{-i\pi\langle\xi,\chi\rangle}\omega=e^{-i\pi\langle\xi,\chi\rangle}g(\eta)d\eta\in\Omega_C^{(1,0)}(X_j)$ with $e^{-i\pi\langle\xi,\chi\rangle}\in C^{\infty}(X_j).$ Hence,
\begin{equation*}
\begin{split}
&[[\partial X_j]](\theta)=\int_{\partial X_j}\theta=\int_{\partial X_j}e^{-i\pi\langle\xi,\chi\rangle}\omega=\int_{\partial X_j}e^{-i\pi\langle\xi,\chi\rangle}g(\eta)d\eta \\ &=\int_{f_j^{-1}(\partial X_j)}f_j^*\left(e^{-i\pi\langle\xi,\chi\rangle}g(\eta)d\eta\right)=\int_{I_j^*}e^{-i\pi\langle f_j,\chi\rangle}g(f_j)df_j:=\mathcal{F}_{I_j^*}\{g(f_j)\}
\end{split}
\end{equation*} for $[[\partial X_j]](\theta)\in\mathcal{D}_{dim(X_j)}(X_j)\subset{\mathbb{R}}^n$ and a unique $f_j:I_j^*\to\partial X_j$. It follows that $[[\partial X_j]](\theta)=Re[[\partial X_j]](\theta),$ or more precisely $\mathcal{F}_{I_j^*}\{g(f_j)\}=Re\mathcal{F}_{I_j^*}\{g(f_j)\}$ meaning that $g$ is real-valued and even. Notice, if $g$ and $f_j$ commute, namely $g\circ f_j=f_j\circ g,$ then $\mathcal{F}_{I_j^*}\{g(f_j)\}=Re\mathcal{F}_{I_j^*}\{g(f_j)\}=\mathcal{F}_{I_j^*}\{Reg(f_j)\}\equiv\mathcal{F}_{I_j^*}\{(Ref_j)\circ g\}=\mathcal{F}_{I_j^*}\{g\circ (Ref_j)\}$ for $g$ and $Ref_j=t_j$ commutable functions. As such, $[[\partial X_j]](\theta)=\mathcal{F}_{I_j^*}\{g\circ t_j\}=\int_{I_j^*}e^{-i\pi\langle t_j,\chi\rangle}g(t_j)dt_j.$ Without invoking such commutative properties, if $T^{-1}:=f_j:I_j\to X_j$, $T(X_j)=I_j,$ then 
\begin{equation*}
\begin{split}
&\sigma_j(\chi)=\int_{\partial X_j}e^{-i\pi\langle\xi,\chi\rangle}g(\eta)d\eta=\int_{f_j(I_j^*)}e^{-i\pi\langle f_j,\chi\rangle}g(f_j)df_j=\int_{f_j(I_j^*)}e^{-i\pi\langle f_j,\chi\rangle}(g\circ f_j)df_j\\ &\equiv\int_{f_j(I_j^*)}e^{-i\pi\langle f_j,\chi\rangle}(g\circ f_j)\frac{df_j}{dt}dt=\int_{I_j^*}e^{-i\pi \langle t,\chi\rangle}g(t)dt\equiv \int_{K_j}e^{-i\pi t\chi}g(t)dt
\end{split}
\end{equation*} for compact $I_j^*:=K_j\subset\mathbb{R}$ uniquely defined by $f_j:I_j^*\to \partial X_j$ and $df_j$ a Radon measure, namely the pushforward $df_j=(f_j)_*:T^*I_j\to T^*X_j,I_j\to T^*X_j$ with trivial kernel. This modified Fourier transform can be specialized to $\sigma_j(z)=\int_{K_j}e^{-i\pi zt}g(t)dt,$ by applying holomorphy under the integral, for $z\in\mathbb{C}.$ If we assume the contrapositive, i.e. if the modified Fourier-like transform $\sigma_j(z)=\int_{K_j}e^{-i\pi zt}g(t)dt$ is zero on a compact subset of $K_j\subset\mathbb{R}$ (that is, it has compact support), it has an accumulation point. Therefore, it is zero everywhere on $\mathbb{C}$ by the isolated zeros theorem, i.e. $\sigma_j(z)\equiv 0.$ By the above computation, $\sigma_j(\chi)=\int_{\partial X_j}e^{-i\pi\langle\xi,\chi\rangle}g(\eta)d\eta=\int_{\partial X_j} \theta=\int_{K_j}e^{-i\pi t\chi}g(t)dt,$ the compact support of $\omega=g(\eta)d\eta$ on $X_j$ implies the compact support of $g(t)dt$ on $K_j$, and conversely. Thus if $res(\omega)|_{X_j}$ has compact support, then $res(g)|_{X_j}$ must have compact support. It follows that $\sigma_j(z)$ and, by identification, $\int_{K_j}e^{-i\pi t\chi}g(t)dt$ cannot have compact support, meaning that $\int_{\partial X_j} e^{-i\pi \langle \xi,\chi\rangle}\omega$ cannot have compact support, unless $\omega\equiv 0$. This proves the first case.

(Case 2). The generalized $m$-current is a linear functional from $\Omega_C^{(m,0)}(X)$ to $\mathbb{C}.$
To prove the general case, we recall the definition of a cobordism. 
\begin{definition}
A cobordism between manifolds $M$ and $N$ is a compact manifold $W$ whose boundary is the disjoint union of $M$ and $N,$ $\partial W=M\sqcup N.$
\end{definition}
As before, let $\sigma_j(\chi):=\int_{\partial X_j}e^{-i\pi\langle \xi,\chi\rangle}\omega$ for $\omega\in\mathcal{M}\Omega^{(1,0)}(X)$ and let $\eta:\mathbb{CP}^1\to X$ be a diffeomorphism. That is, $[[\partial X_j]](\theta)=\int_{\partial X_j}\theta=\int_{\partial X_j}e^{-i\pi\langle \xi,\chi\rangle}\omega$ is a $1$-current in $\left(\Omega_C^{(m,0)}(X)\right)^*_{\mathbb{C}}.$ Therefore, since $dim(X)=dim(\mathbb{CP}^1)$ and because $\eta$ is a diffeomorphism, for $U_j\subset\mathbb{CP}^1$ compact and $\xi\in X_j,$ $\sigma_j(\chi)=\int_{\eta^{-1}(\partial X_j)}\eta^*\left(e^{-i\pi\langle \xi,\chi\rangle}\omega\right)=\int_{\partial U_j}e^{-i\pi\langle \xi\circ\eta,\chi\rangle}\eta^*\omega.$ Let $\eta^*\omega=\Psi(z)dz\in\mathcal{M}\Omega^{(1,0)}(\mathbb{CP}^1)$, then $e^{-i\pi\langle\xi\circ\eta,\chi\rangle}=e^{-i\pi\langle z,\chi\rangle}$ for $z\in\mathbb{CP}^1.$ It follows that $\sigma_j(\chi)=\int_{\partial U_j}e^{-i\pi\langle z,\chi\rangle}\Psi(z)dz.$ We now perform a classical surgery construction of a manifold with boundary $\mathbb{R}_j:=I_j.$ Recall that if $X,Y$ are manifolds, both with boundaries, then the boundary of the product manifold is given by $\partial(X\times Y)=(\partial X\times Y)\cup(X\times\partial Y).$ Likewise, the cobordism of two identical manifolds is $Y\sqcup Y=\partial(Y\times[0,1]).$ For $U_j\subset\mathbb{CP}^1,\text{ } U_j\sqcup U_j=\partial(U_j\times[0,1])=(\partial U_j\times[0,1])\cup(U_j\times\partial[0,1])\equiv(\partial U_j\times[0,1])\cup(U_j\times\emptyset).$ Consider the lift $\widehat{\phi}:\partial U_j\to I_j\subset\mathbb{R},$ then $I_j=\widehat\phi({(\partial U_j\times[0,1])\cup(U_j\times\emptyset)}).$ In the standard topology, $int\mathbb{Q}=\emptyset$ because there exists no open set (that is, an open interval of the form $(a,b)$) inside $\mathbb{Q}.$ Similarly, $cl\mathbb{Q}=\mathbb{R}$ because every real can be realized as the limit of a sequence of rational numbers. Therefore, we have $\partial\mathbb{Q}=cl\mathbb{Q}\setminus int\mathbb{Q}=\mathbb{R},$ and $\mathbb{Q}:=(\partial U\times[0,1])\cup(U\times\emptyset).$ Let $\partial\mathbb{Q}_j=cl\mathbb{Q}_j\setminus int\mathbb{Q}_j=\mathbb{R}_j\equiv I_j.$ Furthermore let $\gamma:U_j\to\mathbb{Q}_j$ be a smooth map, then $\gamma$ induces a lift $\partial U_j=\gamma^{-1}(\partial\mathbb{Q}_j)=\gamma^{-1}(\mathbb{R}_j)=\gamma^{-1}(I_j)$ and 
\begin{equation*}
\sigma_j(\chi)=\int_{\partial U_j}e^{-i\pi\langle z,\chi\rangle}\Psi(z)dz=\lim_{n\to\infty}\int_{\gamma^{-1}(I_j)=\partial U_j}e^{-i\pi\langle \gamma^{-1}(\epsilon_n),\chi\rangle}\Psi(\gamma^{-1}(\epsilon_n))d\gamma^{-1}(\epsilon_n)
\end{equation*} for a uniformly convergent sequence of points $\{\epsilon_n\}\in\mathbb{Q}_j.$ Consider a pointed map $\gamma_{\epsilon_n}:\left(U_j\right)_{\epsilon_n}\to\left(\mathbb{Q}_j\right)_{\epsilon_n},\left(U_j\right)_{\epsilon_n}\to\left(U_j\right)_{\epsilon_n}, \left(\mathbb{Q}_j\right)_{\epsilon_n}\to\left(\mathbb{Q}_j\right)_{\epsilon_n},$ whereby the two spaces $U_j$ and ${\mathbb{Q}}_j$ can be identified at $\epsilon_n.$ Hence, for $\gamma^{-1}_{\epsilon_n}(\epsilon_n)\in\mathbb{Q}_j,$ the pullback $\left(\gamma^{-1}_{\epsilon_n}(\epsilon_n)\right)^*$ transforms $\gamma^{-1}(\epsilon_n)$ to another point $q_n$ in $\mathbb{Q}_j.$ As such, for $\gamma_{\epsilon_n}=\gamma^{-1}_{\epsilon_n}$ bijective, we obtain the following Lebesgue integral 
\begin{equation*}
\begin{split}
&\sigma_j(\chi)=\lim_{n\to\infty}\int_{\gamma_{\epsilon_n}\circ{\gamma}^{-1}_{\epsilon_n}(I_j)={\gamma}^{-1}_{\epsilon_n}\circ{\gamma_{\epsilon_n}}(I_j)\equiv I_j}\left(\gamma^{-1}_{\epsilon_n}\right)^*\left(e^{-i\pi\langle \gamma^{-1}_{\epsilon_n}(\epsilon_n),\chi\rangle}\Psi(\gamma^{-1}_{\epsilon_n}(\epsilon_n))d\gamma^{-1}_{\epsilon_n}(\epsilon_n)\right)\\&=\lim_{n\to\infty}\int_{I_j}e^{-i\pi\langle q_n,\chi\rangle}\Psi(q_n)dq_n=\int_{I_j}e^{-i\pi\langle t,\chi\rangle}\Psi(t)dt,
\end{split}
\end{equation*} for the sequence of rational numbers $\{q_n\}$ converging uniformly to $\{t\}\in\mathbb{R}.$ But $\{q_n\}$ is arbitrary, and therefore the equality holds for all $\{q_n\}$ converging uniformly, whereby there is an uncountable infinity of such points (i.e. the set of points is dense). Specialized to $\sigma_j(z)$ for $z\in\mathbb{C},$ the linear functional becomes $\sigma_j(z)=\int_{I_j}e^{-i\pi\langle t,z\rangle}\Psi(t)dt:=\int_{I_j}e^{-i\pi\langle t,z\rangle}\hat\theta.$ If we assume the contrapositive, i.e. if $\int_{I_j}e^{-i\pi\langle t,z\rangle}\hat\theta$ is zero on a compact subset of $I_j\subset\mathbb{R}$ with accumulation point, then it will be zero on all of $\mathbb{C}$ by the isolated zeros theorem. However, because the Fourier transform is injective, this implies $\Psi(t)=0$ or $\hat{\theta}=0.$ Hence, using the former equality $\sigma_j(\chi)=\int_{\partial X_j}e^{-i\pi\langle\xi,\chi\rangle}\omega=\int_{I_j}e^{-i\pi\langle t,z\rangle}\Psi(t)dt$ and invoking the converse of the above argument, if $\omega$ has compact support then $\hat{\theta}$ has compact support. Then $\sigma_j(z)$ cannot have compact support, and thus, neither can $\sigma_j(\chi),$ i.e., $supp(\sigma_j(\chi))=\overline{\{\chi\in X_j|\sigma_j(\chi)\ne 0\}}$ is not compact, unless $\omega\equiv 0.$ The proof the theorem is now complete.
\qed

The following theorem, due to Forster \cite[Pg. 108]{Forster}, 
will be invoked in the latter discourse to prove several conditions that induce an isomorphism of complex manifolds, namely Riemann surfaces, $X:=\tilde X/\Gamma=\tilde X/Aut_C(X)\cong \mathbb{CP}^1.$
\begin{theorem}\label{theorem: solution compact support}
Suppose that $g\in\mathcal{E}(\mathbb{C})$ is of compact support. Then, there is a solution $f\in\mathcal{E}(\mathbb{C})$ of the equation $\partial f/\partial \bar z=g$ having compact support if and only if 
\begin{equation}
\int\int_{\mathbb{C}}z^ng(z)dz\wedge d\bar z=0
\end{equation}  for all $n\in\mathbb{N}$ and $f(\zeta)=\int\int_{\mathbb{C}}\frac{g(z)}{z-\zeta}dz\wedge d\bar z$, where $\mathcal{E}(\mathbb{C})$ denotes the $\mathbb{C}$-algebra of functions differentiable with respect to the local coordinates $x$ and $y$ for $z=x+iy.$
\end{theorem}
\proof
By Stokes' theorem, this condition is necessary.  The converse direction follows from Serre duality.

Consider  $\omega_k=g(z)\,d\bar z\in\Gamma(\mathbb{CP}^1,\bar K\otimes \mathcal O(-k))$ as well-defined forms with values in the line bundle $\mathcal O(-k)$, by using the standard trivialization of $\mathcal O(-k)$ and the assumption that $g$ has compact support. By Serre duality, there exists a solution of
\begin{equation*}
\bar\partial f_k=\omega_k
\end{equation*}
if and only if the integral condition is satisfied for all $n$ up to order $k-2$, as the $z^kdz$, $k=0$, $\ldots$ , $k-2$ span the holomorphic section $H^0(\mathbb{CP}^1,K\otimes \mathcal O(k)).$
The section $f_k$ is unique, for $k\geq1$, and identifying $f_k$ with a function by using the aforementioned trivialization of $\mathcal O(-k)$, all $f_k$ give rise to the same function, denoted by $f$. Therefore, there exists $f\colon\mathbb{CP}^1\to\mathbb C$ which vanishes up to arbitrary order at $\infty$ and which satisfies
\begin{equation*}
\bar\partial f=g(z)\,d\bar z.
\end{equation*}
Let $U\subset \mathbb{CP}^1$, $\infty\in U$ be an open connected set such that $g(z)\,d\bar z$ is identically zero for $z\in U$. Solutions of $\bar\partial \tilde f=g(z)\,d\bar z$ on $U$ are unique up to holomorphic $1$-forms on $U$. Hence, as $f$ vanishes up to arbitrary order at $\infty$, $f$ must vanish identically on $U$, which means that $f$ has compact support. The proof of the theorem is now complete.
\qed

The following lemma is due to Dolbeault \cite[Pg. 104-105]{Forster} and will be used to prove the existence of solutions to the inhomogeneous Cauchy-Riemann differential equation $\partial f/\partial \bar z=g$ [See Appendix \ref{sec:Dolbeault's Lemma} for a complete proof].
\begin{lemma}[Dolbeault's Lemma]\label{lemma:Forster,dolbeault} 
Suppose $g\in\mathcal{E}(\mathbb{C})$ has compact support. Then there exists a function $f\in\mathcal{E}(\mathbb{C})$ such that $\frac{\partial f}{\partial \bar z}=g.$
\end{lemma} 	

\section{Lie Group Interpretation\label{sec: lie group}}	
	For a complex manifold $(M,\Sigma)$ endowed with a complex structure $\Sigma$ and a metric $g,$ the set of diffeomorphisms of $M$ that preserve the structure (the symmetry group of the structure) is given by 
\begin{equation*}
Aut(M):=\{f:M\to M, f\text{ biholomorphic}\}.
\end{equation*}
In this context we consider the $n$-by-$n$ generalized-linear group
\begin{equation*}
GL(n,\mathbb{C}):=\{A\in M(n,\mathbb{C}):det(A)\ne 0\},
\end{equation*} specialized to $n=2,$
which is an open submanifold of $M(n,\mathbb{C})$ and thus a $2n^2$ real-dimensional manifold. It is a Lie group because matrix products and inverses are smooth functions of the real and imaginary parts of the matrix entries \cite[Pg. 152]{Lee}.

Suppose $[z:w]$ are homogeneous coordinates on $\mathbb{CP}^1,$ then the diffeomorphism $\sigma_0:\mathbb{CP}^1\to \mathbb{S}^2$ is given by $\sigma_0[z:w]=\left(\frac{Re(w\bar{z})}{|w|^2+|z|^2},\frac{Im(w\bar{z})}{|w|^2+|z|^2},\frac{|w|^2-|z|^2}{|w|^2+|z|^2}\right)$ and $\sigma_0^{-1}:\mathbb{S}^2\to\mathbb{CP}^1.$ Furthermore, suppose $\left(\begin{array}{cc}\alpha&\beta\\ \gamma&\delta\end{array}\right)\in GL(2,\mathbb{C}),$ then the linear fractional transformation $f(z)=\frac{\alpha z+\beta}{\gamma z+\delta},$ which is holomorphic for $\{z\in\mathbb{C}: \gamma z+\delta\ne 0\},$ can be extended to a meromorphic function on ${\mathbb{CP}}^1.$ The automorphism $f:\mathbb{CP}^1\to\mathbb{CP}^1$ is obviously biholomorphic. Let $SO(3)$ be the group of orthogonal $3$-by-$3$ matrices having determinant one, i.e.
$SO(3):=\{A\in GL(3,\mathbb{\mathbb{R}}): det(A)=1,A^TA=I\}.$ By identifying $\mathbb{CP}^1$ with the $2$-sphere under the diffeomorphism $\sigma_0,$ then for $A\in SO(3),$ $\sigma_0^{-1}\circ A\circ\sigma_0:\mathbb{CP}^1\to\mathbb{CP}^1$ is biholomorphic \cite[Pg. 8-9]{Forster}. Every matrix $A\in SO(3)$ may be written as a product $A=\prod_{j=1}^kA_j,$ for $A_j=\left(\begin{array}{ccc}0&1&0\\0&0&1\\1&0&0\end{array}\right)$  or as $A_j=\left(\begin{array}{cc}B_j&0\\0&1\end{array}\right).$ Therefore, for $A\in SO(3),$
$
A=\prod_{j=1}^k\left(\begin{array}{cc}B_j&0\\0&1\end{array}\right).
$ Hence $\sigma_0^{-1}\circ A\circ\sigma_0\in Aut_{\rho}(\hat{\mathbb{C}}),$ the group of deck transformations on $\mathbb{CP}^1$ for the universal covering map $\rho:\widetilde{\mathbb{CP}^1}\to\mathbb{CP}^1$ the identity. In particular, $\mathbb{CP}^1$ is realized as a universal cover by the $Aut_{\rho}(\hat{\mathbb{C}})$-action of deck transformations, an orbit $\widetilde{\mathbb{CP}^1}/Aut_{\rho}(\hat{\mathbb{C}})$ where $Aut_{\rho}(\hat{\mathbb{C}})=\pi_1(\hat{\mathbb{C}})$ is trivial. Furthermore, with every invertible $2$-by-$2$ matrix $\frak h=\left(\begin{array}{cc}a&b\\c&d\end{array}\right)$ we can associate a M\"{o}bius transformation $f(z)=\frac{az+b}{cz+d}$ such that $\frak h$ is non-singular to gaurantee conformality, i.e. $det(\frak h)\ne 0.$ For $\hat{\pi}:GL(2,\mathbb{C})\to Aut(\hat{\mathbb{C}}),$ $X$ is a Lie group and $\hat{\pi}\circ\left(\begin{array}{cc}\alpha&\beta\\ \gamma&\delta\end{array}\right)\in Aut(\mathbb{\hat C})\cong PGL(2,\mathbb{C}).$

Consequently, for $\eta^*\omega\in\mathcal{M}\Omega^{(1,0)}(X)$ a logarithmic $(1,0)$-form, $\eta:X\to\mathbb{CP}^1$ a map, $\gamma_j=\partial U_j$ a path in $\mathbb{CP}^1$ (i.e. a continuous map from the unit interval $[0,1]$ into $\mathbb{CP}^1$), and $c\in X$ a point ``lying over'' $\gamma_j(0)$ (i.e. for $p:\widetilde{\mathbb{CP}^1}\to\mathbb{CP}^1$ an identity cover, $p(c)=\gamma_j(0)$), then there exists a unique path $\Gamma_j=\partial X_j$ lying over $\gamma_j$ (for $p\circ\Gamma_j=\gamma_j$) such that $\Gamma_j(0)=c.$ The curve $\Gamma_j$ is the lift of $\gamma_j$ by $p.$ Let $\xi,\chi\in X_j$ and let $\Pi$ be an action of $\frak{aut}_{\rho}(X),$ the Lie subalgebra associated with $Aut_{\rho}(X),$ such that $Lie(Aut_{\rho}(X))\cong \frak{aut}_{\rho}(X)$ for $Aut_{\rho}(\hat{\mathbb{C}}):=\{f|f(z)=\frac{az+b}{cz+d},\text{ }ad-bc\ne 0\}\cong PGL(2,\mathbb{C}).$ Assuming the action of $\Pi$ on $X$ is transitive, then there exists a $\Pi\in\frak{aut}(X)$ such that $\Pi_{\frak{aut}(X)}\xi=\chi.$ Therefore,
\begin{equation}\label{equation:lie group formulation}
\sigma_j(\chi)=\int_{\partial X_j}e^{-i\pi\langle\xi,\Pi\xi\rangle}\eta^*\omega=\int_{\partial X_j}e^{-i\pi\left\langle\xi,log\left(e^{\Pi\xi}\right)\right\rangle}\eta^*\omega=\int_{\Gamma_j}e^{-i\pi\langle\xi,\xi log\bar\Pi\rangle}\eta^*\omega
\end{equation} for $\bar\Pi=e^{\Pi}\in Aut(\hat{\mathbb{C}}),$ generated by the exponential map. It follows that $\hat{\pi}\circ\left(\begin{array}{cc}\alpha&\beta\\ \gamma&\delta\end{array}\right)\in Aut(\hat{\mathbb{C}})$ for $\hat{\pi}:GL(2,\mathbb{C})\to Aut(\hat{\mathbb{C}}).$ By letting $\bar\Pi=\hat{\pi}\circ\left(\begin{array}{cc}\alpha&\beta\\ \gamma&\delta\end{array}\right),$ then 
\begin{equation*}
\sigma_j(\chi)=\int_{\Gamma_j}e^{-i\pi\left\langle\xi,\xi log\left(\hat{\pi}\circ\left(\begin{array}{cc}\alpha&\beta\\ \gamma&\delta\end{array}\right)\right)\right\rangle}\eta^*\omega=\int_{\Gamma_j}e^{-i\pi|\xi|^2log(\hat{\pi}\circ\frak h)}\eta^*\omega=\int_{\Gamma_j}\left(\hat{\pi}\circ\frak h\right)^{-i\pi|\xi|^2}\eta^*\omega
\end{equation*} for $\frak h=\left(\begin{array}{cc}\alpha&\beta\\ \gamma&\delta\end{array}\right)\in GL(2,\mathbb{C}).$ The exponential map restricts to a diffeomorphism from some neighborhood $0$ in $\frak g:=\frak{so}(3)$ to a neighborhood of $e$ in $G:=SO(3).$ If $Aut(X)$ is the Lie group of automorphisms of $X,$ then $\frak{aut}(X)$ is its Lie algebra and $\Phi: SO(3)\to Aut(X)$ is a Lie group homomorphism such that the following diagram commutes: 

\[ \begin{tikzcd}
\frak{so}(3) \arrow{r}{\Phi_*} \arrow[swap]{d}{exp} & \frak{aut}(X) \arrow{d}{exp} \\%
SO(3) \arrow{r}{\Phi}& Aut(X)
\end{tikzcd}. 
\] Then $Aut(X)$ inherits the manifold structure of $SO(3)$ under the mapping $\sigma_0^{-1}\circ A\circ\sigma_0:\mathbb{CP}^1\to\mathbb{CP}^1$ for the diffeomorphism $\sigma_0:\mathbb{CP}^1\to\mathbb{S}^2$ and $A\in SO(3).$  
Note that $\mathbb{S}^n$ is a homogeneous space for $O(n+1),$ which induces the fiber bundle $O(n)\to O(n+1)\to \mathbb{S}^n,$ whereby the orthogonal group $O(n+1)$ acts transitively on the unit sphere, and the stabilizer of a point on $\mathbb{S}^n$ is $O(n).$ 
We now construct the group actions associated with the rotation group $SO(3)$ and obtain an explicit homomorphism between $SO(3)$ and $SU(2),$ the universal covering group of $SO(3).$ This general construction is due to Gelfand, Minlos and Shapiro (1963)\cite{Gelfand}. Let $g\in SO(3)$ be a rotation, then the action $\Sigma_s(g):\mathbb{S}^2\to \mathbb{S}^2$ on the embedding space ${\mathbb{R}}^3$ maps points on $\mathbb{S}^2$ to other points on $\mathbb{S}^2.$ By forming the composition $\sigma_0^{-1}\circ\Sigma_s(g)\circ \sigma_0$ of $\mathbb{CP}^1$ with $P'\in\mathbb{CP}^1$ and $P\in\mathbb{S}^2,$ $\sigma_0^{-1}\circ\Sigma_s(g)\circ \sigma_0:\zeta=P'\mapsto P\mapsto \Sigma_s(g)P\equiv gP\mapsto \sigma_0^{-1}(gP)=\sigma_0^{-1}(g)\sigma_0^{-1}(P)=\sigma_0^{-1}(\Sigma_s(g))\zeta:=\Sigma_u(g)\zeta=\zeta'.$ In this regard, $\Sigma_u(g)$ is a transformation (that is, an automorphism) of $\mathbb{CP}^1$ associated with the transformation $\Sigma_s(g)$ on the embedding space ${\mathbb{R}}^3.$ As such, the two rotations $g_{\phi}$ and $g_{\theta}$ through an angle of $\phi$ about the $z$-axis and $\theta$ about the $x$-axis, respectively, correspond to automorphisms of $\mathbb{CP}^1.$ In fact, these rotations generate all of $SO(3)$ where the composition of $g_{\phi}$ and $g_{\theta}$ corresponds to the composition of M\"{o}bius transformations where a M\"{o}bius transformation,
\begin{equation*}
\zeta'=\frac{\alpha\zeta+\beta}{\gamma\zeta+\delta} \text{ with } \alpha\delta-\beta\gamma\ne 0
\end{equation*} to ensure conformality, can be represented by a matrix transformation
\begin{equation*}
\left(
\begin{array}{ccc}
\alpha & \beta \\ \gamma & \delta 
\end{array}
\right), \alpha\delta-\beta\gamma=1.
\end{equation*} However, these matrices are not uniquely determined, for multiplication by $-I$ still corresponds to the same fractional linear transformation. Therefore, a given M\"{o}bius transformation corresponds to two matrix representations $g,-g\in SL(2,\mathbb{C}).$ Explicitly, the actions on $\mathbb{S}^2$ become
\begin{equation*}
\begin{split}
& \Sigma_u(g_\phi)=\Sigma_u\left[\left(
\begin{array}{ccc}
cos\phi & -sin\phi & 0 \\ sin\phi & cos\phi & 0 \\ 0 & 0 & 1
\end{array}
\right)\right]=\pm \left(
\begin{array}{ccc}
e^{i\frac{\phi}{2}} & 0 \\ 0 & e^{-i\frac{\phi}{2}}
\end{array}
\right),\\ &
\Sigma_u(g_\theta)=\Sigma_u\left[\left(
\begin{array}{ccc}
1 & 0 & 0 \\ 0 & cos\theta & -sin\theta \\ 0 & sin\theta & cos\theta
\end{array}
\right)\right]=\pm \left(
\begin{array}{ccc}
cos\frac{\theta}{2} & isin\frac{\theta}{2} \\ isin\frac{\theta}{2} & cos\frac{\theta}{2}
\end{array}
\right),
\end{split}
\end{equation*} which are in fact unitary matrices $\Sigma_u(SO(3))\subset SU(2)\subset SL(2,\mathbb{C}).$ It follows that for a general rotation $g_{\phi}g_{\theta}g_{\psi},$ 
\begin{equation*}
g(\phi,\theta,\psi)\equiv g_{\phi}g_{\theta}g_{\psi}= \left(
\begin{array}{ccc}
cos\phi & -sin\phi & 0 \\ sin\phi & cos\phi & 0 \\ 0 & 0 & 1
\end{array}
\right) \left(
\begin{array}{ccc}
1 & 0 & 0 \\ 0 & cos\theta & -sin\theta \\ 0 & sin\theta & cos\theta
\end{array}
\right)
\left(
\begin{array}{ccc}
cos\psi & -sin\psi & 0 \\ sin\psi & cos\psi & 0 \\ 0 & 0 & 1
\end{array}
\right),
\end{equation*} the group action is \cite[Ch. 3 \textsection 16]{van der Waerden}
\begin{equation*}
\begin{split}
& \Sigma_u(g(\phi,\theta,\psi))=\pm \left(
\begin{array}{ccc}
e^{i\frac{\phi}{2}} & 0 \\ 0 & e^{-i\frac{\phi}{2}}
\end{array}
\right) 
\left(\begin{array}{ccc}
cos\frac{\theta}{2} & isin\frac{\theta}{2} \\ isin\frac{\theta}{2} & cos\frac{\theta}{2}
\end{array} \right)
\left(
\begin{array}{ccc}
e^{i\frac{\psi}{2}} & 0 \\ 0 & e^{-i\frac{\psi}{2}}
\end{array}
\right) \\ &
\pm 
 \left(
\begin{array}{ccc}
cos\frac{\theta}{2}e^{i\frac{\phi+\psi}{2}} & isin\frac{\theta}{2}e^{i\frac{\phi-\psi}{2}} \\ isin\frac{\theta}{2}e^{-i\frac{\phi-\psi}{2}} & cos\frac{\theta}{2}e^{-i\frac{\phi+\psi}{2}}
\end{array}
\right).
\end{split}
\end{equation*} Thus $\rho:SU(2)\to SO(3)$ is an onto group homomorphism that completely characterizes the universal covering map of the rotation group $SO(3).$ As before, in the Lie group formulation, Eq. \ref{equation:lie group formulation}, we let $\Pi=\sigma_0^{-1}\circ\Sigma_s(g)\circ \sigma_0$ act transitively on the $SO(3)$-space $\mathbb{S}^2,$ for $g\in SO(3).$ Let $\xi$ be chosen such that $\chi=\Pi\xi\equiv \sigma_0^{-1}\circ\Sigma_s(g(d\phi,d\theta,d\psi))\circ \sigma_0\xi,$ for $g(d\phi,d\theta,d\psi)\in SO(3)$ an infinitesimal rotation. Hence, $\bar\Pi=e^{\Pi}\in G:=Aut(\hat{\mathbb{C}})$ inherits a manifold structure from $SO(3),$ and thus we can treat $\Pi$ as an element of a Lie algebra $\frak g$ such that the Lie bracket of $\frak g$ is given by the matrix commutator $[X,Y]=XY-YX$ with the Lie group generated by the exponential map $exp:\frak g\to G$ for $G=Aut(\hat{\mathbb{C}}),$
\begin{equation*}
\begin{split}
& \bar\Pi=exp\left({\sigma_0^{-1}\circ\Sigma_s(g(d\phi,d\theta,d\psi))\circ \sigma_0}\right)=\sum_{k=0}^{\infty}\frac{1}{k!}\left(\sigma_0^{-1}\circ\Sigma_s(g(d\phi,d\theta,d\psi))\circ \sigma_0\right)^k \\ & =I+\sigma_0^{-1}\circ\Sigma_s(g(d\phi,d\theta,d\psi))\circ \sigma_0+\mathcal{O}\left(\left(\sigma_0^{-1}\circ\Sigma_s(g(d\phi,d\theta,d\psi))\circ \sigma_0\right)^2\right)\\&\sim I+\sigma_0^{-1}\circ\Sigma_s(g(d\phi,d\theta,d\psi))\circ \sigma_0,
\end{split}
\end{equation*} to first order linear approximation for $\Sigma_s(g)\subset SU(2),$ such that Eq. \ref{equation:lie group formulation} becomes
\begin{equation}
\begin{split}
&\sigma_j(\chi)=\int_{\partial X_j}e^{-i\pi\langle\xi,\xi log\bar\Pi\rangle}\eta^*\omega\sim\int_{\partial X_j}e^{-i\pi\langle\xi,\xi log(I+\sigma_0^{-1}\circ\Sigma_s(g(d\phi,d\theta,d\psi))\circ \sigma_0)\rangle}\eta^*\omega\\&=\int_{\partial X_j}e^{-i\pi|\xi|^2log(I+\sigma_0^{-1}\circ\Sigma_s(g(d\phi,d\theta,d\psi))\circ \sigma_0)}\eta^*\omega=\int_{\partial X_j}(I+\sigma_0^{-1}\circ\Sigma_s(g(d\phi,d\theta,d\psi))\circ \sigma_0)^{-i\pi|\xi|^2}\eta^*\omega.
\end{split}
\end{equation} Note that the exponential mapping restricts to a diffeomorphism from a neighborhood centered at $0$ of $\frak g$  to a neighborhood centered at $e$ of $G.$ In particular, the differential $(d exp)_0:T_0\frak g\to T_eG$ is the identity map, under the canonical identifications of both $T_0\frak g$ and $T_eG$ with $\frak g$ itself.

\section{Generalized Theory\label{sec:general}}
In what follows, we elaborate further on this condition for a more general class of compact Riemann surfaces $\tilde X/\Gamma,$ identified up to a conformal equivalence class of metrics $[g],$ and prove a stronger theorem. As a topological introduction, the weak topology with respect to the base field $K=\mathbb{C}$ on $X=\tilde X/\Gamma$ is denoted $\sigma(X,\mathbb{C}).$ A subspace for the weak topology is a collection of sets of the form $\phi^{-1}(U)$ where $\phi\in X^*$ and $U$ is an open subset of the basefield $\mathbb{C}.$ Thus, a subset $V$ of $X$ is open in the weak topology if and only if it can be written as $V=\bigcup_i V_i$ each of which is an intersection of finitely many sets, $V_i=\bigcap_j \phi^{-1}(U_{i,j}),$ of the form $\phi^{-1}(U_{i,j}).$

Let $(z_i)$ be smooth local coordinates on the compact Riemann surface $(X, g)$ endowed with a complex structure, i.e. an atlas $\bigcup_j (X_j,\phi_j)$ for $V_j\subset X$ an open subset $V_j$ of $X,$ with $\phi_j:V_j\to U_j,$ $U_j\subset\mathbb{C}$ an open subset of $\mathbb{C}$ and the globally defined chart $\phi:X\to\mathbb{C}.$ We define the map $\eta:X\to\mathbb{CP}^1$ and show it to be a diffeomorphism in Section \ref{sec:biholomorphicity}. Here, $X$ is endowed with a strong topology such that it can be covered by $X=\bigcup_j \eta^{-1}(U_j)$ for $U_j\subsetneq\mathbb{CP}^1,$ henceforth, a strict compact subset of $\mathbb{CP}^1$ such that $\bigcup_jU_j=\mathbb{CP}^1$ is compact. Likewise, we let $X_j\subsetneq X$ be a strict compact subset of $X$ with $\bigcup_jX_j=X$ compact. 
If $(U,\phi=(z_1,z_2))$ is a chart of $X$ for $(z_1,z_2)$ identified with $z_1+iz_2,$ then a local expression for $g$ can be computed as follows. In the local coordinate system on $X,$ given by the real-valued functions $z_1,z_2$ (we use covariant indices only for later convenience, where this notation bears no geometric meaning in this context), for the coordinate vector fields $\{\frac{\partial}{\partial z_1},\frac{\partial}{\partial z_2}\},$ $\{dz_1,dz_2\}$ are the dual $1$-forms. For $p\in U$ and $u,v\in T_pX,$ we can write
\begin{equation*}
u=\sum_iu_i\frac{\partial}{\partial z_i}\biggr |_{p} \text{ and } v=\sum_jv_j\frac{\partial}{\partial z_j}\biggr |_p.
\end{equation*} By letting $g^{ij}(p)=g_p\left(\frac{\partial}{\partial z_i},\frac{\partial}{\partial z_j}\right)$ and using bilinearity,
\begin{equation*}
g_p(u,v)=\sum_{i,j}u_iv_ig_p\left(\frac{\partial}{\partial z_i},\frac{\partial}{\partial z_j}\right)=\sum_{i,j}g^{ij}(p)u_iv_i.
\end{equation*} It follows that the metric is simply 
\begin{equation*}
g=\sum_{i,j}g^{ij}dz_i\otimes dz_j=\sum_{i\le j}\tilde{g}^{ij}dz_idz_j\equiv Edz_1^2+2Fdz_1dz_2+Gdz_2^2,
\end{equation*} the fundamental form on $X,$ where $E=g_{11},$ $F=g_{12},$ and $G=g_{22}$ for
\begin{equation*}
\left[g_{ij}\right]=\left(\begin{array}{ccc}
E & F \\ F & G  
\end{array}
\right)
\end{equation*} the rank 2 metric tensor defined on the complex manifold. 

Let $\tilde X$ denote the universal cover of $X,$ for which $\rho:\tilde X\to X$ is a universal covering map. In particular, the local coordinate $w\in X$ is the linear combination $w=z_1+iz_2$ and $z=x_1+ix_2\in\mathbb{CP}^1.$ As before, consider the functional $\sigma_j(\chi)\in\left(\Omega^{(1,0)}(X)\right)^*$ defined by
\begin{equation*}
\sigma_j(\chi):=\int_{\partial X_j} e^{-i\pi\langle\xi,\chi\rangle}\eta^*\omega
\end{equation*} for the globally defined smooth map $\eta:X\to\mathbb{CP}^1$ (we assume this to be a diffeomorphism by hypothesis) with $\eta^*\omega\in\mathcal{M}\Omega^{(1,0)}(X)$ and $\omega\in\mathcal{M}\Omega^{(1,0)}(\mathbb{CP}^1)\supset\Omega^{(1,0)}_{\mathbb{CP}^1}(log D)$ a meromorphic $(1,0)$-form $\omega=gdz$ (where $g:\mathbb{C}\to\mathbb{C}$ is a meromorphic function and $g:\mathbb{C}\to\mathbb{CP}^1$ is analytic). Here $\partial X_j$ is a rectifiable Jordan curve, lifted onto $X$ by the inverse mapping $\eta^{-1},$ such that $\sup_{\mathcal{P}}\sum_{k=1}^n|z_{j,k}-z_{j,k-1}|< +\infty,$ for $\partial X_j$ given parametrically by $z_j(t)$ on the interval $a\le t\le b$ where $\mathcal{P}=\{t_0,...,t_n\}$ is a partition of the interval $[a,b].$ Likewise, $(TX\times TX)^*\ni\langle\xi,\chi\rangle:TX\times TX\to\mathbb{C}$ is the standard Hermitian inner product, i.e. a positive definite symmetric bilinear form, defined locally on $X$ (or rather pointwise ${\langle\xi,\chi\rangle}_p:T_pX\times T_pX\to\mathbb{C}$) and $D$ is a divisor of $\mathbb{C}$ with simple, normal crossings for which $D:=\sum_{\nu}D_{\nu}$. Here the $D_{\nu}$ are smooth, irreducible, and mutually transverse components. By Stokes' theorem, 
\begin{equation*}
\begin{split}
&\sigma_j(\chi)=\int_{\partial X_j}e^{-i\pi\langle\xi,\chi\rangle}\eta^*\omega=\int_{\partial X_j}g(\eta)e^{-i\pi\langle\xi,\chi\rangle}\frac{\partial \eta}{\partial z_1}dz_1+g(\eta)e^{-i\pi\langle\xi,\chi\rangle}\frac{\partial \eta}{\partial z_2}dz_2 \\ & =\int_{\partial X_j}M(z_1,z_2)dz_1+N(z_1,z_2)dz_2=\int_{X_j}\left(\frac{\partial N}{\partial z_1}-\frac{\partial M}{\partial z_2}\right)dz_1\wedge dz_2:=\int_{\partial X_j}h_{\bar\Pi}dz_1\wedge dz_2
\end{split}
\end{equation*} where $h_{\bar\Pi}=\frac{\partial N}{\partial z_1}-\frac{\partial M}{\partial z_2}$ for obvious notational reasons subject to the discussion of Lie representations. Furthermore, let $Fdz_1\wedge dz_2=d(z\circ \eta)\wedge d(\bar z\circ\eta)=d\eta\wedge d\bar\eta$ for $F$ defined sectionally on $X,$ which behaves locally like the total product space $\eta^{-1}\left(Z_1\times Z_2\right)\cong \eta^{-1}(V)\subset X,$ a local fiber bundle, and $(z_1,z_2)\in Z_1\times Z_2.$ Thus, $\int_{X_j}h_{\bar\Pi}dz_1\wedge dz_2=\int_{X_j}\eta^n \hat g(\eta)d(z\circ\eta)\wedge d(\bar z\circ\eta)$ if $\hat g(z)=\frac{h_{\bar\Pi}}{Fz^n}.$ Under such a substitution, $\int_{X_j}\frac{h_{\bar\Pi}}{F}d(z\circ\eta)\wedge d(\bar z\circ\eta)=\int_{X_j}h_{\bar\Pi}dz_1\wedge dz_2$ or
\begin{equation*}
\int_{X_j}\eta^n\hat g(\eta)d(z\circ\eta)\wedge d(\bar z\circ\eta)=\int\int_{\eta^{-1}(U_j)}\eta^*(z^n\hat g(z)dz\wedge d\bar z)=\int\int_{U_j}z^n\hat g(z)dz\wedge d\bar z.
\end{equation*} By Theorem \ref{theorem: solution compact support}, it follows that if $\hat g(\eta)\in\mathcal{E}(X_j)\subset\mathcal{E}(X)$ does not have compact support, there exists no $f\in\mathcal{E}(X_j)$ which satisfies $\partial f/\partial\bar \eta=\hat g$ having compact support and $\int\int_{X_j}{\eta}^n\hat g(\eta)d\eta\wedge d\bar \eta\ne 0$ for all $n\in\mathbb{N}.$ Therefore, if $\sigma_j(\chi)=\int\int_{U_j}z^n\hat g(z)dz\wedge d\bar z=\int\int_{X_j}{\eta}^n\hat g(\eta)d\eta\wedge d\bar \eta\ne 0$ then there exists no solution $f(\eta)=\int\int_{X_j}\frac{\hat g({\eta}^*)}{{\eta}^*-{\eta}}d{\eta}^*\wedge d\overline{{\eta}^*}$ having compact support. As such, if this hypothesis is true, that is if $\hat g\in\mathcal{E}(X_j)$ does not have compact support by $\hat g(\eta)=1/{\eta}^n$ under the assumption $h_{\bar\Pi}=F$ in the below, then $f\in\mathcal{E}(X_j)$ cannot have compact support. In particular, let $\hat g(z)=\frac{h_{\bar\Pi}}{z^nF(z)}$ for the section $F$ defined in terms of $z\in\mathbb{CP}^1.$ Invoking Lemma \ref{lemma:Forster,dolbeault}, 
\begin{equation*}
f(\eta)=\iint\limits_{X_j}\frac{\hat g({\eta}^*)}{{\eta}^*-\eta}d{\eta}^*\wedge d\overline{{\eta}^*}=\iint\limits_{X_j}d\omega
\end{equation*} is a solution to $\frac{\partial f}{\partial \bar \eta}=g,$ where $\omega:=\left(\int \frac{\hat g(\eta^*)}{\eta^*-\eta}d\eta^*\right)d\overline{\eta^*}$. Moreover, $d\omega=\frac{\partial}{\partial \eta^*}\left(\int\frac{\hat g(\eta^*)}{\eta^*-\eta}d\eta^*\right)d\eta^*\wedge d\overline{\eta^*}+\frac{\partial}{\partial \overline{\eta^*}}\left(\int\frac{\hat g(\eta^*)}{\eta^*-\eta}d\eta^*\right)d\overline{\eta^*}\wedge d\overline{\eta^*}\equiv\frac{\hat g(\eta^*)}{\eta^*-\eta}d\eta^*\wedge d\overline{\eta^*}.$ It follows that 
\begin{equation*}
\begin{split}
&f(\eta)=\int_{X_j}d\omega=\int_{\partial X_j}\omega=\int_{\partial X_j}\left(\int \frac{\hat g(\eta^*)}{\eta^*-\eta}d\eta^*\right)d\overline{\eta^*}\\&=\int\int_{\partial X_j}\frac{\hat g(\eta^*)}{\eta^*-\eta}d\overline{\eta^*}d\eta^*.
\end{split}
\end{equation*} Using the property $z^*\overline{z^*}=|z^*|^2$ on the $\mathbb{C}$-algebra, $d\overline{\eta^*}=\left(2-\frac{|\eta^*|^2}{{\eta^*}^2}\right)d\eta^*=\phi(\eta^*)d\eta^*,$ or the inner complex contour integral becomes $\int_{\partial X_j}\phi(\eta^*)\frac{\hat g(\eta^*)}{\eta^*-\eta}d\eta^*.$ Let $G(\eta^*)$ be defined as $G(\eta^*)=\phi(\eta^*)\hat g(\eta^*),$ then the integral becomes $\int_{\partial X_j}\frac{G(\eta^*)}{\eta^*-\eta}d\eta^*=2\pi iG(\eta)$ for $G(\eta^*)$ defined locally in $X_j.$ As such, $f(\eta)=2\pi i\int G(\eta)d\eta^*.$ To ensure $\int\int_{X_j}\eta^n\hat g(\eta)d\eta\wedge d\bar \eta\ne 0,$ we let $h_{\bar\Pi}=F,$ then $\hat g(\eta)=1/\eta^n$ and $f(\eta)=2\pi i\frac{\eta^*\phi(\eta)}{\eta^n}.$ In particular, we can compute $h_{\bar\Pi}$ and $F$ as follows,
$d(z\circ\eta)\wedge d(\bar z\circ\eta)=d\eta\wedge d\bar\eta\equiv\left[\frac{\partial\eta}{\partial z_1}dz_1+\frac{\partial \eta}{\partial z_2}dz_2\right]\wedge\left[\frac{\partial\bar\eta}{\partial z_1}dz_1+\frac{\partial \bar\eta}{\partial z_2}dz_2\right]=\left[\frac{\partial\eta}{\partial z_1}\frac{\partial\bar\eta}{\partial z_2}-\frac{\partial\eta}{\partial z_2}\frac{\partial\bar\eta}{\partial z_1}\right]dz_1\wedge dz_2$ for $F:=\frac{\partial\eta}{\partial z_1}\frac{\partial\bar\eta}{\partial z_2}-\frac{\partial\eta}{\partial z_2}\frac{\partial\bar\eta}{\partial z_1},$ where $\eta:=\eta(z_1,z_2).$ The integral of the meromorphic $(1,0)$-form $\omega$ over the 1-cycle $\partial U_j$ is
\begin{equation*}
\int_{\partial U_j}\omega=\int_{\partial X_j}\eta^*\omega=\int_{\partial X_j}g(\eta)\left(\frac{\partial\eta}{\partial z_1}dz_1+\frac{\partial\eta}{\partial z_2}dz_2\right).
\end{equation*} Then, by Stokes' theorem, the linear functional $\sigma_j(\chi)$ becomes
\begin{equation*}
\begin{split}
&\int_{\partial X_j}g(\eta(z_1, z_2))e^{-i\pi\langle\xi,\chi\rangle}\left(\frac{\partial\eta}{\partial z_1}\right)dz_1+g(\eta(z_1,z_2))e^{-i\pi\langle\xi,\chi\rangle}\left(\frac{\partial\eta}{\partial z_2}\right)dz_2\\&=\int_{X_j}\frac{\partial}{\partial z_1}\left(g(\eta(z_1, z_2))e^{-i\pi\langle\xi,\chi\rangle}\left(\frac{\partial\eta}{\partial z_2}\right)\right)-\frac{\partial}{\partial z_2}\left(g(\eta(z_1, z_2))e^{-i\pi\langle\xi,\chi\rangle}\left(\frac{\partial\eta}{\partial z_1}\right)\right)dz_1\wedge dz_2.
\end{split}
\end{equation*} A simple computation shows that this expression reduces to
\begin{equation*}
\begin{split}
&\int_{X_j}\left[\left(\frac{\partial\eta}{\partial z_2}\right)\left(\frac{\partial}{\partial z_1}\left(ge^{-i\pi\langle\xi,\chi\rangle}\right)\right)-\left(\frac{\partial\eta}{\partial z_1}\right)\left(\frac{\partial}{\partial z_2}\left(ge^{-i\pi\langle\xi,\chi\rangle}\right)\right)\right]dz_1\wedge dz_2\\
& =\int_{X_j}\left[\left(\frac{\partial\eta}{\partial z_1}\right)\left(-\frac{\partial}{\partial z_2}\left(ge^{-i\pi\langle\xi,\chi\rangle}\right)\right)-\left(\frac{\partial\eta}{\partial z_2}\right)\left(-\frac{\partial}{\partial z_1}\left(ge^{-i\pi\langle\xi,\chi\rangle}\right)\right)\right]dz_1\wedge dz_2
\end{split}
\end{equation*} for $\hat g(z)=\frac{h_{\bar\Pi}}{Fz^n}.$ Thus, by comparison of $F$ and $h_{\bar\Pi}$, to impose the condition of equality, we let $\frac{\partial\bar\eta}{\partial z_1}=-\frac{\partial}{\partial z_1}\left(ge^{-i\pi\langle\xi,\chi\rangle}\right)$ and $\frac{\partial\bar\eta}{\partial z_2}=-\frac{\partial}{\partial z_2}\left(ge^{-i\pi\langle\xi,\chi\rangle}\right).$ For $\frac{\partial\bar\eta}{\partial z_1}=-\frac{\partial}{\partial z_1}\left(ge^{-i\pi\langle\xi,\chi\rangle}\right),$ $\bar\eta:=-ge^{-i\pi\langle\xi,\chi\rangle}+\beta(z_2)$ and for $\frac{\partial\bar\eta}{\partial z_2}=-\frac{\partial}{\partial z_2}\left(ge^{-i\pi\langle\xi,\chi\rangle}\right),$ the solution is $\bar\eta:=-ge^{-i\pi\langle\xi,\chi\rangle}+\alpha(z_1).$ The two solutions must coincide, meaning that $\bar\eta:=-ge^{-i\pi\langle\xi,\chi\rangle}$ or $\eta=-\bar ge^{i\pi\langle\overline{\xi,\chi}\rangle}=-\bar ge^{i\pi\langle\chi,\xi\rangle},$ by the property of the Hermitian inner product. Such a diffeomorphism uniquely defines an equivalence and by extension, the complex manifold, $X,$ itself.

\section{Program for Genera and Isothermal Coordinates\label{sec:genera}}
We first prove that the Riemann surface $X$ is of genus $g=0$ by invoking the Gauss-Bonnet theorem from which it follows, by the uniformization theorem, that $\tilde X$ is holomorphically isomorphic to the complex projective line $\mathbb{CP}^1.$ For $X$ a compact two-dimensional (complex one-dimensional) Riemannian manifold with boundary $\partial X,$ let $K$ be the Gaussian curvature of $X$ and $k_g$ the geodesic curvature of $\partial X.$ Then, by Gauss-Bonnet,
\begin{equation*}
\int_XKdS+\int_{\partial X}k_gds=2\pi\chi(X)
\end{equation*} where $\chi(X)$ is the Euler characteristic of $X,$ $dS$ is the area surface element, and $ds$ is the line element. Note, it is assumed that $X$ is a compact manifold; however, this property  is proved in the latter geometric analysis [Section \ref{sec:Projective Geometry and Compactness}], namely $\mu(X)=\mu\left(\bigcup_j X_j\right)=\sum_j\mu(X_j)<+\infty$ is a finite measure and $X$ is the disjoint union of complex non-singular algebraic curves (for $X\cong\mathbb{CP}^1$ is shown to be true). Given that $X$ is a compact Riemann surface without boundary, the latter integral can be omitted such that 
$2\pi\chi(X)=\int_{X}KdS.$ Consider a Monge patch $\bm{F}:U\to{\mathbb{R}}^3\supset X$ for $X$ a compact Riemann surface, embedded in the ambient space ${\mathbb{R}}^3$, equipped with a Hermitian inner product and $U\subset\mathbb{CP}^1$ compact. More generally, let $\bm{F}=\bm{F}(z_1,z_2)$ for $\eta(z_1,z_2)=-\bar ge^{-i\pi\langle\xi,\chi\rangle}$ and the Hermitian inner product defined as $\langle\xi,\chi\rangle=\bar\xi\chi.$ Thus, $\bm{F}(z_1,z_2)=(Re\text{ }  w, Im\text{ }  w, w)=(x_1,x_2,z_1+iz_2)$ for $Re\text{ } w=x_1,$ $Im\text{ }  w=x_2$ and $w:=z_1+iz_2.$ Let $-g:=\Phi(z_1)\Theta(z_2),$ where $\langle\xi,\chi\rangle$ is independent of the local coordinates $(z_1,z_2),$ and $\Phi(z_1)=\omega(z_1)e^{i\alpha(z_1)},$ $\Theta(z_2)=\kappa(z_2)e^{i\beta(z_2)}.$ Let $\Gamma:\mathbb{R}\to\mathbb{R}$ be a linear transformation such that $z_2:=\Gamma z_1,$ inducing a linear transformation of coordinate vector fields on the tangent space $T_pX,$ $\hat{\Gamma}\frac{\partial}{\partial z_1}=\frac{\partial}{\partial z_2}.$ Thus, $\kappa(z_2)=\kappa(\Gamma z_1)$ and $\eta(z_1,z_2)=\Phi(z_1)\Theta(z_2)e^{i\pi\langle\chi,\xi\rangle}=\omega(z_1)e^{i\alpha(z_1)}\kappa(z_2)e^{i\beta(z_2)}e^{C\pi i}:=\phi(z_1)e^{i[\alpha(z_1)+\beta(z_2)+C\pi]}$ where $|\eta|=\phi(z_1)$ is a real-valued function expressed purely in terms of $z_1.$ As such, $\eta(z_1,z_2)=\phi(z_1)cos(\alpha(z_1)+\beta(z_2)+C\pi)+i\phi(z_1)sin(\alpha(z_1)+\beta(z_2)+C\pi)$ with $x_1=Re\text{ } \eta=\frac{\eta+\bar\eta}{2}=\phi(z_1)cos(\alpha(z_1)+\beta(z_2)+C\pi)$ and $x_2=Im\text{ } \eta=\frac{\eta-\bar\eta}{2i}=\phi(z_1)sin(\alpha(z_1)+\beta(z_2)+C\pi).$ For every coordinate system $\varphi:X\to\mathbb{C}$ and $a,b\in X$ the relation $\left[\varphi_{*}((e_1)_a),\varphi_{*}((e_2)_a)\right]=\mu_{\varphi(a)}$ holds if and only if $\left[\varphi_{*}((e_1)_b),\varphi_{*}((e_2)_b)\right]=\mu_{\varphi(b)}.$ If $x=\varphi(a)$ then for every $a\in X,$ the orientation $\mu_x=\left[\varphi_{*}((e_1)_a),\varphi_{*}((e_2)_a)\right]$ can be chosen consistently such that $X$ is orientable. Consequently, the Riemann surface $X$ can be parameterized in terms of $z_i,$ $i=1,2$ as 
\begin{equation*}
\bm{F}(z_1,z_2)=\begin{pmatrix}\phi(z_1)cos(\alpha(z_1)+\beta(z_2)+C\pi)\\ \phi(z_1)sin(\alpha(z_1)+\beta(z_2)+C\pi)\\z_1+iz_2\end{pmatrix}.
\end{equation*} Similarly, note that $|\eta|=\phi(z_1)$ such that $z_1=\phi^{-1}(|\eta|)$ for $\phi^{-1}$ admitting a continuous inverse, and likewise $arg\eta=\alpha(\phi^{-1}(|\eta|))+\beta(z_2)+C\pi$ or $\beta(z_2)=arg\eta-\alpha(\phi^{-1}(|\eta|))-C\pi$ for $z_2=\beta^{-1}(arg\eta-\alpha(\phi^{-1}(|\eta|))-C\pi)$ if $\beta$ admits a continuous inverse. Therefore, with $x_1=\frac{\eta+\bar\eta}{2}$ and $x_2=\frac{\eta-\bar\eta}{2i},$ the parameterization can be given globally in terms of $\eta$ as
\begin{equation}\label{equation: global parameterization}
\bm{F}(\eta)=\begin{pmatrix}(\eta+\bar\eta)/2\\(\eta-\bar\eta)/2i\\ \phi^{-1}(|\eta|)+i\beta^{-1}(arg\text{ }\eta-\alpha(\phi^{-1}(|\eta|))-C\pi) \end{pmatrix}.
\end{equation} The coefficients of the first fundamental form may be given as $E=g_{11}=\langle\bm{F}_{z_1},\bm{F}_{z_1}\rangle,$ $F=g_{12}=\langle\bm{F}_{z_1},\bm{F}_{z_2}\rangle,$ and $G=g_{22}=\langle\bm{F}_{z_2},\bm{F}_{z_2}\rangle$ where
\begin{equation*}
\bm{F}_{z_1}=\begin{pmatrix}\frac{d\phi}{dz_1}cos(\alpha(z_1)+\beta(z_2)+C\pi)-\phi(z_1)sin(\alpha(z_1)+\beta(z_2)+C\pi)\frac{d\alpha}{dz_1}\\ \frac{d\phi}{dz_1}sin(\alpha(z_1)+\beta(z_2)+C\pi)+\phi(z_1)cos(\alpha(z_1)+\beta(z_2)+C\pi)\frac{d\alpha}{dz_1}\\1\end{pmatrix},
\bm{F}_{z_2}=\begin{pmatrix} -\phi(z_1)sin(\alpha(z_1)+\beta(z_2)+C\pi)\frac{d\beta}{dz_2} \\ \phi(z_1)cos(\alpha(z_1)+\beta(z_2)+C\pi)\frac{d\beta}{dz_2} \\ i  \end{pmatrix}
\end{equation*} such that 
$E=g_{11}=\langle\bm{F}_{z_1},\bm{F}_{z_1}\rangle=1+\left(\frac{d\phi}{dz_1}\right)^2+\phi^2(z_1)\left(\frac{d\alpha}{dz_1}\right)^2, F=g_{12}=\langle\bm{F}_{z_1},\bm{F}_{z_2}\rangle=i+\phi^2(z_1)\frac{d\beta}{dz_2}\frac{d\alpha}{dz_1} \text{, and } G=g_{22}=\langle\bm{F}_{z_2},\bm{F}_{z_2}\rangle=\phi^2(z_1)\left(\frac{d\beta}{dz_2}\right)^2-1.$ It follows that after expansion of the terms $E,G,F,$ the Riemannian metric can be given locally as 
\begin{equation*}
\begin{split}
&g=ds^2=\left(1+\left(\frac{d\phi}{dz_1}\right)^2+\phi^2(z_1)\left(\frac{d\alpha}{dz_1}\right)^2\right)dz_1^2+2\left(i+\phi^2(z_1)\frac{d\beta}{dz_2}\frac{d\alpha}{dz_1}\right)dz_1dz_2\\&+\left(\phi^2(z_1)\left(\frac{d\beta}{dz_2}\right)^2-1\right)dz_2^2.
\end{split}
\end{equation*} Likewise, the Gaussian curvature is given by $K=\frac{eg-f^2}{EG-F^2}$ where $E=g_{11}=\langle\bm{F}_{z_1},\bm{F}_{z_1}\rangle,$ $F=g_{12}=\langle\bm{F}_{z_1},\bm{F}_{z_2}\rangle,$ and $G=g_{22}=\langle\bm{F}_{z_2},\bm{F}_{z_2}\rangle$ are coefficients of the first fundamental form and $e=\langle \bm{N}, \bm{F}_{z_1z_1} \rangle,$ $f=\langle \bm{N}, \bm{F}_{z_1z_2} \rangle,$ and $g=\langle \bm{N}, \bm{F}_{z_2z_2} \rangle$ are coefficients of the second fundamental form, where $\bm{N}$ is the normal vector $\bm{N}=\bm{F}_{z_1}\times\bm{F}_{z_2}.$ A simple calculation in the standard basis gives the second order partial derivatives,
\begin{equation*}
\begin{split}
&\bm{F}_{z_1z_1}=\biggl(\frac{d^2\phi}{dz_1^2}cos(\alpha(z_1)+\beta(z_2)+C\pi)-2\frac{d\phi}{dz_1}\frac{d\alpha}{dz_1}sin(\alpha(z_1)+\beta(z_2)+C\pi)-\phi(z_1)cos(\alpha(z_1)\\&+\beta(z_2)+C\pi)\left(\frac{d\alpha}{dz_1}\right)^2-\phi(z_1)sin(\alpha(z_1)+\beta(z_2)+C\pi)\frac{d^2\alpha}{dz_1^2},\frac{d^2\phi}{dz_1^2}sin(\alpha(z_1)+\beta(z_2)+C\pi)\\&+2\frac{d\phi}{dz_1}\frac{d\alpha}{dz_1}cos(\alpha(z_1)+\beta(z_2)+C\pi)-\phi(z_1)sin(\alpha(z_1)+\beta(z_2)+C\pi)\left(\frac{d\alpha}{dz_1}\right)^2+\phi(z_1)\\&cos(\alpha(z_1)+\beta(z_2)+C\pi)\frac{d^2\alpha}{dz_1^2},0\biggr),
\end{split}
\end{equation*}
\begin{equation*}
\bm{F}_{z_2z_2}=\begin{pmatrix} -\phi(z_1)cos(\alpha(z_1)+\beta(z_2)+C\pi)\left(\frac{d\beta}{dz_2}\right)^2-\phi(z_1)sin(\alpha(z_1)+\beta(z_2)+C\pi)\frac{d^2\beta}{dz_2^2} \\ \phi(z_1)cos(\alpha(z_1)+\beta(z_2)+C\pi)\frac{d^2\beta}{dz_2^2}-\phi(z_1)sin(\alpha(z_1)+\beta(z_2)+C\pi)\left(\frac{d\beta}{dz_2}\right)^2 \\ 0 \end{pmatrix},
\end{equation*} and 
\begin{equation*}
\bm{F}_{z_1z_2}=\begin{pmatrix} -\frac{d\phi}{dz_1}\frac{d\beta}{dz_2}sin(\alpha(z_1)+\beta(z_2)+C\pi)-\phi(z_1)\frac{d\alpha}{dz_1}\frac{d\beta}{dz_2}cos(\alpha(z_1)+\beta(z_2)+C\pi) \\ \frac{d\phi}{dz_1}\frac{d\beta}{dz_2}cos(\alpha(z_1)+\beta(z_2)+C\pi)-\phi(z_1)\frac{d\alpha}{dz_1}\frac{d\beta}{dz_2}sin(\alpha(z_1)+\beta(z_2)+C\pi)\\  0 \end{pmatrix}.
\end{equation*} Furthermore, the normal vector $\bm{N}$ to the orientable Riemann surface $X$ with induced orientation $\partial\mu$ is given by $\bm{N}=\frac{\partial \bm{F}}{\partial z_1}\times\frac{\partial\bm{F}}{\partial z_2}$
\begin{equation*}
\begin{split}
&=\biggl(i\frac{d\phi}{dz_1}sin(\alpha(z_1)+\beta(z_2)+C\pi)+i\phi(z_1)cos(\alpha(z_1)+\beta(z_2)+C\pi)\frac{d\alpha}{dz_1}-\phi(z_1)cos(\alpha(z_1)\\&+\beta(z_2)+C\pi)\frac{d\beta}{dz_2},i\phi(z_1)sin(\alpha(z_1)+\beta(z_2)+C\pi)\frac{d\alpha}{dz_1}-i\frac{d\phi}{dz_1}cos(\alpha(z_1)+\beta(z_2)+C\pi)\\&-\phi(z_1)sin(\alpha(z_1)+\beta(z_2)+C\pi)\frac{d\beta}{dz_2},\phi(z_1)\frac{d\phi}{dz_1}\frac{d\beta}{dz_2} \biggr).
\end{split}
\end{equation*} Lastly, we calculate the coefficients of the second fundamental form, which can in fact be expressed as $e=\langle \bm{N}, \bm{F}_{z_1z_1} \rangle,$ $f=\langle \bm{N}, \bm{F}_{z_1z_2} \rangle,$ and $g=\langle \bm{N}, \bm{F}_{z_2z_2} \rangle$ with
\begin{equation*}
\begin{split}
&e=\phi^2(z_1)\left(\frac{d\beta}{dz_2}\right)\left(\frac{d\alpha}{dz_1}\right)^2-\phi(z_1)\left(\frac{d\beta}{dz_2}\right)\left(\frac{d^2\phi}{dz_1^2}\right)-i\phi^2(z_1)\left(\frac{d\alpha}{dz_1}\right)^3+i\phi(z_1)\left(\frac{d\alpha}{dz_1}\right)\left(\frac{d^2\phi}{dz_1^2}\right)\\&-i\phi(z_1)\left(\frac{d\phi}{dz_1}\right)\left(\frac{d^2\alpha}{dz_1^2}\right)-i\phi(z_1)\left(\frac{d\phi}{dz_1}\right)\left(\frac{d\alpha}{dz_1}\right)^2sin(2\alpha(z_1)+2\beta(z_1)+2\pi C)\\&-2i\left(\frac{d\phi}{dz_1}\right)^2\left(\frac{d\alpha}{dz_1}\right),
\end{split}
\end{equation*}
\begin{equation*}
\begin{split}
g=-i\phi(z_1)\left(\frac{d\phi}{dz_1}\right)\left(\frac{d^2\beta}{dz_2^2}\right)-i\phi^2(z_1)\left(\frac{d\alpha}{dz_1}\right)\left(\frac{d\beta}{dz_2}\right)^2+\phi^2(z_1)\left(\frac{d\beta}{dz_2}\right)^3,
\end{split}
\end{equation*} and 
\begin{equation*}
\begin{split}
f=-i\left(\frac{d\phi}{dz_1}\right)^2\left(\frac{d\beta}{dz_2}\right)-i\phi^2(z_1)\left(\frac{d\alpha}{dz_1}\right)^2\left(\frac{d\beta}{dz_2}\right)+\phi^2(z_1)\left(\frac{d\alpha}{dz_1}\right)\left(\frac{d\beta}{dz_2}\right)^2.
\end{split}
\end{equation*} Thus, the Gaussian curvature $K$ is given by $\frac{eg-f^2}{EG-F^2}$ [see Appendix \ref{sec:Gaussian Curvature} for an explicit calculation], for 
\begin{equation}
EG-F^2=\phi^2(z_1)\left[\left(\frac{d\beta}{dz_2}\right)^2+\left(\frac{d\phi}{dz_1}\right)^2\left(\frac{d\beta}{dz_2}\right)^2-\left(\frac{d\alpha}{dz_1}\right)^2\right]-\left(\frac{d\phi}{dz_1}\right)^2-2i\frac{d\beta}{dz_2}\frac{d\alpha}{dz_1}.
\end{equation}

Consequently, if we impose the condition $\iint\limits_X KdS=2\pi\chi(X)\equiv4\pi$ or 
\begin{equation*}
\int_{X} \frac{K}{\bm{\hat{n}}}dz_1\wedge dz_2\equiv\int_{X} \frac{K}{\bm{\hat{n}}}dz_1dz_2=4\pi
\end{equation*} for $\bm{\hat{n}}$ the unit normal vector to $X,$ then the Riemann surface $X$ must have genus $g=0,$ implying that it is a simply connected, complex manifold. If the map $\eta$ is in fact biholomorphic, then we have proven that $X$ is diffeomorphic (i.e. holomorphically isomorphic) to $\mathbb{CP}^1.$ In particular, every Riemann surface is the quotient group of the universal covering by a discrete group, such that $X=\tilde X/\Gamma$ for $\tilde X=\mathbb{CP}^1$ the universal covering of $X$ and $\Gamma$ the set of deck transformations $Aut_{\rho}(X)$ acting transitively on $X,$ for a universal covering map $\eta^{-1}=\rho:\tilde X\to X$ the identity. Therefore for all purposes, conformally equivalent Riemann surfaces can be identified with one another (i.e. they are identical). 

\section{Biholomorphicity\label{sec:biholomorphicity}}
The map $\eta:X\to\mathbb{CP}^1$ given by $\eta(z_1,z_2)=\phi(z_1)e^{i[\alpha(z_1)+\beta(z_2)+\pi \langle \chi,\xi\rangle]}$ is a bijective holomorphic function (a biholomorphism) if and only if the Wirtinger derivative with respect to the conjugate local coordinate vanishes, that is if $\frac{\partial\eta}{\partial \bar w}=0.$ In particular, let $\eta\in\Gamma(\Lambda^0(T^*X))$ for $\Gamma(\Lambda^0)$ sections over the structure sheaf $\mathcal{O}_X$ of holomorphic functions on $X.$ Then $\eta\in\mathcal{O}_X$ if and only if the Wirtinger derivative with respect to the conjugate of the local complex coordinate $w=z_1+iz_2$ vanishes. Recall that the Wirtinger derivatives are defined as linear partial differential operators of first order, such that $\frac{\partial}{\partial w}=\frac{1}{2}\left(\frac{\partial}{\partial z}_1-i\frac{\partial}{\partial z_2}\right)$ and $\frac{\partial}{\partial \bar w}=\frac{1}{2}\left(\frac{\partial}{\partial z_1}+i\frac{\partial}{\partial z_2}\right).$ Then
\begin{equation*}
\frac{\partial\eta}{\partial z_1}=\frac{d\phi}{dz_1}e^{i[\alpha(z_1)+\beta(z_2)+\pi\langle\chi,\xi\rangle]}+i\phi(z_1)\frac{d\alpha}{dz_1}e^{i[\alpha(z_1)+\beta(z_2)+\pi\langle\chi,\xi\rangle]} 
\end{equation*} and
\begin{equation*}
\frac{\partial\eta}{\partial z_2}=i\phi(z_1)\frac{d\beta}{dz_2}e^{i[\alpha(z_1)+\beta(z_2)+\pi\langle\chi,\xi\rangle]}.
\end{equation*} The Wirtinger condition then becomes
\begin{equation*}
\frac{d\phi}{dz_1}+i\phi(z_1)\left(\frac{d\alpha}{dz_1}+i\frac{d\beta}{dz_2}\right)=0,
\end{equation*} which yields a solution $|\eta|=e^{\lambda_0}e^{-i\alpha(z_1)}e^{z_1\frac{d\beta}{dz_2}}.$ But $\phi(z_1)=|\bar g|>0$ a priori, and thus the chart is biholomorphic if and only if $\phi(z_1)=\lambda e^{z_1\frac{d\beta}{dz_2}-i\alpha(z_1)}$ for $\lambda\in \mathbb{C},$ choosing the positive branch.

The calculation for Gaussian curvature $K$ in the above is rather tedious, when indeed under certain isothermal coordinate conditions this expression may be simplified considerably. We now develop the metric theory of the complex manifold to obtain precise expressions for the Gaussian curvature $K.$ The Riemannian metric on $X$ is given locally as $ds^2=Edz_1^2+2Fdz_1dz_2+Gdz_2^2.$ In complex (local) coordinates, it assumes the form $ds^2=\lambda |dz+\mu d\bar z|^2,$ for $\lambda=\frac{1}{4}\left(E+G+2\sqrt{EG-F^2}\right),$ $\mu=(E-G+2iF)/4\lambda$ such that $\lambda$ and $\mu$ are smooth, satisfying the conditions $\lambda>0$ and $|\mu|<1.$ For isothermal coordinates $(u,v)$ the metric has the form $ds^2=\rho(du^2+dv^2)$ for $\rho>0$ smooth. Likewise, the complex coordinate $w=u+iv$ satisfies $\rho|dw|^2=\rho\left|\frac{\partial w}{\partial z}\right|^2\left|dz+\frac{w_{\bar z}}{w_z}d\bar z\right|^2.$ That is, $(u,v)$ are local isothermal coordinates if and only if they satisfy the Beltrami equation $\frac{\partial w}{\partial \bar z}=\mu\frac{\partial w}{\partial z},$ i.e. the equation gives rise to a diffeomorphic solution. Recall that the metric is given locally on $X$ by 
\begin{equation*}
\begin{split}
&g=ds^2=\left(1+\left(\frac{d\phi}{dz_1}\right)^2+\phi^2(z_1)\left(\frac{d\alpha}{dz_1}\right)^2\right)dz_1^2+2\left(i+\phi^2(z_1)\frac{d\beta}{dz_2}\frac{d\alpha}{dz_1}\right)dz_1dz_2\\&+\left(\phi^2(z_1)\left(\frac{d\beta}{dz_2}\right)^2-1\right)dz_2^2.
\end{split}
\end{equation*} If $(z_1,z_2)$ are local isothermal coordinates then the Riemannian metric assumes the form $g=ds^2=\varphi(dz_1^2+dz_2^2),$ where $\varphi$ is a smooth $C^{\infty}$-differentiable function; that is, the metric is conformal to the Euclidean metric. Therefore, we consider the case for which $F=g_{12}=g_{21}=0$ and $E=G,$ for $g$ conformally equivalent to the Euclidean metric, such that the metric tensor is a positive definite symmetric rank two tensor
\begin{equation*}
\left[g_{ij}\right]=\left(\begin{array}{ccc}
E & 0 \\ 0 & E  
\end{array}
\right)=E\left[\delta_{ij}\right]
\end{equation*} for $det[g_{ij}]=E^2.$ Then $F=0$ if and only if $\phi^2(z_1)\frac{d\beta}{dz_2}\frac{d\alpha}{dz_1}=-i.$ Similarly, $E=G$ if and only if $2+\left(\frac{d\phi}{dz_1}\right)^2+\phi^2(z_1)\left[\left(\frac{d\alpha}{dz_1}\right)^2-\left(\frac{d\beta}{dz_2}\right)^2\right]=0$ with $E=1+\left(\frac{d\phi}{dz_1}\right)^2+\phi^2(z_1)\left(\frac{d\alpha}{dz_1}\right)^2.$ As such, if these conditions are satisfied and if we define $\Gamma:=E=1+\left(\frac{d\phi}{dz_1}\right)^2+\phi^2(z_1)\left(\frac{d\alpha}{dz_1}\right)^2$ then $(z_1,z_2)$ are necessarily isothermal coordinates and $g=\Gamma(dz_1^2+dz_2^2)\in[g]=\{g|h=\lambda^2h\text{ for }\lambda\text{ a smooth real-valued function}\}.$ In isothermal coordinates, the Gaussian curvature takes the form
\begin{equation*}
K=-\frac{1}{2}e^{-log\Gamma}\Delta log\Gamma.
\end{equation*} More precisely, if we let $\rho=log\Gamma,$ then
\begin{equation*}
K=-\frac{1}{2}e^{-\rho}\left(\frac{\partial^2\rho}{\partial z_1^2}+\frac{\partial^2\rho}{\partial z_2^2}\right)
\end{equation*} given that the following set of differential equations are satisfied:
\begin{equation}\label{equation: coupled differential equations}
\phi^2(z_1)\frac{d\beta}{dz_2}\frac{d\alpha}{dz_1}=-i\text{  , }2+\left(\frac{d\phi}{dz_1}\right)^2+\phi^2(z_1)\left[\left(\frac{d\alpha}{dz_1}\right)^2-\left(\frac{d\beta}{dz_2}\right)^2\right]=0.
\end{equation} Therefore the condition of genus $g=0$ reduces to
\begin{equation}
\chi(X)=-\frac{1}{4\pi}\iint\limits_X\frac{\Delta\rho}{e^{\rho}}dS=-\frac{1}{4\pi}\iint\limits_X\frac{\Delta log\Gamma}{\Gamma}dS=2.
\end{equation} 
In isothermal coordinates, the calculation of Gaussian curvature is reduced significantly to \begin{equation*}K=-\frac{1}{2}e^{-log\Gamma}\left(\frac{\partial^2}{\partial z_1^2}log\Gamma+\frac{\partial^2}{\partial z_2^2}log\Gamma\right)\end{equation*} [see Appendix \ref{sec:Gaussian Curvature}]. The passage of isothermal coordinates $(z_1,z_2)$ can also be realized for the global parameterization
\begin{equation*}
\bm{F}(\eta)=\begin{pmatrix}(\eta+\bar\eta)/2\\(\eta-\bar\eta)/2i\\ \phi^{-1}(|\eta|)+i\beta^{-1}(arg\text{ }\eta-\alpha(\phi^{-1}(|\eta|))-C\pi) \end{pmatrix}
\end{equation*} by computing the coefficients of the first fundamental form $E,F,G$ in terms of $\eta.$ Under the assumption that $(z_1,z_2)$ are local isothermal coordinates for $X,$ the Beltrami equation is $\frac{\partial \eta}{\partial\bar\eta}=\mu$ where $\frac{\partial \eta}{\partial\bar\eta}=\frac{1}{2}\left(\frac{\partial\eta}{\partial z_1}+i\frac{\partial\eta}{\partial z_2}\right)=0$ for $\mu=(E-G+2iF)/4\lambda.$ It follows that $(z_1,z_2)=z_1+iz_2$ are isothermal coordinates provided $E-G+2iF=0.$ This is equivalent to the coupled set of differential equations Eq. \ref{equation: coupled differential equations},
\begin{equation*}
\left(\frac{d\phi}{dz_1}\right)^2+\phi^2(z_1)\left[\left(\frac{d\alpha}{dz_1}\right)^2-\left(\frac{d\beta}{dz_2}\right)^2+2i\frac{d\beta}{dz_2}\frac{d\alpha}{dz_1}\right]=0,
\end{equation*} obtained without invoking the Beurling transform to solve the Beltrami equation. More precisely, the partial differential equation $\frac{\partial w}{\partial\bar z}=\mu\frac{\partial w}{\partial z}$ is the the Beltrami equation for $w$ a complex distribution of the complex variable $z$ defined in an open subset $U\subset\mathbb{C}$ for $\mu$ a complex $L^{\infty}(U)$-function of norm less than unity \cite[Pg. 314-317]{Spivak}. In particular for $U\subset\mathbb{C},$ the Riemannian metric $X$ assumes the form 
\begin{equation*}
g=Edz_1^2+2Fdz_1dz_2+Gdz_2^2
\end{equation*} with the Beltrami coefficient 
\begin{equation*}
\mu=\frac{E-G+2iF}{E+G+2\sqrt{EG-F^2}}
\end{equation*} satisfying the norm property
\begin{equation}
|\mu|^2=\frac{E+G-2\sqrt{EG-F^2}}{E+G+2\sqrt{EG-F^2}}<1. 
\end{equation} We now obtain a property on sections $g\in\mathcal{M}\Omega^{(0,0)}(\mathbb{CP}^1)$ by imposing the condition of isothermal coordinates $(z_1,z_2)$ by the following constuct: With respect to the global parameterization Eq. \ref{equation: global parameterization} for the local complex coordinate $w$ on $X,$ $w:=\phi^{-1}(|\eta|)+i\beta^{-1}(arg \text{ }\eta-\alpha(\phi^{-1}(|\eta|))-C\pi)=\phi^{-1}(|\bar\eta|)+i\beta^{-1}(-arg \text{ }\bar\eta-\alpha(\phi^{-1}(|\bar\eta|))-C\pi)$ for $|\eta|=|\bar\eta|$ and $arg\text{ }\bar\eta=-arg\text{ }\eta.$ Therefore, $\frac{\partial w}{\partial\bar\eta}=\frac{d\phi^{-1}}{d|\bar\eta|}\frac{\bar\eta}{|\bar\eta|}+i\frac{d\beta^{-1}}{d\beta}\left(-\frac{d}{d\bar\eta}arg\text{ }\bar\eta-\frac{d\alpha}{d\phi^{-1}}\frac{d\phi^{-1}}{d|\bar\eta|}\frac{\bar\eta}{|\bar\eta|}\right)=\frac{d\phi^{-1}}{d|\eta|}\frac{\bar\eta}{|\eta|}+i\frac{d\beta^{-1}}{d\beta}\left(\frac{i}{\bar\eta}-\frac{i\bar\eta}{|\eta|^2}-\frac{d\alpha}{d\phi^{-1}}\frac{d\phi^{-1}}{d|\eta|}\frac{\bar\eta}{|\eta|}\right)$ for $\frac{d}{d\bar\eta}arg\text{ }\bar\eta=i\left(\frac{\bar\eta}{|\bar\eta|^2}-\frac{1}{\bar\eta}\right).$ Likewise, $\frac{\partial w}{\partial\eta}=\frac{d\phi^{-1}}{d|\eta|}\frac{\eta}{|\eta|}+i\frac{d\beta^{-1}}{d\beta}\left(\frac{i\eta}{|\eta|^2}-\frac{i}{\eta}-\frac{d\alpha}{d\phi^{-1}}\frac{d\phi^{-1}}{d|\eta|}\frac{\eta}{|\eta|}\right).$ Then the Beltrami equation $\frac{\partial w}{\partial\bar\eta}=\mu\frac{\partial w}{\partial \eta}$ reduces to 
\begin{equation*}
\frac{1}{|\eta|}\frac{d\phi^{-1}}{d|\eta|}(\bar\eta-\mu\eta)+i\frac{d\beta^{-1}}{d\beta}\left(i\left(\frac{1}{\bar\eta}+\frac{\mu}{\eta}\right)-\frac{i}{|\eta|^2}\left(\bar\eta-\mu\eta\right)-\frac{d\alpha}{d\phi^{-1}}\frac{d\phi^{-1}}{d|\eta|}\frac{1}{|\eta|}\left(\bar\eta-\mu\eta\right)\right)=0,
\end{equation*} which is satisfied by the sufficient condition $-\frac{\eta}{\bar\eta}=\frac{\bar\eta}{\eta},$ namely $-\eta^2=\bar\eta^2$ if and only if $g^2=-{\bar g}^2.$ If the sufficient condition $g^2=-{\bar g}^2$ is satisfied then the Beltrami equation holds, meaning that $(z_1,z_2)$ are isothermal coordinates. This simplifies the condition of $g=0$ to
\begin{equation}
\chi(X)=-\frac{1}{4\pi}\iint\limits_X\frac{\Delta\rho}{e^{\rho}}dS=-\frac{1}{4\pi}\iint\limits_X\frac{\Delta log\Gamma}{\Gamma}dS=2,
\end{equation} yielding the following theorem:
\begin{theorem}
If the sufficient condition $g^2=-{\bar g}^2$ is satisfied for $g\in\mathcal{M}\Omega^{(0,0)}(\mathbb{CP}^1),$ then the Beltrami equation is satisfied and $(z_1,z_2)$ must necessarily be local isothermal coordinates on $X,$ such that \begin{equation}
\chi(X)=-\frac{1}{4\pi}\iint\limits_X\frac{\Delta\rho}{e^{\rho}}dS=-\frac{1}{4\pi}\iint\limits_X\frac{\Delta log\Gamma}{\Gamma}dS=2
\end{equation} implies that the Riemann surface has genus zero.
\end{theorem} 

\section{Quasidistributional Analysis\label{sec:quasidsitribution}}
We now give an interpretation of the complex local coordinate distribution $w$ in a quasi-distributional sense. In particular \cite[Pg. 65]{Ahlfors1},
\begin{definition}
If $f: D\to E$ is an orientation preserving homeomorphism between two open sets in the complex plane $\mathbb{C},$ and if f is continuously differentiable then it is said to be a $K$-quasiconformal mapping if the derivative of $f$ at every point maps circles to ellipeses with bounded eccentricity $K.$
\end{definition} That is, for the complex local coordinate $w=z_1+iz_2$ and for $z'=x_1+ix_2,$ suppose $w: D\to E$ for $z_i,$ $i=1,2$ coordinate functions variably dependent on other complex coordinates $z=\phi$ for $\phi:X\to\mathbb{C}.$ Then if $w$ is assumed to have continuous partial derivatives, $w$ is a $\mu$-quasiconformal mapping provided it satisfies the Beltrami equation
\begin{equation}\label{equation: distributional Beltrami}
\frac{\partial w}{\partial \bar z}=\mu(z)\frac{\partial w}{\partial z}
\end{equation} for a complex-valued Lebesgue measurable $\mu$ satisfying the norm condition $|\mu|^2<1$ or $sup\text{ } |\mu|<1.$ 
Such a quasiconformal theory has a geometric interpretation. If we equip $D$ (identified with an open subset of $X$) with the metric
\begin{equation}\label{equation: complex metric}
g=ds^2=Edz_1^2+2Fdz_1dz_2+Gdz_2^2=\lambda(\phi)|d\phi+\mu(\phi)d\bar\phi|^2,
\end{equation} for $\lambda>0$ and $|\mu|<1,$ then $w$ satisfies Eq. \ref{equation: distributional Beltrami} if and only if it is a conformal transformation from $D,$ equipped with the metric given by Eq. \ref{equation: complex metric}, to the domain $E$ equipped with the standard Euclidean metric. In the above analysis, we treated $w=\phi=z$ to verify the assumption that $(z_1,z_2)$ were isothermal coordinates, such that $D$ was, itself, equipped with a Euclidean metric if and only if one of the following conditions was satisfied: $g^2=-{\bar g}^2$ for $g\in\mathcal{M}\Omega^{(0,0)}(\mathbb{CP}^1)$ or $\left(\frac{d\phi}{dz_1}\right)^2+\phi^2(z_1)\left[\left(\frac{d\alpha}{dz_1}\right)^2-\left(\frac{d\beta}{dz_2}\right)^2+2i\frac{d\beta}{dz_2}\frac{d\alpha}{dz_1}\right]=0$; that is, the Beltrami equation reduced to $\frac{\partial \eta}{\partial \bar \eta}=\mu(z).$ For isothermal coordinates $(z_1,z_2)$ the metric assumes the form $ds^2=\rho(dz_1^2+dz_2^2),$ where the complex variable $w=z_1+iz_2$ satisfies $\rho|dw|^2=\rho|w_w|^2\left|dw+\frac{w_{\bar w}}{w_w}d\bar w\right|^2=\rho|dw+w_{\bar w}d\bar w|^2,$ which relative to the Euclidean metric $dzd\bar z=dwd\bar w$ has eigenvalues
\begin{equation}
(1+|\mu|^2)^2\left|\frac{\partial w}{\partial w}\right|^2=1, \text{  } \left(1-|\mu|\right)^2\left|\frac{\partial w}{\partial w}\right|^2=1.
\end{equation} The repeated eigenvalues correspond to the squared distance of the major and minor axis of the ellipse realized by pulling back the unit circle along $w.$ Thus, $w$ is $\mu$-quasiconformal, where $\mu=0$ is the condition for $(z_1,z_2)$ to be isothermal coordinates. It follows that the dilatation of $w$ at a point $z=w$ is given by 
\begin{equation}
K(z)=\frac{1+|\mu(z)|}{1-|\mu(z)|}=1.
\end{equation} Furthermore, the essential supremum norm of $K$ is
\begin{equation}
K=\sup_{z\in D}|K(z)|=\frac{1+||\mu||_{\infty}}{1-||\mu||_{\infty}}=1.
\end{equation} This complex distribution satisfies the Measurable Riemann Mapping theorem, 
\begin{theorem}[Measurable Riemann Mapping Theorem]\label{theorem: measurable riemann}
For $\mu(z)$ a bounded, measurable function with $z\in\mathbb{C}$ and the $L^{\infty}$-norm $||\mu||_{\infty}<1,$ there exists a quasiconformal homeomorphism that is a solution of the Beltrami equation
\begin{equation}
\frac{\partial f}{\partial\bar z}=\mu(z)\frac{\partial f}{\partial z},
\end{equation} where $f$ has fixed points $0, 1,$ and $\infty.$
\end{theorem}

\section{Hodge Theory\label{sec: hodge}}
In this section we analyze the genus condition of $\int_X e^{-log\Gamma}\Delta\rho dS=-8\pi$ or $\int_Xe^{-log\Gamma}\Delta\rho+\frac{\kappa}{\Sigma}dS=0,$ whereby we let $\frac{1}{\Sigma}\int_X dS=1$ and $\kappa=8\pi.$ 
Recall that the area of the Riemann surface is given in terms of the coefficients of the first fundamental form; that is, $\Sigma=\iint\limits_X\sqrt{EG-F^2}dz_1dz_2.$ Let the Hodge star operator $\star:\Lambda^{k}(V)\to\Lambda^{n-k}(V)$ be an isomorphism from the $k$-th exterior power space to the $(n-k)$-th exterior power space, such that the composition with itself maps $\star\circ\star:\Lambda^{k}(V)\to\Lambda^{k}(V).$ The former space $\Lambda^{k}(V)$ has dimension $\binom{n}{k}$, while the latter $\Lambda^{n-k}(V)$ has dimension $\binom{n}{n-k},$ which are equal, inducing an isomorphism. The Hodge dual operator induces a natural nondegenerate symmetric bilinear form (i.e. an inner product) on the exterior power space of $k$-vectors, i.e. on the exterior algebra $\Lambda(V).$ 
Let $\eta$ and $\zeta$ be two $k$-forms, then $\zeta\wedge\star\eta=\langle\zeta,\eta\rangle\omega$ for $\omega$ the normalized $n$-form (whereby $\omega\wedge\star\omega=\omega$) \cite{Darling}. The normalized $n$-form is simply the volume $n$-form for a Riemannian manifold, i.e. $\omega=\sqrt{|det[g_{ij}]|}dx^1\wedge...\wedge dx^n$ in local coordinates $x^i$ for $[g_{ij}]$ the metric tensor on the manifold. For the Hodge star applied twice $\star\circ\star:\Lambda^{k}(V)\to\Lambda^{k}(V),$ we define the identity element up to sign, $\star\star\eta=(-1)^{k(n-k)}s\eta$ for $s=sgn(det[g_{ij}])$ with $\eta\in\Lambda^k(V)$ and $V$ an $n$-dimensional space. Therefore, the inverse of the Hodge star operator can be defined as 
\begin{equation}
\begin{cases} \star^{-1}:\Lambda^{k}(V)\to\Lambda^{(n-k)}(V) \\ \eta\mapsto (-1)^{k(n-k)}s\star\eta\end{cases}.
\end{equation} We hence define the Hodge dual of a differential $k$-form in the cotangent space of an $n$-dimensional Riemannian manifold.
\begin{definition}
For $\eta$ and $\zeta$ both differential $k$-forms of the cotangent space of an $n$-dimensional, orientable (pseudo-)Riemannian manifold, $\star\eta$ is the unique differential $(n-k)$-form such that
\begin{equation}
\eta\wedge\star\zeta=\langle\eta,\zeta\rangle\omega
\end{equation} for $\omega$ the differential $n$-form. The $L^2$ inner product of differential $k$-forms can be defined in terms of the Hodge dual as
\begin{equation}
(\eta,\zeta)=\int_M\eta\wedge\star\zeta
\end{equation} for $\eta$ and $\zeta$ both $k$-forms in the section $\Omega^k(M):=\Gamma(\Lambda^k(T^*M)).$
\end{definition} Lastly, to fully develop the theory, we define the codifferential of $k$-forms. In particular for $d:\Omega^k(M)\to\Omega^{k+1}(M),$ the codifferential $\delta:\Omega^k(M)\to\Omega^{k-1}(M)$ is defined as $\delta=(-1)^{n(k-1)+1}s\star d\star=(-1)^k\star^{-1}d\star.$ In fact, for $\zeta$ a $(k+1)$-form and $\eta$ a $k$-form,
\begin{equation} 
\int_Md(\eta\wedge\star\zeta)=\int_M(d\eta\wedge\star\zeta-\eta\wedge\star(-1)^{(k+1)}\star^{-1}d\star\zeta)=(d\eta,\zeta)-(\eta,\delta\zeta)=0,
\end{equation}
such that the codifferential is the adjoint of the exterior derivative, i.e. $(\eta,\delta\zeta)=(d\eta,\zeta).$ The above integral is zero for $M$ of empty boundary or with $\eta$ and $\star\zeta$ assuming zero values on the boundary $\partial M.$ Consider the condition on genera, $\int_Xe^{-log\Gamma}\Delta\rho+\frac{\kappa}{\Sigma}dS=0.$ Then the area element is given by the modulus of the exterior product, $dS=|dz_1\wedge dz_2|\equiv\frac{|dz_1\wedge dz_2|}{dz_1\wedge dz_2}dz_1\wedge dz_2=\frac{|dz_1\wedge dz_2|}{dz_1\otimes dz_2+I}(dz_1\otimes dz_2+I)$ where the exterior algebra is generated by the tensor algebra $T(V)$ quotiented by the two-sided ideal, formed by elements $\alpha\otimes\alpha$ for $\alpha\in V,$ $\Lambda(V)=T(V)/I$ for $\alpha\wedge\beta=\alpha\otimes\beta+I.$ As such, let the unit normal vector of the Riemann surface $X$ at $(z_1,z_2)$ be given by $\bm{\hat n}:=\Phi(z_1,z_2),$ for which 
\begin{equation}
\int_Xe^{-log\Gamma}\Delta\rho+\frac{\kappa}{\Sigma}dS=\int_X\frac{e^{-log\Gamma}\Delta\rho}{\Phi(z_1,z_2)}+\frac{\kappa}{\Sigma\Phi(z_1,z_2)}dz_1\wedge dz_2=\int_{X}\frac{\Sigma\Delta\rho+\Gamma\kappa}{\Gamma\Sigma\Phi(z_1,z_2)}dz_1\wedge dz_2=0.
\end{equation} Likewise, let 
\begin{equation*}
\int_{X}\frac{\Sigma\Delta\rho+\Gamma\kappa}{\Gamma\Sigma\Phi(z_1,z_2)}dz_1\wedge dz_2=\int_Xd(\eta\wedge\star\zeta)=0.
\end{equation*} Then let $\omega=\eta\wedge\star\zeta=\langle\eta,\zeta\rangle\varpi$ for $\eta$ a $k$-form, $\zeta$ a $(k+1)$-form, $\star\zeta$ a $(n-k-1)$-form and $\varpi\in\Gamma(\Lambda^{dimV}(T^*X))$ where $dimV=dim(\Lambda^{n-k-1}(V)\wedge\Lambda^k(V))=n-1=1$ for $n=dim_{\mathbb{R}}X.$ Suppose
\begin{equation*}
d\omega=d(\eta\wedge\star\zeta)=\frac{\Sigma\Delta\rho+\Gamma\kappa}{\Gamma\Sigma\Phi(z_1,z_2)}dz_1\wedge dz_2,
\end{equation*} then let $\omega=fdz_2$ such that the exterior derivative becomes $\frac{\partial f}{\partial z_1}dz_1\wedge dz_2=\frac{\Sigma\Delta\rho+\Gamma\kappa}{\Gamma\Sigma\Phi(z_1,z_2)}dz_1\wedge dz_2$ or 
\begin{equation*}
\frac{\partial f}{\partial z_1}=\frac{\Sigma\Delta\rho+\Gamma\kappa}{\Gamma\Sigma\Phi(z_1,z_2)}.
\end{equation*} Consequently, the solution is of the form
\begin{equation}
f(z_1,z_2)=\int\frac{\Sigma\Delta\rho+\Gamma\kappa}{\Gamma\Sigma\Phi(z_1,z_2)}dz_1+\Psi(z_2),
\end{equation} the potential corresponding to $d\omega\in\Gamma(\Lambda^n(T^*X))=\Gamma(\Lambda^2(T^*X)),$ where
\begin{equation}
\omega=\eta\wedge\star\zeta=\left(\int\frac{\Sigma\Delta\rho+\Gamma\kappa}{\Gamma\Sigma\Phi(z_1,z_2)}dz_1+\Psi(z_2)\right)dz_2\equiv\left(\int\frac{\Sigma\Delta\rho+\Gamma\kappa}{\Gamma\Sigma\Phi(z_1,z_2)}dz_1+\Psi(z_2)\right)\wedge dz_2
\end{equation} Recall that $\eta\wedge\star\zeta=\langle\eta,\zeta\rangle\varpi$ for $\eta$ a $k$-form, $\zeta$ a $(k+1)$-form and $\varpi$ a $(n-1)$-form for $\star\zeta$ a $(n-k-1)$-form. As such, $n=2,$ and for $k=0,$ $\zeta$ and $\star\zeta$ are both $1$-forms. Therefore, we let $\omega=f\wedge dz_2=\eta\wedge\star\zeta,$ $\star\zeta=dz_2.$ The Hodge dual of $\zeta$ is $\star\zeta=-i\zeta$ or $\zeta=i\star\zeta$ and $\zeta=idz_2.$ To invoke the adjointness property $(d\eta,\zeta)=(\eta,\delta\zeta),$ we first calculate the corresponding terms. Note that 
\begin{equation*}
f(z_1,z_2)=\eta:=\int\frac{\Sigma\Delta\rho+\Gamma\kappa}{\Gamma\Sigma\Phi(z_1,z_2)}dz_1+\Psi(z_2),
\end{equation*} such that the exterior derivative of $\eta$ is
\begin{equation*}
\begin{split}
&d\eta=\frac{\Sigma\Delta\rho+\Gamma\kappa}{\Gamma\Sigma\Phi(z_1,z_2)}dz_1+\left(\frac{\partial\Psi}{\partial z_2}+\int\frac{\partial}{\partial z_2}\frac{\Sigma\Delta\rho+\Gamma\kappa}{\Gamma\Sigma\Phi(z_1,z_2)}dz_1\right)dz_2\\&=\frac{\Sigma\Delta\rho+\Gamma(z_1)\kappa}{\Gamma(z_1)\Sigma\Phi(z_1,z_2)}dz_1+\left(\frac{\partial\Psi}{\partial z_2}+\int-\frac{1}{\Gamma\Phi^2(z_1,z_2)}\frac{\partial^2\rho}{\partial z_1^2}\frac{\partial \Phi}{\partial z_2}-\frac{\kappa}{\Sigma\Phi^2(z_1,z_2)}\frac{\partial\Phi}{\partial z_2}dz_1\right)dz_2
\end{split}
\end{equation*} for $\Gamma=\Gamma(z_1), \Delta\rho=\Delta log\Gamma(z_1)$ solely dependent on the first local coordinate $z_1,$ and $\Sigma,\kappa\in\mathbb{R}$ constants. Likewise, for $\zeta=idz_2,$ $\delta\zeta=\delta(idz_2)=i\delta(dz_2)\equiv-\star^{-1}d(dz_2)=0.$ Noting that the $L^2$-norm for differential $k$-forms is $(\eta,\zeta)=\int_X\eta\wedge\star\zeta,$ the inner product
\begin{equation*}
\begin{split} 
&\left(\left(\frac{\Sigma\Delta\rho+\Gamma(z_1)\kappa}{\Gamma(z_1)\Sigma\Phi(z_1,z_2)}dz_1+\left(\frac{\partial\Psi}{\partial z_2}-\int\frac{1}{\Gamma\Phi^2(z_1,z_2)}\frac{\partial^2\rho}{\partial z_1^2}\frac{\partial \Phi}{\partial z_2}+\frac{\kappa}{\Sigma\Phi^2(z_1,z_2)}\frac{\partial\Phi}{\partial z_2}dz_1\right)dz_2\right),idz_2\right)\\&=0
\end{split}
\end{equation*} becomes
\begin{equation*}
\begin{split}
&\int_X\biggl(\frac{\Sigma\Delta\rho+\Gamma(z_1)\kappa}{\Gamma(z_1)\Sigma\Phi(z_1,z_2)}dz_1+\biggl(\frac{\partial\Psi}{\partial z_2}-\int\frac{1}{\Gamma\Phi^2(z_1,z_2)}\frac{\partial^2\rho}{\partial z_1^2}\frac{\partial \Phi}{\partial z_2}+\frac{\kappa}{\Sigma\Phi^2}\frac{\partial\Phi}{\partial z_2}dz_1\biggr)dz_2\biggr)\wedge\star\left(idz_2\right)=0,
\end{split}
\end{equation*} or for $\star(idz_2)=dz_2$
\begin{equation*}
\begin{split}
&\int_X\left(\frac{\Sigma\Delta\rho+\Gamma(z_1)\kappa}{\Gamma(z_1)\Sigma\Phi(z_1,z_2)}dz_1+\left(\frac{\partial\Psi}{\partial z_2}-\int\frac{1}{\Gamma\Phi^2(z_1,z_2)}\frac{\partial^2\rho}{\partial z_1^2}\frac{\partial \Phi}{\partial z_2}+\frac{\kappa}{\Sigma\Phi^2}\frac{\partial\Phi}{\partial z_2}dz_1\right)dz_2\right)\wedge dz_2\\&=\int_X\frac{\Sigma\Delta\rho+\Gamma(z_1)\kappa}{\Gamma(z_1)\Sigma\Phi(z_1,z_2)}dz_1\wedge dz_2=0,
\end{split}
\end{equation*} which induces a tautology, as asserted. A weaker condition can be imposed by simply observing that if $\Delta\rho=\frac{\partial^2\rho}{\partial z_1^2}>0$ for $\rho=\rho(z_1),$ then  
\begin{equation}
\begin{split}
&\int_Xe^{-log\Gamma}\Delta\rho+\frac{\kappa}{\Sigma}dS=\int_Xe^{-log\Gamma}\Delta\rho+\frac{\kappa}{\Sigma}\left|dz_1\wedge dz_2\right|\\&=\int_X\left|e^{-log\Gamma}\Delta\rho+\frac{\kappa}{\Sigma}dz_1\wedge dz_2\right|\ge \left|\int_Xe^{-log\Gamma}\Delta\rho+\frac{\kappa}{\Sigma}dz_1\wedge dz_2\right|=0.
\end{split}
\end{equation}
Therefore, $g=0$ means that the Riemann surface $X$ is simply connected, corresponding to the first singular homology group being
\begin{equation}
H_k(X;\mathbb{C})=\begin{cases} \mathbb{C} & \text{ if } k=0,2,\\ 0 & \text{ otherwise. }\end{cases}
\end{equation} We obtain the following theorem.
\begin{theorem}
The singular homology groups assume the form \begin{equation}
H_k(X;\mathbb{C})=\begin{cases} \mathbb{C} & \text{ if } k=0,2,\\ 0 & \text{ otherwise } \end{cases}
\end{equation} if and only if the condition for $g=0,$ 
\begin{equation} \int_X\frac{\Sigma\Delta\rho+\Gamma(z_1)\kappa}{\Gamma(z_1)\Sigma\Phi(z_1,z_2)}dz_1\wedge dz_2=0, \end{equation}
in isothermal coordinates $(z_1,z_2)$ is satisfied.
\end{theorem} 

\section{Projective Geometry and Compactness\label{sec:Projective Geometry and Compactness}}
Furthermore, the Riemann surface $X$ can be realized as the universal cover by the quotient of a free, proper action of a discrete group, say $\Gamma.$ In particular $X=\tilde X/\Gamma$ for $\Gamma$ the group of deck transformations $Aut_{\rho}(X)$ acting transitively on the space $X$ under the universal covering map $\rho:\tilde X\to X$ (the identity). We have now shown that if $\phi(z_1)=\lambda e^{z_1\frac{d\beta}{dz_2}-i\alpha(z_1)}$ for $\lambda\in \mathbb{C}$ (condition of biholomorphicity) and $\int_X\frac{\Sigma\Delta\rho+\Gamma(z_1)\kappa}{\Gamma(z_1)\Sigma\Phi(z_1,z_2)}dz_1\wedge dz_2=0$ for $g^2=-{\bar g}^2$ (condition of genus zero), then, by the uniformization theorem, the universal cover $\tilde X$ (and by extension the Riemann surface $X$ for $\rho$ the identity) is necessarily holomorphically isomorphic to the complex projective line, that is $X\cong\mathbb{CP}^1.$ Thus, the Riemann surface $X$ and the Riemann sphere $\mathbb{CP}^1$ can be identified with each other. By the definition of $X,$ $\hat g(\eta)=1/{\eta}^n$ and $f(\eta)=2\pi i\frac{\eta^*\phi(\eta)}{{\eta}^n}.$ As such, $\hat g$ does not have compact support on a compact subset $X_j\subset X$ (unless $\eta$ has a dense set of poles in $X_j,$ of which there are only finitely many for otherwise $\eta$ would be constant, contrary to hypothesis), and by extension $f$ is never of compact support for the biholomorphic map $\eta$, defined locally on $X$ satisfying the Wirtinger derivative condition. To guarantee that $f$ never has compact support, we impose the condition that $|\eta(z_1,z_2)|=|\phi(z_1)|=|\bar g|\equiv|g|<M$ for all real $M>0$ sufficiently large (where $\phi(z_1)$ was defined in Section \ref{sec:genera}) and the globally defined chart $\phi(z_1,z_2)=z_1+iz_2\in\mathbb{C}.$ In particular, $sup_{\phi:=(z_1,z_2)}|g|<+\infty.$ It then follows that $\iint\limits_{X_j}\eta^n\hat g(\eta)d\eta\wedge d\bar \eta=\iint\limits_{U_j}z^n\hat g(z)dz\wedge d\bar z=\sigma_j(\chi)=\int_{\partial X_j}e^{-i\pi\langle\xi,\chi\rangle}\eta^*\omega\ne 0;$ that is, $\sigma_j(\chi)$ is not of compact support unless $\omega\equiv 0$, meaning that $\int_{\partial U_j}\omega=0$ by Theorem \ref{theorem: compact support}, as was to be shown. Let $\Omega^{(1,0)}_X(logD)\ni\omega=\frac{f'(z)}{f(z)}dz$ be a logarithmic $(1,0)$-form for $D\in X$ a principal divisor, $\omega$ a holomorphic $(1,0)$-form on $X-D,$ $\omega$ singular on $D,$ and $f$ analytic on $\overline{I(\partial X_j)}$ except for poles in $I(\partial X_j)$ at finitely many points $b_1,...,b_n. $ Moreover, suppose $f$ has zeros at finitely many points $a_1,...,a_m$ in $I(\partial X_j)$ but none on $\partial X_j$ itself; that is, $f:X\to\mathbb{CP}^1$ is a non-constant meromorphic function. Then since $\partial X_j$ is a closed rectifiable Jordan curve, by the argument principal, \begin{equation*}\frac{1}{2\pi i}\int_{\partial X_j}\frac{f'(z)}{f(z)}dz:=\frac{1}{2\pi i}\int_{\partial X_j}\omega=\sum_{k=1}^m\alpha_k-\sum_{k=1}^n\beta_k=0, \end{equation*}where $\alpha_k$ is the order of $a_k$ and $\beta_k$ the order of $b_k.$ That is, on a compact Riemann surface $X$ of genus $g=0,$ isomorphic to $\mathbb{CP}^1,$ every non-constant meromorphic function $f:X\to\mathbb{CP}^1$ has as many zeros as poles, where each is counted according to multiplicities. Therefore, we have proven the following strong conditions.
\begin{theorem}
If $\int_Xe^{-log\Gamma}\Delta\rho dS=0$ for a Riemann surface $X$ of genus $g=0,$ and if $X\cong\mathbb{CP}^1$ for $X$ given parametrically in terms of local coordinates $(z_i)$ by the diffeomorphism $\eta:X\to\mathbb{CP}^1$ with $\eta(z_1,z_2)=-\bar g(\eta(z_1,z_2))e^{i\pi\langle\chi,\xi\rangle}$ and $\sup_{\phi=(z_1,z_2)}|g|<+\infty,$ then $\int_{\partial U_j}\omega=0$ and $\int_{\partial X_j}\eta^*\omega=0$ necessarily for $\omega\in\mathcal{M}\Omega^{(1,0)}(\mathbb{CP}^1)$ and $\eta^*\omega\in\mathcal{M}\Omega^{(1,0)}(X)$ meromorphic $(1,0)$-forms. 
\end{theorem}
\begin{theorem}\label{theorem: zeros and poles}
For a compact Riemann surface $X$ of genus $g=0,$ isomorphic to $\mathbb{CP}^1,$ every non-constant meromorphic function $f:X\to\mathbb{CP}^1$ has as many zeros as poles, where each is counted according to multiplicities.
\end{theorem}

In performing these calculations and applying the uniformization theorem, we have assumed that the Riemann surface is strictly compact. Recall that the Riemannian metric can be given locally as 
\begin{equation*}
\begin{split}
&g=ds^2=\left(1+\left(\frac{d\phi}{dz_1}\right)^2+\phi^2(z_1)\left(\frac{d\alpha}{dz_1}\right)^2\right)dz_1^2+2\left(i+\phi^2(z_1)\frac{d\beta}{dz_2}\frac{d\alpha}{dz_1}\right)dz_1dz_2\\&+\left(\phi^2(z_1)\left(\frac{d\beta}{dz_2}\right)^2-1\right)dz_2^2.
\end{split}
\end{equation*} The Riemann surface $X$ is a complex, conformal manifold equipped with an equivalence class of Riemannian metrics, for which two metrics $g$ and $h$ are identified if and only if $h=\lambda^2 g$ for $\lambda$ a real-valued smooth function. In particular, two metrics $g$ and $h$ on the Riemann surface are equivalent if and only if 
\begin{equation}
\begin{split}
&h=\lambda^2\biggl[\left(1+\left(\frac{d\phi}{dz_1}\right)^2+\phi^2(z_1)\left(\frac{d\alpha}{dz_1}\right)^2\right)dz_1^2+2\left(i+\phi^2(z_1)\frac{d\beta}{dz_2}\frac{d\alpha}{dz_1}\right)dz_1dz_2\\&+\left(\phi^2(z_1)\left(\frac{d\beta}{dz_2}\right)^2-1\right)dz_2^2\biggr].
\end{split}
\end{equation} An equivalence of such metrics is the conformal class. The standard metric $g=ds^2$  is the restriction of the Euclidean metric to the Riemann surface $X$ with $w=z_1+iz_2$ for the chart $\phi:X\to\mathbb{C}.$ By invoking isothermal coordinates $(z_1,z_2)$ above, we showed that such a conformal metric is in fact conformally flat. The conformal class of $g$ denoted $[g]=\{\lambda^2g|\lambda>0\},$ which gives a realization of $X$ as a complex manifold, is the collection of such representatives. To prove that $X$ is indeed a compact, complex manifold, we use the following well known theorem of Behnke and Stein [1948] that is stated without proof \cite[Pg. 81]{Forster},
\begin{theorem}\label{theorem: stein manifold}
Let $X$ be a connected non-compact Riemann surface. Then $X$ is a Stein manifold.
\end{theorem} In this regard, $X$ is a Stein manifold if it is holomorphically convex whereby for a compact set $K\subset X,$ the convex-hull $\bar K:=\{z\in X:|f(z)|\le sup_{w\in K}|f(w)|, \forall f\in\mathcal{O}(X)\}$ is a compact subset $\bar K\subset X.$ For the diffeomorphism $\eta(z_1,z_2)=\Phi(z_1)\Theta(z_2)e^{i\pi\langle\chi,\xi\rangle}=\hat\Phi(z_1)\Theta(z_2),$ let $p(z_1,z_2)=\eta(z_1,z_2)-\hat\Phi(z_1)\Theta(z_2)=0.$ Then for $z_0\in X\setminus D$ where $D$ is the divisor on $X$ for which $\eta:X\to\mathbb{CP}^1$ is singular, we homogenize the polynomial as follows: 
\begin{equation}\label{equation: homogenization}
{}^hp(z_0,z_1,z_2)=z_0^{deg(p)}p\left(\frac{z_1}{z_0},\frac{z_2}{z_0}\right) \text{ or } {}^hp(z_0,z_1,z_2)=z_0^{deg(p)}\left[\eta\left(\frac{z_1}{z_0},\frac{z_2}{z_0}\right)-\hat\Phi\left(\frac{z_1}{z_0}\right)\Theta\left(\frac{z_2}{z_0}\right)\right].
\end{equation} By Chow's theorem, complex projective varieties are automatically algebraic as they are defined by the vanishing of homogenous polynomial equations. In particular, any compact Riemann surface is a projective variety, i.e. it can be given by polynomial equations inside a projective space as in the case of ${}^hp(z_0,z_1,z_2)=0,$ which is a projective algebraic curve. 
Hence, given that every affine algebraic curve of vanishing polynomial $p(z_1,z_2)=0$ may be completed into the projective curve of equation ${}^hp(z_0,z_1,z_2)=0,$ then such a completion Eq. \ref{equation: homogenization} implies that the Riemann surface is a algebraic curve (algebraic variety of dimension one), for $z_i=z_i(x_1,x_2).$ We hereby invoke the theorem of Griffiths and Harris \cite[Pg. 215]{Griffiths} that every compact Riemann surface is an algebraic curve. Therefore, the Riemann surface $X$ must be compact, yielding the following theorem.
\begin{theorem}
The Riemann surface $X$ given by the vanishing polynomial equation
\begin{equation}
z_0^{deg(p)}\left[\eta\left(\frac{z_1}{z_0},\frac{z_2}{z_0}\right)-\hat\Phi\left(\frac{z_1}{z_0}\right)\Theta\left(\frac{z_2}{z_0}\right)\right]=0
\end{equation} is an algebraic variety of dimension one such that $X$ is necessarily compact by Griffiths and Harris.
\end{theorem}
For genus $g=0,$ $X$ is biholomorphic to ${\mathbb{CP}}^1.$ Therefore $X$ is a simply connected, compact algebraic variety of dimension one. This completes the proof of the isomorphism $X\cong\mathbb{CP}^1$ for $\tilde X=\mathbb{CP}^1$ the universal cover of $X$ and $\Gamma$ a discrete group acting on $X.$ Hence, $X$ and $\mathbb{CP}^1,$ an elliptic geometry of positive constant curvature, can be identified with one another. It follows that $X$ also inherits an elliptic geometry. Note that we associate to the manifold $X$ its universal cover $\tilde X=\mathbb{CP}^1,$ expressing the original $X$ as the quotient of $\tilde X$ by the group of deck transformations $Aut(\rho)$ for $\rho:\tilde X\to X$ a universal covering map, for which an automorphism of a cover $\rho$ is a homeomorphism $f:\tilde X\to\tilde X$ such that $\rho\circ f=\rho.$ Such a deck transformation permutes the elements of each fiber of $\rho.$  If this action is transitive on some fiber, then it is transitive on all fibers such that the cover $\rho$ is regular. Every universal cover is regular such that the group of deck transformations $Aut(X)=Aut_X(\rho)$ (note, we omit the explicit reference to $X$ when it is clear from context) is isomorphic to the first homotopy group, i.e. $Aut(X)\cong\pi_1(X),$ which is trivial since $X$ is simply connected. Therefore, the topological space $X$ can be expressed as $X:=\tilde X/Aut(\rho)=\tilde X.$ If $\Gamma\subset Aut_X(\rho)$ then $X=\tilde X.$ For an elliptic Riemann surface $X,$ by the uniformization theorem, the universal cover of $X$ has to be (identified with) the complex projective line.

\section{The M\"{o}bius Group and Quotient Topology \label{sec:Mobius}}

The Riemann surface $X$ is realized as a universal cover by the $Aut(\rho)$-action of deck transformations, an orbit $\tilde X/Aut(\rho)$ for $\rho:\tilde X\to X$ a universal covering map. The automorphisms on $\tilde X$ therefore act as automorphisms of $\mathbb{CP}^1=\hat{\mathbb{C}}$ as a complex manifold (i.e. a complex Lie group), whereby $Aut(\hat{\mathbb{C}})=\{\text{meromorphic bijections } f:\hat{\mathbb{C}}\to\hat{\mathbb{C}} \}$ is the M\"{o}bius group. For every invertible $2$-by-$2$ matrix $\frak h=\left(\begin{array}{cc}a&b\\c&d\end{array}\right)$ we can associate a M\"{o}bius transformation $f(z)=\frac{az+b}{cz+d}$ such that $det(\frak h)\ne 0.$ Let $\hat{\pi}: GL(2,\mathbb{C})\to Aut(\hat{\mathbb{C}})$ be a group homomorphism from the general linear group $GL(2,\mathbb{C})$ to the group of deck transformations on $\hat{\mathbb{C}},$ sending $\frak h$ to the transformation $f.$ Note that $\hat{\pi}$ is not injective because all nonzero scalar multiples of a given matrix $\frak h$ are taken to the same automorphism. The kernel of the map is the subgroup of $GL(2,\mathbb{C})$ consisting of all nonzero scalar multiples of the identity matrix $(\mathbb{C}\setminus\{0\})I=\mathbb{C}^*I=\biggl\{\lambda\left(\begin{array}{cc}1&0\\0&1\end{array}\right):\lambda\in\mathbb{C},\lambda\ne 0\biggr\}.$ By the first isomorphism theorem, there is an isomorphism $GL(2,\mathbb{C})/(\mathbb{C}^*I)\cong Aut(\hat{\mathbb{C}}), \biggl\{\lambda\left(\begin{array}{cc}a&b\\c&d\end{array}\right)\biggr\}\mapsto f(z)=\frac{az+b}{cz+d}.$ The quotient of $GL(2,\mathbb{C})$ by the nonzero scalar multiples of the identity matrix is the projective general linear group of $2$-by-$2$ complex matrices, $PGL(2,\mathbb{C})=GL(2,\mathbb{C})/(\mathbb{C}^*I)\cong Aut(\hat{\mathbb{C}})\cong PSL(2,\mathbb{C}).$ The M\"{o}bius group $Aut(\hat{\mathbb{C}})$ can be given the structure of a complex Lie group such that composition and inversion are biholomorphic and $Lie(Aut(\hat{\mathbb{C}}))\cong aut(\hat{\mathbb{C}})$ is a Lie algebra. 

Then $X\cong\mathbb{CP}^1$ is a Lie group for it can be given the structure of a complex manifold in such a way that composition and inversion are holomorphic maps. Note that $\hat{\pi}\circ\left(\begin{array}{cc}\alpha&\beta\\ \gamma&\delta\end{array}\right)\in Aut(\mathbb{\hat C})\cong PGL(2,\mathbb{C}).$ Consequently, for $\eta^*\omega\in\mathcal{M}\Omega^{(1,0)}(X)$ a meromorphic $(1,0)$-form, $p:=\eta$ (chosen for standard notational reasons) a diffeomorphism, $\gamma_j=\partial U_j$ a rectifiable Jordan curve in $\mathbb{CP}^1$ (i.e. a continuous map from the unit interval $[0,1]$ into $\mathbb{CP}^1$) and $c\in X$ a point lying over $\gamma_j(0)=\gamma_j(1)$ (i.e. for $p:X\to\mathbb{CP}^1$ a cover, $c$ is in the fiber over $\gamma_j(0)\in\mathbb{CP}^1,$ and $p(c)=\gamma_j(0)$), then there exists a unique path $\Gamma_j=\partial X_j$ lying over $\gamma_j$ (for $p\circ\Gamma_j=\gamma_j$) such that $\Gamma_j(0)=\Gamma_j(1),$ with trivial monodromy action, where $X$ has no ramification points (i.e. $X$ is not ramified). In particular, the degree of the cover $p:X\to\mathbb{CP}^1$ (that is, the cardinality of any fiber of $p$) is equal to the index $[\pi_1(\mathbb{CP}^1,\gamma_j(0)):p_{\sharp}(\pi_1(X,\Gamma_j(0)))]=1$ for $p_{\sharp}:\pi_1(X,\Gamma_j(0))\to\pi_1(\mathbb{CP}^1,\gamma_j(0))$. For the genus defined as half of the first Betti number, i.e., half of the ${\mathbb{C}}$-dimension of the first singular homology group $H_1(X;\mathbb{C})$  with complex coefficients, we obtain the following result.
\begin{theorem}
The Riemann surface $X$ is simply connected, such that the first singular homology group  $H_1(X;\mathbb{C})$ is trivial and the genus (half of the first Betti number) is, thus, necessarily zero.
\end{theorem}
The curve $\Gamma_j$ is the lift of $\gamma_j$ by $p.$ Similarly, since $\rho:\tilde X\to X,$ the identity, induces the isomorphism $X\cong\mathbb{CP}^1,$ then $Aut(X)\cong Aut(\hat{\mathbb{C}})=\{\text{meromorphic bijections }f| f(z)=\frac{az+b}{cz+d},\text{ }ad-bc\ne 0\}.$ It follows that $\Pi\in \frak{aut}(X)\cong \frak{aut}(\mathbb{CP}^1),$ the Lie algebra of automorphisms of the complex projective line $\mathbb{CP}^1.$ Let $\xi,\chi\in X_j$ and let $\Pi\in \frak{aut}(X)$ be an action on the $\frak{aut}(X)$-space $X,$ the Lie subalgebra associated with $Aut_{\rho}(X),$ such that $Lie(Aut_{\rho}(X))\cong \frak{aut}_{\rho}(X)$ for $Aut_{\rho}(\hat{\mathbb{C}}):=\{f|f(z)=\frac{az+b}{cz+d},\text{ }ad-bc\ne 0\}\cong PGL(2,\mathbb{C})\cong PSL(2,\mathbb{C}).$ Assuming the action of $\Pi$ on $X$ is transitive, then there exists a $\Pi\in \frak{aut}(X)$ such that $\Pi_{\frak{aut}(X)}\xi=\chi.$ Then, it follows that 
\begin{equation*}
\begin{split}
&\sigma_j(\chi)=\int_{\partial X_j}e^{-i\pi\langle\xi,\chi\rangle}\eta^*\omega=\int_{\partial X_j}e^{-i\pi\langle\xi,\Pi\xi\rangle}\eta^*\omega=\int_{\partial X_j}e^{-i\pi\langle\xi,log(e^{\Pi\xi})\rangle}\eta^*\omega\\&=\int_{\Gamma_j,\bar\Pi\in Aut(X)}e^{-i\pi\langle\xi,\xi log\bar\Pi\rangle}\eta^*\omega=\int_{\Gamma_j,\left(\begin{array}{cc}\alpha&\beta\\ \gamma&\delta\end{array}\right)\in GL(2,\mathbb{C})}e^{-i\pi\biggl\langle\xi,\xi log\left(\hat{\pi}\circ\left(\begin{array}{cc}\alpha&\beta\\ \gamma&\delta\end{array}\right)\right)\biggr\rangle}\eta^*\omega\\&=\int_{\Gamma_j,\frak h\in GL(2,\mathbb{C})}e^{-i\pi|\xi|^2log(\hat{\pi}\circ\frak h)}\eta^*\omega=\int_{\Gamma_j,\frak h\in GL(2,\mathbb{C})}\left(\hat{\pi}\circ\frak h\right)^{-i\pi|\xi|^2}\eta^*\omega,
\end{split}
\end{equation*} invoking the identification of $X$ with $\mathbb{CP}^1.$ 
Let $M=\tilde X$ be a topological space, namely the universal covering space of $X,$ and let $\Gamma\subset Aut(M),$ a subgroup of the group of automorphisms on $M$. Thus, let $\Gamma:=\{g(z)=\tilde z|\tilde z\in Aut(M)\}.$ We define the equivalence relation $z\sim\tilde z$ for $z\in M$ if and only if $z-\tilde z\in\Gamma.$ We denote by $[z]$ the equivalence class represented by $z.$ Then from the natural projection $\pi: M\to M/\Gamma=M/\sim,\text{ } z\mapsto [z],$ we obtain the quotient space $M/\sim$ or $M/\Gamma,$ and we define a quotient topology on $M/\Gamma.$ Namely the subset $\hat U\subset M/\Gamma$ is open if and only if $\pi^{-1}(\hat U)$ is open in $X.$ Let $\nu=\{[U]=U/\sim|U \text{ is open in }M\text{ such that }g(U)\cap U=\emptyset\text{ for }g\ne Id, g\in\Gamma\}$ \cite[Pg. 1-5]{Ionel}. Then $\nu$ forms a basis for the topology of $M/\Gamma.$ Observe that $\pi: M\to M/\Gamma$ is a covering quotient map such that for any $p\in M/\Gamma,$ $p$ has a neighborhood $[U_p]\subset\nu$ and 
\begin{equation}
\pi^{-1}([U_p])=\bigcup_{g\in\Gamma}g(U_p) \text{ and }g(U_p)\bigcap g'(U_p)\ne\emptyset \text{ if and only if }g=g'.\end{equation} In particular $\pi|_{g(U_p)}:g(U_p)\to[U_p]$ is a homeomorphism. By interpreting $(\pi|_{g(U_p)})^{-1}$ as a coordinate map, $M/\Gamma=\tilde X/\Gamma=\tilde X$ is a complex manifold where such a topology is equivalent to the previously constructed topology of $\mathbb{CP}^1.$
Note that, in the above, for $\frak{aut}(X)$ a vector space (a Lie algebra endowed with a Lie bracket commutator $[X,Y]=XY-YX$) and $X$ a topological space, the left group action $\phi$ of $\frak{aut}(X)$ on $X$ is a function $\phi:\frak{aut}(X)\times X\to X: (\Pi,x)\mapsto\phi(\Pi,x)$ that satisfies the identity, compatibility axioms (where we denote $\phi(\Pi,x)$ as $\Pi\cdot x$ for $\Pi\in \frak{aut}(X)$). The action of $\Pi$ on $X$ is transitive because $X$ is non-empty, and for $\xi,\chi\in X$ there exists a $\Pi$ in $\frak{aut}(X)$ such that $\chi=\Pi\cdot \xi$ for $X$ an $\frak{aut}(X)$-space equipped with an action of $\Pi$ on $X.$ Therefore, $\frak{aut}(X)$ automatically acts by automorphisms $\Pi$ on the set (topological space). If $X$ in addition belongs to some category, then the elements of $\frak{aut}(X)$ are structure preserving. Thus, $X$ is a homogeneous $\frak{aut}(X)$-space on which $\frak{aut}(X)$ acts transitively. If $X$ is an object of the category $C,$ then the structure of a $\frak{aut}(X)$-space is a homomorphism $\rho:\frak{aut}_C(X)\to Aut_C(X)$ into the group automorphisms of the object $X=\tilde X/\Gamma=\tilde X/\pi_1(X)=\tilde X/Aut_C(X)\equiv\tilde X$ in the category $C.$ The pair $(X,\rho)$ defines the homogenous space provided $\rho(\frak{aut}_C(X))$ is a transitive group of symmetries of the underlying set $X$ with $\rho: \frak{aut}_C(X)\to Aut_C(X)$ evidently the exponential map. For the Riemann surface $X=\tilde X\cong \mathbb{CP}^1,$ the Burnside's lemma gives the cardinality
\begin{equation}
|\tilde X/\Gamma|=\frac{1}{|\Gamma|}\sum_{g\in\Gamma}\left|{\tilde X}^g\right|
\end{equation} where ${\tilde X}^g$ is the set of points fixed by $\Gamma.$

\section{Cohomology Theory\label{sec:cohomology}}
Consider the de Rham complex, i.e. the cochain complex of differential forms on the Riemann surface $X$ with the exterior derivative $d^p$ as the coboundary operator $d^p:\Omega^p(X)\to\Omega^{p+1}(X),$
\begin{equation}
0 \to \Omega^0(X) \xrightarrow{d^0} \Omega^1(X) \xrightarrow{d^1}...\xrightarrow{d^{p-1}}\Omega^p(X) \xrightarrow{d^p}...
\end{equation}
for $\Omega^p(X):=\Omega^{(p,0)}(X).$
Closed forms on $X$ are classified by requiring that two closed forms $\alpha,\beta\in\Omega^p(X)$ are cohomologous if they differ by an exact form, i.e. $\alpha-\beta$ is an exact form. Such a classification gives rise to an equivalence class on the space of closed forms in $\Omega^p(X),$ for the $p$-th de Rham cohomology $H^p_{dR}(X)$ the set of equivalence classes. For later application, we begin by extending de Rham's theorem to the case in which the coefficient field of cohomology is $\mathbb{C}$ and the manifold under consideration is the Riemann surface $X$ analyzed above. Consider the map $I:H^p_{dR}(X)\to H^p(X;\mathbb{C})$ defined in the following manner: For any $[\omega]\in H_{dR}^p(X),$ by assumption, let $I(\omega)$ be the element of $Hom(H_p(X;\mathbb{C}),\mathbb{C})\cong H^p(X;\mathbb{C})$ that acts as
\begin{equation}
H_p(X;\mathbb{C})\ni[c]\mapsto\int_c\omega\in Hom(H_p(X;\mathbb{C}),\mathbb{C})\cong H^p(X;\mathbb{C})
\end{equation} for $I(\omega):[\omega]\mapsto \int_c\omega$ or $I(\omega):H^p_{dR}(X)\to H^p(X;\mathbb{C}),$ where $c$ is a $p$-cycle representing the homology class $[c]\in H_p(X;\mathbb{C}).$ The theorem of de Rham asserts that such a map is in fact an isomorphism between de Rham cohomology and singular cohomology. To construct singular cohomology, consider the set of all possible $n$-simplices $\sigma_n(\Delta^n)$ on a topological space for the continuous mapping $\sigma_n:\Delta^n\to X.$ This may be used as the basis of a free abelian group such that each $\sigma_n(\Delta^n)$ is a generator of the group. Note that the set of generators is usually infinite as there are many ways of mapping any one simplex into the topological space. The $n$-simplex $\Delta^n$ is the convex hull of $n+1$ vertices. More precisely for $n+1$ points $u_0,...,u_n\in{\mathbb{R}}^n$ affinely independent we can define the $n$-simplex by
\begin{equation}
\Delta^n:=\left\{\theta_0u_0+...+\theta_nu_n\biggr\rvert\sum_{i=0}^n\theta_i=1 \text{ where }\theta_i\ge 0,\forall i\right\}.
\end{equation} The free abelian group generated by this basis is denoted $C^n(X)$ for $\sum_in_i\sigma_i,$ with $n_i\in\mathbb{Z},$ an element of $C^n(X).$ The coboundary operator $\partial^n:C^n(X)\to C^{n+1}(X)$ is defined to act on singular cochains. The coboundary operator together with the free abelian groups $C^n$ form a cochain complex $C^*$, namely the singular cochain complex. Hence, we define the $n$-th cohomology group as the quotient $H^n(X)=ker(\partial^n)/im(\partial^{n-1}):=Z^n(X)/B^n(X)$ for the coboundary operator satisfying $\partial^n\circ\partial^{n-1}=0_{n-1,n+1}.$ Let $C^*$ and $\Omega^*$ be the singular cochain and de Rham cochain complexes, respectively. Let $f$ be a map between the two cochain complexes $\Omega^*:=(\Omega^{\bullet}(X),d^{\bullet})$ and $C^*:=(C^{\bullet}(X),\partial^{\bullet})$ whereby $f_n:\Omega^n(X)\to C^n(X)$ is a sequence of homomorphisms, for each $n,$ that commutes with the coboundary operators on the two cochain complexes $\partial^n\circ f_n=f_{n+1}\circ d^n.$ Such a map sends cocycles to cocycles and coboundaries to coboundaries, and thus descends to a map on cohomology $(f_{\bullet})^{*}:H^{\bullet}(\Omega^{\bullet}(X),d^{\bullet})\to H^{\bullet}(C^{\bullet}(X),\partial^{\bullet}).$ Therefore with $f_n:\Omega^n(X)\to C^n(X),$ the following diagram commutes:
\begin{equation}
\begin{tikzcd}
0 \arrow{r}
& \Omega^0(X) \arrow{d}{f_0} \arrow{r}{d^0} & \Omega^1(X)\arrow{d}{f_1}\arrow{r}{d^1}&\dots
\arrow{r}{d^{p-2}}&\Omega^{p-1}(X)\arrow{d}{f_{p-1}} \arrow{r}{d^{p-1}}&\Omega^p(X)\arrow{d}{f_p}\arrow{r}{d^p}&\dots\\
0 \arrow{r}
& C^0(X) \arrow{r}{\partial^0} & C^1(X)\arrow{r}{\partial^1}& \dots
\arrow{r}{\partial^{p-2}}&C^{p-1}(X)  \arrow{r}{\partial^{p-1}}&C^p(X)\arrow{r}{\partial^p}&\dots.
\end{tikzcd}.
\end{equation}
In the notation of Hodge theory, the coboundary operator $d$ coincides with the Dolbeaut operator $\partial$ (different from the coboundary operator associated with the singular cohain complex) given by $\partial=\pi^{(p+1,q)}\circ d:\Omega^{(p,q)}\to\Omega^{(p+1,q)}$ with $k:=p+q,$ $E^k$ the the total degree space of complex differential forms, and the canonical projection of vector bundles $\pi^{(p,q)}:E^k\to\Omega^{(p,q)}$ for $q=0,$ such that this commutative diagram becomes
\begin{equation}
\begin{tikzcd}
0 \arrow{r}
& \Omega^{(0,0)}(X) \arrow{d}{f_0} \arrow{r}{d^0} & \Omega^{(1,0)}(X)\arrow{d}{f_1}\arrow{r}{d^1}&\dots
\arrow{r}{d^{p-2}}&\Omega^{(p-1,0)}(X)\arrow{d}{f_{p-1}} \arrow{r}{d^{p-1}}&\Omega^{(p,0)}(X)\arrow{d}{f_p}\arrow{r}{d^p}&\dots\\
0 \arrow{r}
& C^0(X) \arrow{r}{\partial^0} & C^1(X)\arrow{r}{\partial^1}& \dots
\arrow{r}{\partial^{p-2}}&C^{p-1}(X)  \arrow{r}{\partial^{p-1}}&C^p(X)\arrow{r}{\partial^p}&\dots.
\end{tikzcd}.
\end{equation} The homomorphisms $f_n:\Omega^{(n,0)}(X)\to C^n(X)$ descend onto cohomology and induce a sequence of homomorphisms $(f_n)^*:H_{dR}^n(X)\to H^n(X;\mathbb{C})$ such that $(d^n)^*$ and $(\partial^n)^*$ are induced coboundary operators associated with the respective cochain complexes $(H_{dR}^{\bullet},(d^{\bullet})^*)$ and $(H^{\bullet}(X;\mathbb{C}),(\partial^{\bullet})^*).$ Therefore, recalling that $X$ is a smooth complex manifold and $H^0(X;\mathbb{C})\cong\mathbb{C}$ with $H_{dR}^0(X)\cong\mathbb{C}$ for $X$ simply connected, we obtain the following commutative diagram:
\begin{equation}
\begin{tikzcd}
0 \arrow{r}
& H_{dR}^0(X) \arrow{d}{(f_0)^*} \arrow{r}{(d^0)^*}  & H_{dR}^1(X)\arrow{d}{(f_1)^*}\arrow{r}{(d^1)^*}&\dots
\arrow{r}{(d^{p-2})^*}&H_{dR}^{p-1}(X)\arrow{d}{(f_{p-1})^*} \arrow{r}{(d^{p-1})^*}&H_{dR}^p(X)\arrow{d}{(f_p)^*}\arrow{r}{(d^p)^*}&\dots\\
0 \arrow{r}
& H^0(X;\mathbb{C}) \arrow{r}{(\partial^0)^*} & H^1(X;\mathbb{C})\arrow{r}{(\partial^1)^*}& \dots
\arrow{r}{(\partial^{p-2})^*}&H^{p-1}(X;\mathbb{C})  \arrow{r}{(\partial^{p-1})^*}&H^p(X;\mathbb{C})\arrow{r}{(\partial^p)^*}&\dots
\end{tikzcd}.
\end{equation} The induced sequence of homomorphisms $(f_n)^*:H_{dR}^n(X)\to H^n(X;\mathbb{C})$ is precisely $I(\omega^n)=\int_{c^n}\omega^n\in H^n(X,\mathbb{C})$ for $c^n$ an $n$-cycle representing the homology class $[c^n]\in H_n(X;\mathbb{C}).$ If $I(\omega^n)=\int_{c^n}\omega^n\ne 0,$ then $H^n(X;\mathbb{C})$ is necessarily nontrivial if and only if the $n$-th de Rham cohomology group $H_{dR}^n(X)$ is nontrivial. In particular, 
\begin{theorem}
For suppose $[\theta^n]\in H_{dR}^n(X)$ then $\int_{c^n}\theta^n\ne 0$ implies that $H^n(X;\mathbb{C})\ne 0$ if and only if the $n$-th de Rham cohomology group $H_{dR}^n(X)$ is nontrivial for $c^n$ an $n$-cycle representing the homology class $[c^n]\in H_n(X;\mathbb{C}).$
\end{theorem} For the case of $n=1,$ let $\theta^1:=e^{-i\pi\langle\xi,\chi\rangle}\eta^*\omega\in\Omega^{(1,0)}(X)$ for $[\theta^1]\in H_{dR}^1(X),$ then $\int_{c^1}\theta^1\ne 0$ means that the first singular cohomology group $H^1(X;\mathbb{C})$ is nontrivial if and only if $H_{dR}^1(X)$ is nontrivial. The condition $\int_{c^1}\theta^1=\int_{\partial X_j}e^{-i\pi\langle\xi,\chi\rangle}\eta^*\omega\ne 0$ follows from $X:=\tilde X/\Gamma=\tilde X/\pi_1(X)=\tilde X\cong\mathbb{CP}^1$ for the Riemann surface $X$ belonging to the category $C.$ This result leads naturally to the following theorem:
\begin{theorem}
If $X$ is a compact Riemann surface belonging to the category $C,$ and if $X\cong\mathbb{CP}^1$ for $\Gamma=Aut_C(X)\equiv\pi_1(X)$ trivial where $X$ is simply connected, then the first singular cohomology group is nontrivial if and only if the first de Rham cohomology group is nontrivial. However, since $X$ is simply connected, the first singular cohomology group $H^1(X;\mathbb{C})$ vanishes, which implies that for $(f_1)^*=\int_{c^1}\theta^1$
\begin{equation*}
ker\left((f_1)^*\right)=\{[\omega]\in H_{dR}^1(X):\text{ } (f_1)^*([\omega])=e_{H^1(X;\mathbb{C})}=0\}=H_{dR}^1(X),
\end{equation*} by the first isomorphism theorem since $H^1(X;\mathbb{C})$ is the trivial group, i.e. $ker\left(\int_{c^1}\theta^1\right)=H_{dR}^1(X)$ where $c^1$ denotes a $1$-cycle in $[c^1].$
\end{theorem} Hurewicz's theorem states that the abelianization of the fundamental group (i.e. the first homotopy group) is isomorphic to the first homology group $H_1(X)\cong\pi_1(X)/[\pi_1(X),\pi_1(X)].$ That is, the canonical abelianization map $h_{*}:\pi_1(X)\to\pi_1(X)/[\pi_1(X),\pi_1(X)]$ is an isomorphism. In this particular case, for $X$ a compact Riemann surface, the first cohomology group vanishes because $X$ is path connected and $\pi_1(X)$ is a perfect group. Assuming that the homology groups with $\mathbb{C}$-coefficients are finitely generated, then this means that $H^n(X;\mathbb{C})\cong H_n(X;\mathbb{C})$ for the dimension of the dual space of a finite-dimensional vector space is the same as the dimension of the vector space, inducing an isomorphism. Consequently, by Hurewicz's theorem, $H_1(X;\mathbb{C})\cong H^1(X;\mathbb{C})\cong \pi_1(X)/[\pi_1(X),\pi_1(X)].$ However, $H^1(X;\mathbb{C})$ vanishes and therefore, for $\pi_1(X)=Aut_C(X)\cong Aut(\hat{\mathbb{C}}),$ $[\pi_1(X),\pi_1(X)]\cong [Aut(\hat{\mathbb{C}}),Aut(\hat{\mathbb{C}})]$ is the normal commutator subgroup. Recall that $Aut(\hat{\mathbb{C}}):=\{f|f(z)=\frac{az+b}{cz+d},\text{ }ad-bc\ne 0\}$ and if $\rho: Aut(\hat{\mathbb{C}})\to GL(2,\mathbb{C})$ is a representation, then for $\frak h=\left(\begin{array}{cc}\alpha&\beta\\ \gamma&\delta\end{array}\right)\in GL(2,\mathbb{C}),$ $\rho^{-1}\circ\left(\begin{array}{cc}\alpha&\beta\\ \gamma&\delta\end{array}\right)\in Aut(\hat{\mathbb{C}}).$ Let $det(\frak h)\ne 0$ for $\frak h\in GL(2,\mathbb{C})$ be an equivalence class on $GL(2,\mathbb{C}),$ then
\begin{equation}
Aut(\hat{\mathbb{C}}):=\left\{\rho^{-1}\circ\left(\begin{array}{cc}\alpha&\beta\\ \gamma&\delta\end{array}\right)/\sim\biggr\rvert\left(\begin{array}{cc}\alpha&\beta\\ \gamma&\delta\end{array}\right)\in GL(2,\mathbb{C})\right\}=\{\text{meromorphic bijections }f:\hat{\mathbb{C}}\to\hat{\mathbb{C}}\}.
\end{equation} Since a M\"{o}bius transformation determines $\frak h$ only up to scalar multiples $\lambda\in\mathbb{C}^*,$ then $Aut(\hat{\mathbb{C}})\cong GL(2,\mathbb{C})/(\mathbb{C}^*I).$ Thus an element of $[Aut(\hat{\mathbb{C}}),Aut(\hat{\mathbb{C}})]$ assumes the form 
\begin{equation*}
\begin{split}
&\left[\rho^{-1}\circ \left(\begin{array}{cc}\alpha_1&\beta_1\\ \gamma_1&\delta_1\end{array}\right),\rho^{-1}\circ\left(\begin{array}{cc}\lambda\alpha_1&\lambda\beta_1\\ \lambda\gamma_1&\lambda\delta_1\end{array}\right)\right]...\left[\rho^{-1}\circ \left(\begin{array}{cc}\alpha_n&\beta_n\\ \gamma_n&\delta_n\end{array}\right),\rho^{-1}\circ\left(\begin{array}{cc}\lambda\alpha_n&\lambda\beta_n\\ \lambda\gamma_n& \lambda\delta_n\end{array}\right)\right]\\&=\prod_{i=1}^n\left[\rho^{-1}\circ \left(\begin{array}{cc}\alpha_i&\beta_i\\ \gamma_i&\delta_i\end{array}\right),\rho^{-1}\circ\left(\begin{array}{cc}\lambda\alpha_i&\lambda\beta_i\\ \lambda\gamma_i&\lambda\delta_i\end{array}\right)\right]
\end{split}
\end{equation*} and 
\begin{equation*}
\begin{split}
&\left[\rho^{-1}\circ \left(\begin{array}{cc}\alpha_i&\beta_i\\ \gamma_i&\delta_i\end{array}\right),\rho^{-1}\circ\left(\begin{array}{cc}\lambda\alpha_i&\lambda\beta_i \\ \lambda\gamma_i&\lambda\delta_i\end{array}\right)\right]\\&= \left(\begin{array}{cc}\alpha_i&\beta_i\\ \gamma_i&\delta_i\end{array}\right)^{-1}\circ\rho  \left(\begin{array}{cc}\lambda\alpha_i&\lambda\beta_i\\ \lambda\gamma_i& \lambda\delta_i\end{array}\right)^{-1}\circ\rho\circ\rho^{-1}\circ  \left(\begin{array}{cc}\alpha_i&\beta_i\\ \gamma_i& \delta_i\end{array}\right)\rho^{-1}\circ  \left(\begin{array}{cc}\lambda\alpha_i&\lambda\beta_i\\ \lambda\gamma_i& \lambda\delta_i\end{array}\right)\\&= \left(\begin{array}{cc}\alpha_i&\beta_i\\ \gamma_i&\delta_i\end{array}\right)^{-1}\circ\rho\left(\begin{array}{cc}\lambda\alpha_i&\lambda\beta_i\\ \lambda\gamma_i& \lambda\delta_i\end{array}\right)^{-1}\left(\begin{array}{cc}\alpha_i&\beta_i\\ \gamma_i&\delta_i\end{array}\right)\rho^{-1}\circ\left(\begin{array}{cc}\lambda\alpha_i&\lambda\beta_i\\ \lambda\gamma_i&\lambda\delta_i\end{array}\right),
\end{split}
\end{equation*} which implies that an element of the normal commutator subgroup is given by 
\begin{equation*}
\begin{split}  
&\prod_{i=1}^n\left(\begin{array}{cc}\alpha_i&\beta_i\\ \gamma_i&\delta_i\end{array}\right)^{-1}\circ\rho\left(\begin{array}{cc}\lambda\alpha_i&\lambda\beta_i\\ \lambda\gamma_i&\lambda\delta_i\end{array}\right)^{-1}\left(\begin{array}{cc}\alpha_i&\beta_i\\ \gamma_i&\delta_i\end{array}\right)\rho^{-1}\circ\left(\begin{array}{cc}\lambda\alpha_i&\lambda\beta_i\\ \lambda\gamma_i&\lambda\delta_i\end{array}\right) \\&=\prod_{i=1}^n\left(\begin{array}{cc}\alpha_i&\beta_i\\ \gamma_i&\delta_i\end{array}\right){\rho}^{-1}\circ\left(\begin{array}{cc}\alpha_i&\beta_i\\ \gamma_i&\delta_i\end{array}\right)
\end{split}
\end{equation*}
for the representation $\rho: Aut(\hat{\mathbb{C}})\to GL(2,\mathbb{C})$ and for
\begin{equation*}
\frak h_i^{+}:=\left(\begin{array}{cc}\alpha_i&\beta_i\\ \gamma_i&\delta_i\end{array}\right),\text{ }\frak h_i^{-}:=\left(\begin{array}{cc}\lambda\alpha_i&\lambda\beta_i\\ \lambda\gamma_i& \lambda\delta_i\end{array}\right)\in GL(2,\mathbb{C})
\end{equation*} with $\lambda\in\mathbb{C}^*.$ It follows that an element of the commutator subgroup has the form
\begin{equation}
\prod_{i=1}^n\left(\frak h_i^{+}\right)^{-1}\circ\rho\left(\frak h_i^{-}\right)^{-1}\frak h_i^{+}\rho^{-1}\circ \frak h_i^{-1}=\prod_{i=1}^n \frak h_i^+{\rho}^{-1}\circ\frak h_i^+\in [\pi_1(X),\pi_1(X)].
\end{equation}
\section{Degree Theory\label{sec: degree theory}} 
For any meromorphic function $f,$ there exists a divisor $D,$ a finite linear combination of points on the Riemann surface $X$ with integer coefficients, defined as \cite[Pg. 116-117]{Griffiths}
\begin{equation*}
(f):=\sum_{z_{\nu}\in R(f)}s_{\nu}z_{\nu}
\end{equation*}
where $R(f)$ denotes the set of all zeros and poles of $f,$ and $s_{\nu}$ is defined as 
\begin{equation*}
s_{\nu}:=\begin{cases} \alpha & \text{if } z_{\nu} \text{ is a zero of order } \alpha, \\ -\alpha & \text{if } z_{\nu} \text{ is a pole of order }\alpha \end{cases}
\end{equation*} for $z_{\nu}\in X.$ The divisor $(f):=\sum_{z_{\nu}\in R(f)}s_{\nu}z_{\nu}$ is equivalent to the integral $\frac{1}{2\pi i}\int_{\partial U_j}z\frac{f'(z)}{f(z)}dz,$ modulo pullback for $\phi: X\to\mathbb{C};$ that is,
\begin{equation*}
(f):=\sum_{z_{\nu}\in R(f)}s_{\nu}z_{\nu}=\frac{1}{2\pi i}\int_{\partial X}\phi^*\left(z\frac{f'(z)}{f(z)}dz\right)=\frac{1}{2\pi i}\int_{\partial X}(z\circ\phi)\frac{f'(z\circ\phi)}{f(z\circ\phi)}d(z\circ\phi).
\end{equation*} Recall, it was shown that the simply connected Riemann surface $X$ is compact. Therefore, invoking the Riemann-Roch theorem for a compact Riemann surface of genus $g$ with canonical divisor $K,$ which states $\ell(D)-\ell(K-D)=deg(D)-g+1,$ it follows that $deg(f)=0$ for any principal divisor $(f):=D$ on $X$ since a meromorphic function has as many zeros as poles (see Theorem \ref{theorem: zeros and poles}). As such,
\begin{equation*}
\ell(D)-\ell(K-D)=deg(D)+1\equiv \frac{1}{2\pi i}\int_{\partial X}\phi^*\left(\frac{f'(z)}{f(z)}dz\right)+1=1
\end{equation*} or
\begin{equation*}
\ell(D)-\ell(K-D)-1=0
\end{equation*} for $\int_{\partial X}\phi^*\left(\frac{f'(z)}{f(z)}dz\right)$ well defined since $X$ is compact and orientable, and $g=0$ follows from the first singular cohomology group being trivial, i.e. $H_1(X;\mathbb{C})=0$. If $X\cong\mathbb{CP}^1$ and $\sup_{\phi=(z_1,z_2)}|g|<+\infty,$ then $\int_{\partial X}\phi^*\left(\frac{f'(z)}{f(z)}dz\right)$ vanishes, and we have the following theorem.
\begin{theorem}
If $X\cong\mathbb{CP}^1$ and $\sup_{\phi=(z_1,z_2)}|g|<+\infty,$ then $X$ is a compact Riemann surface of genus $g=0$ with canonical divisor $K,$ such that $\ell(D)-\ell(K-D)-1=0.$ In fact, more generally, $deg(D)=0$ for any principal divisor $(f)=D$ on a compact Riemann surface since a non-constant meromorphic function $f:X\to\mathbb{CP}^1$ has as many zeros as poles.
\end{theorem}

The homotopy category consists of topological spaces, equipped with morphisms of homotopy equivalence classes of continuous maps. The topological spaces $X$ and $\mathbb{CP}^1$ are isomorphic in this category if and only if they are homotopy equivalent. In fact, more generally, if $X$ and $Y$ are two topological spaces that are homotopy equivalent (of the same homotopy type), then their homology groups are equal $H_n(X)=H_n(Y)$ for all $n\ge 0$ \cite[Pg. 13-18]{Marco}. Thus, the isomorphism $X\cong\mathbb{CP}^1$ induces a homotopy equivalence, which means that $H_n(X;M)=H_n(\mathbb{CP}^1;M)$ with a group coefficient $M.$ This leads to the following theorem:
\begin{theorem}\label{theorem:homotopy} 
A homotopy equivalence between the two topological spaces $X$ and $\mathbb{CP}^1$ is induced by the isomorphism $X\cong\mathbb{CP}^1$, such that they share homology groups $H_n(X;M)=H_n(\mathbb{CP}^1;M)$ for all $n\ge 0.$
\end{theorem}
The fact that if $X$ and $Y$ are two topological spaces that are homotopy equivalent (of the same homotopy type) then the homology groups are equal $H_n(X)=H_n(Y),$ for all $n\ge 0,$ was invoked to prove the above theorem.
\begin{theorem}
If $X$ and $Y$ are two topological spaces that are homotopy equivalent (of the same homotopy type), then the homology groups are equal $H_n(X)=H_n(Y)$ for all $n\ge 0.$
\end{theorem}
\proof
A continuous mapping $f:X\to Y$ induces a homomorphism $f_{\sharp}:C_n(X)\to C_n(Y).$ It follows that $f_{\sharp}$ is a chain map, such that $\partial f_{\sharp}=f_{\sharp}\partial,$ descending to homomorphisms on homology $f_{*}:H_n(X)\to H_n(Y).$ If $f$ and $g$ are homotopically equivalent then $f_{*}=g_{*},$ from which it follows that if $f$ is a homotopy equivalence (i.e. $X$ and $Y$ are homotopy equivalent) then $f_{*}$ must necessarily be an isomorphism. As such, let $F:X\times[0,1]\to Y$ be a homotopy map that takes $f$ to $g.$ We define a homomorphism on the level of chains, $P:C_n(X)\to C_{n+1}(Y),$ that takes a basis element $\sigma:\Delta^n\to X,$ a generator of $C_n(X),$ to the prism $P(\sigma):\Delta^n\times I\to Y.$ The boundary, obtained by the alternating formal sum, is $\partial P(\sigma)=f_{\sharp}(\sigma)-g_{\sharp}(\sigma)+P(\partial \sigma).$ Therefore, if $\alpha\in C_n(X)$ is an $n$-cycle then $f_{\sharp}(\alpha)$ and $g_{\sharp}(\alpha)$ only differ by the boundary$f_{\sharp}(\alpha)-g_{\sharp}(\alpha)=\partial P(\alpha),$ which means that the homomorphisms are homologous. The proof of the theorem is now complete. Theorem \ref{theorem:homotopy} follows at once.
\qed

\section{Homology for The Generalized Result and Concluding Remarks\label{sec:homology}}
We now consider the local exactness of the above condition for compact support. The characterization of connectivity of a region leads to the important idea of homology, in a complex-analytic sense. In particular, we give a contextualized definition \cite[Pg. 141]{Ahlfors}.
\begin{definition}
A cycle $\gamma=\partial X$ in a open set $\Omega$ is said to be homologous to zero with respect to $\Omega$ if $n(\gamma,a)=0$ for all points $a$ in the complement of $\Omega.$ To denote this relation, we write $\gamma\sim 0(mod \Omega).$ The notation $\gamma_1\sim\gamma_2$ is equivalent to $\gamma_1-\gamma_2\sim 0.$
\end{definition} 
The homology $\gamma\sim 0(mod \Omega)$ implies $\gamma\sim 0(mod \Omega')$ for $\Omega\subset\Omega'.$ In this particular case, if $\int_Xe^{-log\Gamma}\Delta\rho dS=0$ for the Riemann surface $X$ of genus $g=0,$ and if $X\cong\mathbb{CP}^1$ for $X$ given parametrically in terms of local coordinates $(z_i)$ by the chart $\phi:X\to\mathbb{C},$ $\phi(z_1,z_2)=r(z_1,z_2)tan\left(\frac{\pi}{4}+\frac{z_2}{2}\right)(cosz_1+isinz_2)$ for $r(z_1,z_2)>0$ and the diffeomorphism $\eta:X\to\mathbb{CP}^1$ with $\eta(z_1,z_2)=-\bar g(\eta(z_1,z_2))e^{i\pi\langle\chi,\xi\rangle},$ for $\omega\in\mathcal{M}\Omega^{(1,0)}(\mathbb{CP}^1),\text{ } \eta^*\omega\in\mathcal{M}\Omega^{(1,0)}(X)$ meromorphic $(1,0)$-forms, then $\partial X_j$ is homologous to zero with respect to $X_j,$ i.e. $\partial X_j\sim 0(mod X_j),$ and $\partial X_j\sim 0(mod X\setminus X_j)$ for $X\setminus X_j\supset X_j.$ Hence, by Cauchy's theorem, if $f(z)$ is analytic in $X_j,$ then $\int_{\partial X_j}f(z)dz=0$ for every cycle $\partial X_j$ which is homologous to zero in $X_j,$ whereby $\int_{\partial X_j}f(z)dz=\int_{\partial X_j}\eta^*\omega\equiv 0$ with $\sup_{\phi=(z_1,z_2)}|g|<+\infty.$ Thus, it follows that $\eta^*\omega$ is necessarily an exact differential form. Note that the homology groups, with coefficients in $\mathbb{C},$ of the topological space $\mathbb{CP}^1$ are given by
\begin{equation*}
H_k(\mathbb{CP}^1,\mathbb{C}):=\begin{cases} \mathbb{C} & \text{for } k=0,2, \\ 0 & \text{otherwise, } \end{cases}
\end{equation*} which follows from the existence of the diffeomorphism $\mathbb{S}^2\cong\mathbb{CP}^1$ via stereographic projection. This simple observation became important in proving $ker\left(\int_{c^1}\theta^1\right)=H_{dR}^1(X)$ for $X\cong\mathbb{CP}^1$ where $c^1$ denotes a $1$-cycle in $[c^1].$
The condition of compact support yields the following theorem [see Appendix \ref{sec:homology} for a complete complex-analytic proof].
\begin{theorem}\label{theorem:exactness}
The integral $\int_{\partial X_j}\eta^*\omega=\int_{\partial X_j}M(z_1,z_2)d\lambda+N(z_1,z_2)d\phi$ is locally exact in $X_j,$ which implies that $\int_{\partial X_j}\eta^*\omega=0$ for every cycle $\partial X_j\sim 0$ in $\partial X_j.$
\end{theorem}
\begin{remark}
We conclude by presenting a particularly revealing example for which the above theory applies, namely the case in which $g\in\Omega_{\mathbb{C}}^0(logD)\subset\mathcal{M}\Omega^{0}(\mathbb{C})$ is the logarithmic derivative of the Riemann Xi-function, i.e. $g=\xi'(z)/\xi(z),$ defined in the sense of Landau as $\xi(z)=\frac{1}{2}z(z-1)\pi^{-z/2}\Gamma(\frac{z}{2})\zeta(z).$ Note that here $D\subset\mathbb{C}$ is a divisor of $\mathbb{C}$ consisting of the set of points for which $g$ is singular; namely, the set of zeros of $\xi$ is a subset of this divisor. Since $\xi$ is entire, satisfying the symmetry $\xi(z)=\xi(1-z),$ then it can be alternatively defined by the Weierstrass product
\begin{equation}
\xi(z):=\xi(0)\prod_{\rho, |Im(\rho)|}\left(1-\frac{z}{\rho}\right)
\end{equation}
where the product extends over the non-trivial zeros of the Riemann zeta function $\zeta(z),$ $\rho,$ in order of increasing $|Im\text{ }\rho|$. Here $N_j=\frac{1}{2\pi i}\int_{\partial U_j}\omega=\frac{1}{2\pi i}\int_{\partial U_j}\frac{\xi'(z)}{\xi(z)}dz$ counts the number of zeros of $\xi$ by the argument principle, and thereby the number of zeros of $\zeta$ in the compact region $U_j\subset\mathbb{C}$ where $\xi$ is entire, i.e. it has no poles. Note that the logarithmic derivative of the Riemann Xi-function can be expressed as
\begin{equation}
\frac{\xi'(z)}{\xi(z)}=\sum_{\rho}\frac{1}{z-\rho}.
\end{equation} Thus, by imposing the condition that $\int_Xe^{-log\Gamma}\Delta\rho dS=0$ for the Riemann surface $X$ of genus $g=0,$ and if $X\cong\mathbb{CP}^1$ for $X$ given parametrically in terms of local coordinates $(z_i)$ by the diffeomorphism $\eta:X\to\mathbb{CP}^1$ with $\eta(z_1,z_2)=-\frac{\bar\xi'(\eta(z_1,z_2))}{\bar\xi(\eta(z_1,z_2))}e^{i\pi\langle\chi,\xi\rangle},$ then we can locate precisely the values of $z\in\mathbb{C}$ for which $N_j=\frac{1}{2\pi i}\int_{\partial U_j}\frac{\xi'(z)}{\xi(z)}dz=0.$
\end{remark}
\begin{remark}
Furthermore, the theory developed in the above analysis was applied to complex meromorphic differential $(1,0)$-forms $\omega\in\mathcal{M}\Omega^{(1,0)}.$ However, a similar theory can be developed for meromorphic $(p,q)$-forms in the space $\Omega^{(p,q)}=\bigwedge_p\Omega^{(1,0)}\bigwedge_q\Omega^{(0,1)},$ stable under a holomorphic change of coordinates. Thus, in local coordinates the $(p,q)$-form may be expressed as $\omega=\sum_{|I|=p,|J|=q}f_{IJ}dz^I\wedge d\bar z^J\in\Omega^{(p,q)}$ for $I, J$ multi-indices, where $\Omega^{(p,q)}$ is equipped with the Dolbeaut operators $\partial:\Omega^{(p,q)}\to\Omega^{(p+1,q)}$ and $\bar\partial:\Omega^{(p,q)}\to\Omega^{(p,q+1)}.$
\end{remark}

\appendix\section{Gaussian Curvature\label{sec:Gaussian Curvature}}

The Riemann surface $X$ can be parameterized in terms of the local coordinates $z_i,$ $i=1,2$ as 
\begin{equation*}
\bm{F}(z_1,z_2)=\begin{pmatrix}\phi(z_1)cos(\alpha(z_1)+\beta(z_2)+C\pi)\\ \phi(z_1)sin(\alpha(z_1)+\beta(z_2)+C\pi)\\z_1+iz_2\end{pmatrix}.
\end{equation*} The coefficients of the first fundamental form may be given as $E=g_{11}=\langle\bm{F}_{z_1},\bm{F}_{z_1}\rangle,$ $F=g_{12}=\langle\bm{F}_{z_1},\bm{F}_{z_2}\rangle,$ and $G=g_{22}=\langle\bm{F}_{z_2},\bm{F}_{z_2}\rangle$ where
\begin{equation*}
\bm{F}_{z_1}=\begin{pmatrix}\frac{d\phi}{dz_1}cos(\alpha(z_1)+\beta(z_2)+C\pi)-\phi(z_1)sin(\alpha(z_1)+\beta(z_2)+C\pi)\frac{d\alpha}{dz_1}\\ \frac{d\phi}{dz_1}sin(\alpha(z_1)+\beta(z_2)+C\pi)+\phi(z_1)cos(\alpha(z_1)+\beta(z_2)+C\pi)\frac{d\alpha}{dz_1}\\1\end{pmatrix},
\end{equation*}
\begin{equation*}
\bm{F}_{z_2}=\begin{pmatrix} -\phi(z_1)sin(\alpha(z_1)+\beta(z_2)+C\pi)\frac{d\beta}{dz_2} \\ \phi(z_1)cos(\alpha(z_1)+\beta(z_2)+C\pi)\frac{d\beta}{dz_2} \\ i  \end{pmatrix} \text{ such that } \end{equation*}  \allowdisplaybreaks  \begin{equation*}
\begin{split}
&E=g_{11}=\langle\bm{F}_{z_1},\bm{F}_{z_1}\rangle=1+\biggl(\frac{d\phi}{dz_1}cos(\alpha(z_1)+\beta(z_2)+C\pi)-\phi(z_1)sin(\alpha(z_1)+\beta(z_2)\\&+C\pi)\frac{d\alpha}{dz_1}\biggr)^2+\left(\frac{d\phi}{dz_1}sin(\alpha(z_1)+\beta(z_2)+C\pi)+\phi(z_1)cos(\alpha(z_1)+\beta(z_2)+C\pi)\frac{d\alpha}{dz_1}\right)^2\\&=1+\left(\frac{d\phi}{dz_1}\right)^2+\phi^2(z_1)\left(\frac{d\alpha}{dz_1}\right)^2,\\& 
F=g_{12}=\langle\bm{F}_{z_1},\bm{F}_{z_2}\rangle=i+\phi(z_1)cos(\alpha(z_1)+\beta(z_2)+C\pi)\frac{d\beta}{dz_2}\biggl(\frac{d\phi}{dz_1}sin(\alpha(z_1)+\beta(z_2)\\&+C\pi)+\phi(z_1)cos(\alpha(z_1)+\beta(z_2)+C\pi)\frac{d\alpha}{dz_1}\biggr)-\phi(z_1)sin(\alpha(z_1)+\beta(z_2)+C\pi)\frac{d\beta}{dz_2}\biggl(\frac{d\phi}{dz_1}\\&cos(\alpha(z_1)+\beta(z_2)+C\pi)-\phi(z_1)sin(\alpha(z_1)+\beta(z_2)+C\pi)\frac{d\alpha}{dz_1}\biggr)=i+\phi^2(z_1)\frac{d\beta}{dz_2}\frac{d\alpha}{dz_1},\text{ and }\\&
G=g_{22}=\langle\bm{F}_{z_2},\bm{F}_{z_2}\rangle=-1+\phi^2(z_1)sin^2(\alpha(z_1)+\beta(z_2)+C\pi)\left(\frac{d\beta}{dz_2}\right)^2\\&+\phi^2(z_1)cos^2(\alpha(z_1)+\beta(z_2)+C\pi)\left(\frac{d\beta}{dz_2}\right)^2=\phi^2(z_1)\left(\frac{d\beta}{dz_2}\right)^2-1.
\end{split}
\end{equation*} Therefore, after simplifying terms $E,G,F,$ the Riemannian metric can be given locally as 
\begin{equation*}
\begin{split}
&g=ds^2=\biggl(1+\left(\frac{d\phi}{dz_1}cos(\alpha(z_1)+\beta(z_2)+C\pi)-\phi(z_1)sin(\alpha(z_1)+\beta(z_2)+C\pi)\frac{d\alpha}{dz_1}\right)^2\\&+\left(\frac{d\phi}{dz_1}sin(\alpha(z_1)+\beta(z_2)+C\pi)+\phi(z_1)cos(\alpha(z_1)+\beta(z_2)+C\pi)\frac{d\alpha}{dz_1}\right)^2\biggr)dz_1^2\\&+2\biggl(i+\phi(z_1)cos(\alpha(z_1)+\beta(z_2)+C\pi)\frac{d\beta}{dz_2}\biggl(\frac{d\phi}{dz_1}sin(\alpha(z_1)+\beta(z_2)+C\pi)+\phi(z_1)cos(\alpha(z_1)\\&+\beta(z_2)+C\pi)\frac{d\alpha}{dz_1}\biggr)-\phi(z_1)sin(\alpha(z_1)+\beta(z_2)+C\pi)\frac{d\beta}{dz_2} \biggl(\frac{d\phi}{dz_1}cos(\alpha(z_1)+\beta(z_2)+C\pi)\\&-\phi(z_1)sin(\alpha(z_1)+\beta(z_2)+C\pi)\frac{d\alpha}{dz_1}\biggr)\biggr)dz_1dz_2+\biggl(\phi^2(z_1)sin^2(\alpha(z_1)+\beta(z_2)+C\pi)\left(\frac{d\beta}{dz_2}\right)^2\\&+\phi^2(z_1)cos^2(\alpha(z_1)+\beta(z_2)+C\pi)\left(\frac{d\beta}{dz_2}\right)^2-1\biggr)dz_2^2\\&
=\left(1+\left(\frac{d\phi}{dz_1}\right)^2+\phi^2(z_1)\left(\frac{d\alpha}{dz_1}\right)^2\right)dz_1^2+2\left(i+\phi^2(z_1)\frac{d\beta}{dz_2}\frac{d\alpha}{dz_1}\right)dz_1dz_2\\&+\left(\phi^2(z_1)\left(\frac{d\beta}{dz_2}\right)^2-1\right)dz_2^2.
\end{split}
\end{equation*}
Thus, the $EG-F^2$ term becomes
\begin{equation}
EG-F^2=\phi^2(z_1)\left[\left(\frac{d\beta}{dz_2}\right)^2+\left(\frac{d\phi}{dz_1}\right)^2\left(\frac{d\beta}{dz_2}\right)^2-\left(\frac{d\alpha}{dz_1}\right)^2\right]-\left(\frac{d\phi}{dz_1}\right)^2-2i\frac{d\beta}{dz_2}\frac{d\alpha}{dz_1}.
\end{equation}
Likewise, the Gaussian curvature is given by $K=\frac{eg-f^2}{EG-F^2}$ where $E=g_{11}=\langle\bm{F}_{z_1},\bm{F}_{z_1}\rangle,$ $F=g_{12}=\langle\bm{F}_{z_1},\bm{F}_{z_2}\rangle,$ and $G=g_{22}=\langle\bm{F}_{z_2},\bm{F}_{z_2}\rangle$ are coefficients of the first fundamental form and $e=\langle \bm{N}, \bm{F}_{z_1z_1} \rangle,$ $f=\langle \bm{N}, \bm{F}_{z_1z_2} \rangle,$ and $g=\langle \bm{N}, \bm{F}_{z_2z_2} \rangle$ are the coefficients of the second fundamental form, where $\bm{N}$ is the normal vector $\bm{N}=\bm{F}_{z_1}\times\bm{F}_{z_2}$. Calculation in the standard Euclidean basis gives the second order partial derivatives,
\begin{equation*}
\begin{split}
&\bm{F}_{z_1z_1}
=\biggl(\frac{d^2\phi}{dz_1^2}cos(\alpha(z_1)+\beta(z_2)+C\pi)-2\frac{d\phi}{dz_1}\frac{d\alpha}{dz_1}sin(\alpha(z_1)+\beta(z_2)+C\pi)-\phi(z_1)cos(\alpha(z_1)\\&+\beta(z_2)+C\pi)\left(\frac{d\alpha}{dz_1}\right)^2-\phi(z_1)sin(\alpha(z_1)+\beta(z_2)+C\pi)\frac{d^2\alpha}{dz_1^2},\frac{d^2\phi}{dz_1^2}sin(\alpha(z_1)+\beta(z_2)+C\pi)\\&+2\frac{d\phi}{dz_1}\frac{d\alpha}{dz_1}cos(\alpha(z_1)+\beta(z_2)+C\pi)-\phi(z_1)sin(\alpha(z_1)+\beta(z_2)+C\pi)\left(\frac{d\alpha}{dz_1}\right)^2\\&+\phi(z_1)cos(\alpha(z_1)+\beta(z_2)+C\pi)\frac{d^2\alpha}{dz_1^2},0\biggr),
\end{split}
\end{equation*}
\begin{equation*}
\bm{F}_{z_2z_2}=\begin{pmatrix} -\phi(z_1)cos(\alpha(z_1)+\beta(z_2)+C\pi)\left(\frac{d\beta}{dz_2}\right)^2-\phi(z_1)sin(\alpha(z_1)+\beta(z_2)+C\pi)\frac{d^2\beta}{dz_2^2} \\ \phi(z_1)cos(\alpha(z_1)+\beta(z_2)+C\pi)\frac{d^2\beta}{dz_2^2}-\phi(z_1)sin(\alpha(z_1)+\beta(z_2)+C\pi)\left(\frac{d\beta}{dz_2}\right)^2 \\ 0 \end{pmatrix},
\end{equation*} and 
\begin{equation*}
\bm{F}_{z_1z_2}=\begin{pmatrix} -\frac{d\phi}{dz_1}\frac{d\beta}{dz_2}sin(\alpha(z_1)+\beta(z_2)+C\pi)-\phi(z_1)\frac{d\alpha}{dz_1}\frac{d\beta}{dz_2}cos(\alpha(z_1)+\beta(z_2)+C\pi) \\ \frac{d\phi}{dz_1}\frac{d\beta}{dz_2}cos(\alpha(z_1)+\beta(z_2)+C\pi)-\phi(z_1)\frac{d\alpha}{dz_1}\frac{d\beta}{dz_2}sin(\alpha(z_1)+\beta(z_2)+C\pi)\\  0 \end{pmatrix}.
\end{equation*} Furthermore, the normal vector $\bm{N}$ to the Riemann surface $X$ is given by $\bm{N}=\frac{\partial \bm{F}}{\partial z_1}\times\frac{\partial\bm{F}}{\partial z_2}$
\begin{equation*}
\begin{split}
&=\biggl(i\frac{d\phi}{dz_1}sin(\alpha(z_1)+\beta(z_2)+C\pi)+i\phi(z_1)cos(\alpha(z_1)+\beta(z_2)+C\pi)\frac{d\alpha}{dz_1}-\phi(z_1)cos(\alpha(z_1)\\&+\beta(z_2)+C\pi)\frac{d\beta}{dz_2},i\phi(z_1)sin(\alpha(z_1)+\beta(z_2)+C\pi)\frac{d\alpha}{dz_1}-i\frac{d\phi}{dz_1}cos(\alpha(z_1)+\beta(z_2)+C\pi)\\&-\phi(z_1)sin(\alpha(z_1)+\beta(z_2)+C\pi)\frac{d\beta}{dz_2},\phi(z_1)\frac{d\phi}{dz_1}\frac{d\beta}{dz_2}\biggr).
\end{split}
\end{equation*} Lastly, we calculate the coefficients of the second fundamental form, which can in fact be expressed as $e=\langle \bm{N}, \bm{F}_{z_1z_1} \rangle,$ $f=\langle \bm{N}, \bm{F}_{z_1z_2} \rangle,$ and $g=\langle \bm{N}, \bm{F}_{z_2z_2} \rangle$ for
\begin{equation*}
\begin{split}
&e=\left\langle\bm{N},\frac{\partial^2\bm{F}}{\partial z_1^2}\right\rangle=\biggl(i\frac{d\phi}{dz_1}sin(\alpha(z_1)+\beta(z_2)+C\pi)+i\phi(z_1)cos(\alpha(z_1)+\beta(z_2)+C\pi)\frac{d\alpha}{dz_1}\\&-\phi(z_1)cos(\alpha(z_1)+\beta(z_2)+C\pi)\frac{d\beta}{dz_2}\biggr)\biggl(\frac{d^2\phi}{dz_1^2}cos(\alpha(z_1)+\beta(z_2)+C\pi)-2\frac{d\phi}{dz_1}\frac{d\alpha}{dz_1}sin(\alpha(z_1)\\&+\beta(z_2)+C\pi)-\phi(z_1)cos(\alpha(z_1)+\beta(z_2)+C\pi)\left(\frac{d\alpha}{dz_1}\right)^2-\phi(z_1)sin(\alpha(z_1)+\beta(z_2)+C\pi)\frac{d^2\alpha}{dz_1^2}\biggr)\\&+\biggl(i\phi(z_1)sin(\alpha(z_1)+\beta(z_2)+C\pi)\frac{d\alpha}{dz_1}-i\frac{d\phi}{dz_1}cos(\alpha(z_1)+\beta(z_2)+C\pi)-\phi(z_1)sin(\alpha(z_1)\\&+\beta(z_2)+C\pi)\frac{d\beta}{dz_2}\biggr)\biggl(\frac{d^2\phi}{dz_1^2}sin(\alpha(z_1)+\beta(z_2)+C\pi)+2\frac{d\phi}{dz_1}\frac{d\alpha}{dz_1}cos(\alpha(z_1)+\beta(z_2)+C\pi)\\&-\phi(z_1)sin(\alpha(z_1)+\beta(z_2)+C\pi)\left(\frac{d\alpha}{dz_1}\right)^2+\phi(z_1)cos(\alpha(z_1)+\beta(z_2)+C\pi)\frac{d^2\alpha}{dz_1^2}\biggr)\\&
=\phi^2(z_1)\left(\frac{d\beta}{dz_2}\right)\left(\frac{d\alpha}{dz_1}\right)^2-\phi(z_1)\left(\frac{d\beta}{dz_2}\right)\left(\frac{d^2\phi}{dz_1^2}\right)-i\phi^2(z_1)\left(\frac{d\alpha}{dz_1}\right)^3+i\phi(z_1)\left(\frac{d\alpha}{dz_1}\right)\left(\frac{d^2\phi}{dz_1^2}\right)\\&-i\phi(z_1)\left(\frac{d\phi}{dz_1}\right)\left(\frac{d^2\alpha}{dz_1^2}\right)-i\phi(z_1)\left(\frac{d\phi}{dz_1}\right)\left(\frac{d\alpha}{dz_1}\right)^2sin(2\alpha(z_1)+2\beta(z_1)+2\pi C)\\&-2i\left(\frac{d\phi}{dz_1}\right)^2\left(\frac{d\alpha}{dz_1}\right),
\end{split}
\end{equation*}
\begin{equation*}
\begin{split}
&g=\left\langle\bm{N}, \frac{\partial^2\bm{F}}{\partial z_2^2}\right\rangle=\biggl(i\phi(z_1)sin(\alpha(z_1)+\beta(z_2)+C\pi)\frac{d\alpha}{dz_1}-i\frac{d\phi}{dz_1}cos(\alpha(z_1)+\beta(z_2)+C\pi)\\&-\phi(z_1)sin(\alpha(z_1)+\beta(z_2)+C\pi)\frac{d\beta}{dz_2}\biggr)\biggl(\phi(z_1)cos(\alpha(z_1)+\beta(z_2)+C\pi)\frac{d^2\beta}{dz_2^2}-\phi(z_1)sin(\alpha(z_1)\\&+\beta(z_2)+C\pi)\left(\frac{d\beta}{dz_2}\right)^2\biggr)-\biggl(i\phi(z_1)sin(\alpha(z_1)+\beta(z_2)+C\pi)\frac{d\alpha}{dz_1}-i\frac{d\phi}{dz_1}cos(\alpha(z_1)+\beta(z_2)\\&+C\pi)-\phi(z_1)sin(\alpha(z_1)+\beta(z_2)+C\pi)\frac{d\beta}{dz_2}\biggr)\biggl(\phi(z_1)cos(\alpha(z_1)+\beta(z_2)+C\pi)\frac{d^2\beta}{dz_2^2}-\phi(z_1)\\&sin(\alpha(z_1)+\beta(z_2)+C\pi)\left(\frac{d\beta}{dz_2}\right)^2\biggr)=-i\phi(z_1)\left(\frac{d\phi}{dz_1}\right)\left(\frac{d^2\beta}{dz_2^2}\right)-i\phi^2(z_1)\left(\frac{d\alpha}{dz_1}\right)\left(\frac{d\beta}{dz_2}\right)^2\\&+\phi^2(z_1)\left(\frac{d\beta}{dz_2}\right)^3,
\end{split}
\end{equation*} and 
\begin{equation*}
\begin{split}
&f=\left\langle\bm{N}, \frac{\partial^2\bm{F}}{\partial z_1\partial z_2}\right\rangle=\biggl(i\phi(z_1)sin(\alpha(z_1)+\beta(z_2)+C\pi)\frac{d\alpha}{dz_1}-i\frac{d\phi}{dz_1}cos(\alpha(z_1)+\beta(z_2)+C\pi)\\&-\phi(z_1)sin(\alpha(z_1)+\beta(z_2)+C\pi)\frac{d\beta}{dz_2}\biggr)\biggl(\frac{d\phi}{dz_1}\frac{d\beta}{dz_2}cos(\alpha(z_1)+\beta(z_2)+C\pi)-\phi(z_1)\frac{d\alpha}{dz_1}\frac{d\beta}{dz_2}\\&sin(\alpha(z_1)+\beta(z_2)+C\pi)\biggr)-\biggl(i\frac{d\phi}{dz_1}sin(\alpha(z_1)+\beta(z_2)+C\pi)+i\phi(z_1)cos(\alpha(z_1)+\beta(z_2)+C\pi)\\&\frac{d\alpha}{dz_1}-\phi(z_1)cos(\alpha(z_1)+\beta(z_2)+C\pi)\frac{d\beta}{dz_2}\biggr)\biggl(\frac{d\phi}{dz_1}\frac{d\beta}{dz_2}cos(\alpha(z_1)+\beta(z_2)+C\pi)-\phi(z_1)\frac{d\alpha}{dz_1}\frac{d\beta}{dz_2}\\&sin(\alpha(z_1)+\beta(z_2)+C\pi)\biggr)=-i\left(\frac{d\phi}{dz_1}\right)^2\left(\frac{d\beta}{dz_2}\right)-i\phi^2(z_1)\left(\frac{d\alpha}{dz_1}\right)^2\left(\frac{d\beta}{dz_2}\right)\\&+\phi^2(z_1)\left(\frac{d\alpha}{dz_1}\right)\left(\frac{d\beta}{dz_2}\right)^2.
\end{split}
\end{equation*} Thus, the Gaussian curvature $K=\frac{eg-f^2}{EG-F^2}$ is given by the following explicit formula

\begin{equation}
\begin{split}
&K=\biggl(\biggl(\phi^2(z_1)\left(\frac{d\beta}{dz_2}\right)\left(\frac{d\alpha}{dz_1}\right)^2-\phi(z_1)\left(\frac{d\beta}{dz_2}\right)\left(\frac{d^2\phi}{dz_1^2}\right)-i\phi^2(z_1)\left(\frac{d\alpha}{dz_1}\right)^3\\&+i\phi(z_1)\left(\frac{d\alpha}{dz_1}\right)\left(\frac{d^2\phi}{dz_1^2}\right)-i\phi(z_1)\left(\frac{d\phi}{dz_1}\right)\left(\frac{d^2\alpha}{dz_1^2}\right)-i\phi(z_1)\left(\frac{d\phi}{dz_1}\right)\left(\frac{d\alpha}{dz_1}\right)^2sin(2\alpha(z_1)+2\beta(z_1)\\&+2\pi C)-2i\left(\frac{d\phi}{dz_1}\right)^2\left(\frac{d\alpha}{dz_1}\right)\biggr)\biggl(-i\phi(z_1)\left(\frac{d\phi}{dz_1}\right)\left(\frac{d^2\beta}{dz_2^2}\right)-i\phi^2(z_1)\left(\frac{d\alpha}{dz_1}\right)\left(\frac{d\beta}{dz_2}\right)^2\\&+\phi^2(z_1)\left(\frac{d\beta}{dz_2}\right)^3\biggr)-\biggl(-i\left(\frac{d\phi}{dz_1}\right)^2\left(\frac{d\beta}{dz_2}\right)-i\phi^2(z_1)\left(\frac{d\alpha}{dz_1}\right)^2\left(\frac{d\beta}{dz_2}\right)\\&+\phi^2(z_1)\left(\frac{d\alpha}{dz_1}\right)\left(\frac{d\beta}{dz_2}\right)^2\biggr)\biggr)\biggr/\biggl(\phi^2(z_1)\left[\left(\frac{d\beta}{dz_2}\right)^2+\left(\frac{d\phi}{dz_1}\right)^2\left(\frac{d\beta}{dz_2}\right)^2-\left(\frac{d\alpha}{dz_1}\right)^2\right]\\&-\left(\frac{d\phi}{dz_1}\right)^2-2i\frac{d\beta}{dz_2}\frac{d\alpha}{dz_1}\biggr).
\end{split}
\end{equation}
As in the above analysis, to make the calculation of Gaussian curvature less tedious, we invoke the isothermal coordinate argument. If either the sufficient condition $g^2=-{\bar g}^2$ for $g\in\mathcal{M}\Omega^{(0,0)}(\mathbb{CP}^1)$ or the coupled differential equation $\left(\frac{d\phi}{dz_1}\right)^2+\phi^2(z_1)\left[\left(\frac{d\alpha}{dz_1}\right)^2-\left(\frac{d\beta}{dz_2}\right)^2+2i\frac{d\beta}{dz_2}\frac{d\alpha}{dz_1}\right]=0$ is satisfied then $(z_1,z_2)$ must necessarily be local isothermal coordinates on $X,$ such that \begin{equation}
\chi(X)=\frac{1}{2\pi}\iint\limits_XKdS=-\frac{1}{4\pi}\iint\limits_X\frac{\Delta\rho}{e^{\rho}}dS=-\frac{1}{4\pi}\iint\limits_X\frac{\Delta log\Gamma}{\Gamma}dS=2
\end{equation} means that the Riemann surface has genus zero. Since $\rho=log\Gamma,$ for $\Gamma=1+\left(\frac{d\phi}{dz_1}\right)^2+\phi^2(z_1)\left(\frac{d\alpha}{dz_1}\right)^2,$ is solely dependent on the first local coordinate $z_1$ then 
\begin{equation*}
\begin{split}
&\Delta\rho=\Delta log\Gamma=\frac{\partial^2}{\partial z_1^2}log\Gamma=
2\frac{\partial}{\partial z_1}\left[\frac{\left(\frac{d\phi}{dz_1}\right)\left(\frac{d^2\phi}{dz_1^2}\right)+\phi(z_1)\left(\frac{d\phi}{dz_1}\right)\left(\frac{d\alpha}{dz_1}\right)^2+\phi^2(z_1)\left(\frac{d\alpha}{dz_1}\right)\left(\frac{d^2\alpha}{dz_1^2}\right)}{1+\left(\frac{d\phi}{dz_1}\right)^2+\phi^2(z_1)\left(\frac{d\alpha}{dz_1}\right)^2}\right]\\&
=\frac{2}{\Gamma}\biggl[\left(\frac{d^2\phi}{dz_1^2}\right)^2+\left(\frac{d\phi}{dz_1}\right)\left(\frac{d^3\phi}{dz_1^3}\right)+\left(\frac{d\phi}{dz_1}\right)^2\left(\frac{d\alpha}{dz_1}\right)^2+\phi(z_1)\left(\frac{d^2\phi}{dz_1^2}\right)\left(\frac{d\alpha}{dz_1}\right)^2+4\phi(z_1)\left(\frac{d\phi}{dz_1}\right)\\&\left(\frac{d\alpha}{dz_1}\right)\left(\frac{d^2\alpha}{dz_1^2}\right)+\phi^2(z_1)\left(\frac{d^2\alpha}{dz_1^2}\right)^2+\phi^2(z_1)\left(\frac{d\alpha}{dz_1}\right)\left(\frac{d^3\alpha}{dz_1^3}\right)\biggr]+\frac{4}{\Gamma^2}\biggl[\biggl(\left(\frac{d\phi}{dz_1}\right)\left(\frac{d^2\phi}{dz_1^2}\right)\\&+\phi(z_1)\left(\frac{d\phi}{dz_1}\right)\left(\frac{d\alpha}{dz_1}\right)^2+\phi^2(z_1)\left(\frac{d\alpha}{dz_1}\right)\left(\frac{d^2\alpha}{dz_1^2}\right)\biggr)\biggl(\left(\frac{d\phi}{dz_1}\right)\left(\frac{d^2\phi}{dz_1^2}\right)+\phi(z_1)\left(\frac{d\phi}{dz_1}\right)\left(\frac{d\alpha}{dz_1}\right)^2\\&+\phi^2(z_1)\left(\frac{d\alpha}{dz_1}\right)\left(\frac{d^2\alpha}{dz_1^2}\right)\biggr)\biggr].
\end{split}
\end{equation*}
To conclude, we impose the condition that $(z_1,z_2)=z_1+iz_2$ are isothermal coordinates provided $E-G+2iF=0.$ This is equivalent to the coupled set of differential equations, Eq. \ref{equation: coupled differential equations},
\begin{equation}
\left(\frac{d\phi}{dz_1}\right)^2+\phi^2(z_1)\left[\left(\frac{d\alpha}{dz_1}\right)^2-\left(\frac{d\beta}{dz_2}\right)^2+2i\frac{d\beta}{dz_2}\frac{d\alpha}{dz_1}\right]=0.
\end{equation} In a similar vein, we impose the condition of biholomorphicity for which the map $\eta:X\to\mathbb{CP}^1$ is biholomorphic if and only if $\phi(z_1)=\lambda e^{z_1\frac{d\beta}{dz_2}-i\alpha(z_1)}$ for $\lambda\in \mathbb{C}.$ To obtain a differential equation purely in terms of $\alpha$ and $\beta,$ we compute the first order ordinary derivative $\frac{d\phi}{dz_1}=\lambda\left(\frac{d\beta}{dz_2}-i\frac{d\alpha}{dz_1}\right)e^{z_1\frac{d\beta}{dz_2}-i\alpha(z_1)}$ such that under substitution the condition for isothermal coordinates is automatically satisfied, inducing a tautology.

 \appendix\section{Homology\label{sec:homology}}
The condition of compact support in the above analysis gave a statement on the exactness of the pullback $\eta^*\omega$ [Theorem \ref{theorem:exactness}].
\begin{theorem}
The integral $\int_{\partial X_j}\eta^*\omega=\int_{\partial X_j}M(z_1,z_2)dz_1+N(z_1,z_2)dz_2$ is locally exact in $X_j,$ which implies that $\int_{\partial X_j}\eta^*\omega=0$ for every cycle $\partial X_j\sim 0$ in $\partial X_j.$
\end{theorem}
\proof [Proof, Ahlfors]
We invoke an elementary complex-analytic proof of Ahlfors  \cite[pages 144--146]{Ahlfors}. To simplify notation, let $\gamma$ denote $\partial X_j$ and $\Omega$ denote $X_j.$ Then we construct $\sigma,$ a polygonal approximation of $\gamma$ with horizontal and vertical sides such that every locally exact differential form has the same integral over $\sigma$ and $\gamma.$ Using the property $n(\sigma,a)=n(\gamma,a)$ for $a\in\Omega',$ i.e. $\sigma\sim 0,$ it will be sufficient to prove the theorem for rectifiable polygonal curves, with sides parallel to the axes.
	
	Hence, for $\sigma$ an polygonal approximation, as described previously, of $\gamma,$ let the Euclidean distance from $\gamma$ to $\Omega'$ be $\rho.$ If $\gamma$ is given parametrically as $z=z(t),$ then the function $z(t)$ is uniformly continuous on the closed interval $[a,b].$ Furthermore, let $\delta>0$ such that $|z(t)-z(t')|<\rho$ for $|t-t'|<\delta$ by subdividing the interval $[a,b]$ into subintervals of length strictly less than $\rho.$ Then, the subarcs $\gamma_i$ of $\gamma$ have the property that each is contained in a disk of radius $\rho$ which is contained entirely in $\Omega.$ The end points of $\gamma_i$ can be joined by a polygon $\sigma_i$ within that disk, consisting of a horizontal and vertical line segment. The exactness of the differential form, that is the path-independence, in the disk implies that 
\begin{equation*}
\int_{\sigma_i}Mdz_1+Ndz_2=\int_{\gamma_i}Mdz_1+Ndz_2\equiv\int_{\partial X_{j,i}=\left(\partial X_j\right)_i}Mdz_1+Ndz_2,
\end{equation*} and for $\sigma:=\sum_i\sigma_i,$ we have
\begin{equation*}
\int_{\sigma}Mdz_1+Ndz_2=\int_{\gamma}Mdz_1+Ndz_2=\int_{\partial X_j}Mdz_1+Ndz_2.
\end{equation*} Proceeding, we invoke the following construction: extend all segments that make up $\sigma$ to infinite lines. They must divide the plane into some finite rectangles $R_i$ and some unbounded regions $R_j',$ thought of as infinite rectangles.

Choosing a point $a_i$ from the interior of each $R_i,$ we form the cycle 
\begin{equation*}
\sigma_0=\sum_in(\sigma,a_i)\partial R_i
\end{equation*} where the sum extends over all finite rectangles. The coefficients $n(\sigma,a_i)$ are well determined for no $a_i$ lies on $\sigma.$ Similarly, let $a_j'$ denote the points chosen from the interior or each $R_j'.$ For $k=i,$ the index is $n(\partial R_i,a_k)=1,$ whereas $n(\partial R_i,a_k)=0$ if $k\ne i,$ and likewise, $n(\partial R_i,a_j')=0$ for all $j.$ From $\sigma_0=\sum_in(\sigma,a_i)\partial R_i,$ it becomes clear that $n(\sigma_0,a_i)=n(\sigma,a_i)$ and $n(\sigma_0,a_j')=0,$ with $R_j'$ an unbounded region. From linearity, it follows that $n(\sigma-\sigma_0,a)=0$ for $a=a_i$ and $a=a_j'.$
Let $\sigma_{ik}$ be the side shared by adjacent rectangles $R_i,R_k,$ where the orientation is such that $R_i$ lies to the left of $\sigma_{ik}.$ Suppose that the expression of $\sigma-\sigma_0$ contains the multiple $c\sigma_{ik};$ then, the cycle $\sigma-\sigma_0-c\partial R_i$ does not contain $\sigma_{ik},$ meaning that $a_i$ and $a_k$ must have the same index with respect to this cycle. The indices are $-c$ and $0,$ respectively, which forces $c=0.$ Similar reasoning applies to an infinite rectangle $R_j'.$ The common side occurs with coefficient zero in $\sigma-\sigma_0,$ which proves
\begin{equation*}
\sigma=\sum_in(\sigma,a_i)\partial R_i
\end{equation*}
by invoking $n(\sigma-\sigma_0,a)=0,$ meaning that $\sigma$ and $\sigma_0$ are equivalent up to cancellation of mutually shared boundaries. 

Lastly, to complete the proof, it must be shown that if $n(\sigma,a_i)\ne 0$ for $a_i\in R_i$ then $R_i$ is contained in $\Omega.$ Suppose that a point $a$ in the closed rectangle $R_i$ were not in $\Omega.$ Then $n(\sigma,a)=0$ for $\sigma\sim 0(mod \Omega).$ The curve joining $a$ and $a_i$ does not intersect $\sigma,$ which implies that $n(\sigma,a_i)=n(\sigma,a)\equiv 0.$ Therefore, by the local exactness of the integral, $\int_{\partial R_j}Mdz_1+Ndz_2$ over $\partial R_j$ is zero by $\sigma=\sum_in(\sigma,a_i)\partial R_i.$ Therefore, $\int_{\sigma}Mdz_1+Ndz_2=\int_{\partial X_j}Mdz_1+Ndz_2=0,$ as was to be shown.
\qed

\appendix\section{Dolbeault's Lemma\label{sec:Dolbeault's Lemma}}
We prove the following lemma that guarantees the existence of a solution to the inhomogeneous Cauchy-Riemann differential equation $\partial f/\partial\bar{z}=g.$
\begin{lemma} 
Suppose $g\in\mathcal{E}(\mathbb{C})$ has compact support. Then there exists a function $f\in\mathcal{E}(\mathbb{C})$ such that $\frac{\partial f}{\partial \bar z}=g.$
\end{lemma}
\proof
Define the function $f:\mathbb{C}\to\mathbb{C}$ by 
\begin{equation*}
f(\zeta)=\frac{1}{2\pi i}\int\int_{\mathbb{C}}\frac{g(z)}{z-\zeta}dz\wedge d\bar z.
\end{equation*} Since the integrand has a singular point when $z=\zeta,$ one has to show that the integral exists and depends differentiably on $\zeta.$ The simplest way is to change variables by translation and then introduce polar coordinates $r,\theta.$ Namely, let $z=\zeta+re^{i\theta}.$ Performing the integration, one treats $\zeta$ as a constant, where the polar transformation yields $dz\wedge d\bar z=-2i dx\wedge dy=-2irdr\wedge d\theta.$ Consequently 
\begin{equation*}
\begin{split}
&f(\zeta)=-\frac{1}{\pi}\int\int\frac{g(\zeta+re^{i\theta})}{re^{i\theta}}rdrd\theta\\&=-\frac{1}{\pi}\int\int g(\zeta+re^{i\theta})e^{-i\theta}drd\theta.
\end{split}
\end{equation*} By hypothesis $g\in\mathcal{E}(\mathbb{C})$ has compact support, meaning that one has to only integrate over the rectangular region $0\le r\le R,0\le\theta\le 2\pi,$ given that $R$ is sufficiently large. One may then differentiate under the integral sign for $f\in\mathcal{E}(\mathbb{C}),$ such that 
\begin{equation*}
\frac{\partial f}{\partial \bar\zeta}(\zeta)=-\frac{1}{\pi}\int\int\frac{\partial g\left(\zeta+re^{i\theta}\right)}{\partial \bar\zeta}e^{-i\theta}drd\theta.
\end{equation*} Transforming back to the original coordinates and letting $B_{\epsilon}:=\{z\in\mathbb{C}:\epsilon\le|z|\le R\}$, 
\begin{equation*}
\frac{\partial f}{\partial \bar\zeta}(\zeta)=\frac{1}{2\pi i}\lim_{\epsilon\to 0}\iint\limits_{B_{\epsilon}}\frac{\partial g(\zeta+z)}{\partial\bar\zeta}\frac{1}{z}dz\wedge d\bar z.
\end{equation*} Assuming that $z\ne 0,$ then $\frac{\partial g(\zeta+z)}{\partial\bar\zeta}\frac{1}{z}=\frac{\partial g(\zeta+z)}{\partial\bar z}\frac{1}{z}=\frac{\partial}{\partial\bar z}\left(\frac{g(\zeta+z)}{z}\right)$ and thus, 
\begin{equation*}
\frac{\partial f}{\partial \bar\zeta}(\zeta)=\frac{1}{2\pi i}\lim_{\epsilon\to 0}\iint\limits_{B_{\epsilon}}\frac{\partial}{\partial\bar z}\left(\frac{g(\zeta+z)}{z}\right)dz\wedge d\bar z=-\lim_{\epsilon\to 0}\iint\limits_{B_{\epsilon}}d\omega,
\end{equation*} where the differential $1$-form $\omega$ is given by $\omega(z)=\frac{1}{2\pi i}\frac{g(\zeta+z)}{z}dz$ for $z$ a variable and $\zeta$ a constant. Thus, by Stokes' Theorem,
\begin{equation*}
\frac{\partial f}{\partial \bar\zeta}(\zeta)=-\lim_{\epsilon\to 0}\iint\limits_{B_{\epsilon}}d\omega=-\lim_{\epsilon\to 0}\int_{\partial B_{\epsilon}}\omega=\lim_{\epsilon\to 0}\int_{|z|=\epsilon}\omega.
\end{equation*} By parameterizing the circle $|z|=\epsilon$ by $z=\epsilon e^{i\theta},$ $0\le\theta \le 2\pi,$ one obtains $\frac{\partial f}{\partial \bar\zeta}(\zeta)=\lim_{\epsilon\to 0}\frac{1}{2\pi}\int\limits_0^{2\pi}g(\zeta+\epsilon e^{i\theta})d\theta,$ which is the average value of the continuous function $g$ over the circle $\zeta+\epsilon e^{i\theta}$ for $0\le\theta \le 2\pi.$ By continuity, the above integral converges to $g(\zeta)$ as $\epsilon\to 0$ such that 
\begin{equation*}
\frac{\partial f}{\partial \bar\zeta}(\zeta)=g(\zeta).
\end{equation*} The proof of the lemma is now complete.
 \qed

\newpage
\bibliographystyle{plain}

\begin{thebibliography}{9}
\bibitem{Ahlfors}
Ahlfors, L. V. (2007). Complex Analysis: An Introduction to The Theory of Analytic Functions of One Complex Variable. New York: McGraw-Hill.
\bibitem{Ahlfors1}
Ahlfors, Lars V. (2006) [1966], Lectures on Quasiconformal Mappings, University Lecture Series, 38 (2nd ed.), Providence, R.I
\bibitem{Darling}
 Darling, R. W. R. (1994). Differential Forms and Connections. Cambridge University Press.
\bibitem{Forster}
Forster, O., and Gilligan, B. (1991). Lectures on Riemann surfaces. New York: Springer.
\bibitem{Gelfand}
Gelfand, I.M.; Minlos, R.A.; Shapiro, Z.Ya. (1963), Representations of the Rotation and Lorentz Groups and their Applications, New York: Pergamon Press
\bibitem{Griffiths}
Griffiths, Phillip; Harris, Joseph (1994), Principles of Algebraic Geometry, Wiley Classics Library, New York: John Wiley and Sons
\bibitem{Ionel}
Ionel, E. (n.d.). Complex Manifolds. Math Stanford, 1-75. Retrieved April 15, 2017.
\bibitem{Lee}
Lee, J. M. (2006). Introduction to Smooth Manifolds. New York: Springer.
\bibitem{Marco}
Marco, A., and Pérez, B. (2010). De Rham's Theorem. UQAM.
\bibitem{Peters}
Peters, C., and Steenbrink, J. H. (2008). Mixed Hodge structures. Berlin: Springer.
\bibitem{Spivak}
Spivak, Michael (1979), A Comprehensive Introduction to Differential Geometry. Vol. IV (2nd ed.), Publish or Perish
\bibitem{van der Waerden}
van der Waerden, B. L. (1952), Group Theory and Quantum Mechanics, Springer Publishing, ISBN 978-3642658624 (translation of the original 1932 edition, Die Gruppentheoretische Methode in Der Quantenmechanik).
\bibitem{Whitney}
Whitney, H. (1957), Geometric Integration Theory, Princeton Mathematical Series, 21, Princeton, NJ and London: Princeton University Press and Oxford University Press, pp. XV+387, MR 0087148, Zbl 0083.28204
\end{thebibliography}


\end{document}